\documentclass[11pt]{amsart} 
\usepackage[english]{babel}
\usepackage[latin1]{inputenc}
\usepackage{color}
\usepackage[colorinlistoftodos]{todonotes}
\usepackage[normalem]{ulem}
\usepackage{amssymb}
\usepackage{amsmath}
\usepackage{amscd}
\usepackage{latexsym}
\usepackage{amsthm}
\usepackage{mathrsfs}
\usepackage{mathabx}
\usepackage{tikz-cd}

\theoremstyle{plain}
\newtheorem{thm}{Theorem}[section]
\newtheorem{cor}[thm]{Corollary}
\newtheorem{lem}[thm]{Lemma}
\newtheorem{prop}[thm]{Proposition}

\newtheorem{defn}[thm]{Definition}
\newtheorem{oss}[thm]{Remark}

\newtheorem{sett}[thm]{Setting}

\title[The second integral cohomology of moduli spaces]{The second integral cohomology of moduli spaces of sheaves on K3 and Abelian surfaces}
\author{Arvid Perego}
\email{perego@dima.unige.it}
\address{Dipartimento di Matematica dell'Universit\`a di Genova, via Dodecaneso, 16148 Genova, Italy}
\author{Antonio Rapagnetta}
\email{rapagnet@axp.mat.uniroma2.it}
\address{Dipartimento di Matematica dell'Universit\`a di Roma II - Tor Vergata, 00133 Roma, Italy}

\begin{document}

\begin{abstract}
In this paper we study the second integral cohomology of moduli spaces of semistable sheaves on projective K3 surfaces. If $S$ is a projective K3 surface, $v$ a Mukai vector and $H$ a $v-$generic polarization on $S$, we show that $H^{2}(M_{v},\mathbb{Z})$ is a free $\mathbb{Z}-$module of rank 23 carrying a pure weight-two Hodge structure and a lattice structure, with respect to which $H^{2}(M_{v},\mathbb{Z})$ is Hodge isometric to the Hodge sublattice $v^{\perp}$ of the Mukai lattice of $S$. Similar results are proved for Abelian surfaces.
\end{abstract}

\maketitle

\tableofcontents

\section{Introduction and main results}

A compact, connected K\"ahler manifold $X$ is an \textit{irreducible symplectic manifold} if it is holomorphically symplectic, simply connected and $\dim(H^{0}(X,\Omega_{X}^{2}))=1$. Here, by holomorphically symplectic manifold we mean a complex manifold $X$ carrying an everywhere nondegenerate, closed, holomorphic $2-$form (called \textit{holomorphic symplectic form}). A complex manifold admitting a holomorphic symplectic form is called \textit{holomorphic symplectic manifold}. In particular, a holomorphic symplectic manifold (and hence an irreducible symplectic manifold) has even complex dimension and trivial canonical bundle. 

Irreducible symplectic manifolds are one of the three building blocks for compact K\"ahler manifolds with numerically trivial canonical bundle: the Bogomolov Decomposition Theorem asserts that if $X$ is a connected, compact K\"ahler manifold with numerically trivial canonical bundle, then there is a finite \'etale cover $Y$ of $X$ such that $$Y=T\times\prod_{i=1}^{n}X_{i}\times\prod_{j=1}^{m}Y_{j},$$where $T$ is a complex torus, the $Y_{j}$'s are irreducible symplectic manifolds, the $X_{i}$'s are irreducible Calabi-Yau manifolds (i. e. compact, connected, simply connected manifolds with trivial canonical bundle and such that $H^{0}(X_{i},\Omega_{X_{i}}^{p})=0$ for all $0<p<dim(X_{i})$).

If $X$ is an irreducible symplectic manifold, then $H^{2}(X,\mathbb{Z})$ is a free $\mathbb{Z}-$module carrying a pure weight-two Hodge structure and a nondegenerate integral quadratic form $b_{X}$ of signature $(3,b_{2}(X)-3)$, the \textit{Beauville form} of $X$ (see \cite{B}, Th\'eor\`eme 5 (a)).

The \textit{Beauville lattice} $(H^{2}(X,\mathbb{Z}),b_{X})$ of $X$ is an important deformation invariant of $X$. First, as shown in \cite{F} there is a positive rational number $C_{X}$, called \textit{Fujiki constant} of $X$, such that for every $\alpha\in H^{2}(X,\mathbb{Z})$ we have $$\int_{X}\alpha^{2n}=C_{X}b_{X}(\alpha)^{n},$$where $2n$ is the complex dimension of $X$ (this implies that the Beauville form and the Fujiki constant are deformation invariant).

Moreover, for irreducible symplectic manifold a Local and a Global Torelli Theorem hold. More precisely, if we let $Def(X)$ be the base of a Kuranishi family of $X$ and $$\Omega_{X}:=\{\alpha\in\mathbb{P}(H^{2}(X,\mathbb{C}))\,|\,b_{X}(\alpha)=0,\,\,b_{X}(\alpha+\overline{\alpha})>0\},$$there is a holomorphic map $p:Def(X)\longrightarrow\Omega_{X}$, called the \textit{period map}, which is a local biholomorphism (see \cite{B}, Th\'eor\`eme 5 (b)). Moreover, two irreducible symplectic manifolds $X$ and $Y$ are bimeromorphic if and only if there is a Hodge isometry between $H^{2}(X,\mathbb{Z})$ and $H^{2}(Y,\mathbb{Z})$ which comes from a parallel transport operator (see \cite{V}, \cite{Mark}, \cite{Huy}).

There are very few known deformation classes of irreducible symplectic manifolds, and for all of them the Beauville lattice is known.
\begin{enumerate}
 \item A compact, connected, smooth complex surface is an irreducible symplectic manifold if and only if it is a K3 surface. The Beauville form of a K3 surface is just the intersection pairing, and the Beauville lattice is isometric to the \textit{K3 lattice} $E_{8}^{\oplus 2}\oplus U(-1)^{\oplus 3}$.
 \item If $S$ is a K3 surface and $n\in\mathbb{N}$, $n\geq 2$, the Hilbert scheme $Hilb^{n}(S)$ of $n$ points on $S$ is an irreducible symplectic manifold of dimension $2n$ and second Betti number $23$ (see \cite{B}, Th\'eor\`eme 3 and Proposition 6). There is an isometry $$H^{2}(Hilb^{n}(S),\mathbb{Z})\simeq H^{2}(S,\mathbb{Z})\oplus^{\perp}\mathbb{Z}\cdot\delta,$$where $2\delta$ is the class of the exceptional divisor of the Hilbert-Chow morphism $\rho:Hilb^{n}(S)\longrightarrow Sym^{n}(S)$, whose square with respect to the Beauville form is $2-2n$.
 \item If $T$ is a $2-$dimensional complex torus and $n\in\mathbb{N}$, $n\geq 2$, then the fibers $Kum^{n}(T)$ of the sum morphism $Hilb^{n+1}(T)\longrightarrow T$ are irreducible symplectic manifolds of dimension $2n$ and second Betti number 7 (see \cite{B}, Th\'eor\`eme 4 and Proposition 8). There is an isometry $$H^{2}(Kum^{n}(T),\mathbb{Z})\simeq H^{2}(T,\mathbb{Z})\oplus^{\perp}\mathbb{Z}\cdot\delta,$$where $2\delta$ is the class of the restriction to $Kum^{n}(T)$ of the exceptional divisor of the Hilbert-Chow morphism, and the square of $\delta$ with respect to the Beauville form is $2-2n$.
 \item There are two more known deformation classes: $OG_{6}$, in dimension 6 and with second Betti number 8, and $OG_{10}$, in dimension 10 and with second Betti number 24 (see \cite{OG2}, \cite{OG3} and \cite{R}). Their Beauville lattices are computed in \cite{R} (for $OG_{6}$) and \cite{R1} (for $OG_{10}$).
\end{enumerate}

The singular analogue of irreducible symplectic manifolds is given by irreducible symplectic varieties, whose definition was introduced in \cite{GKP}. To recall them we need the following notation: if $X$ is a normal complex algebraic variety and $X_{reg}$ is the smooth locus of $X$ whose open embedding in $X$ is $j:X_{reg}\longrightarrow X$, for every $p\in\mathbb{N}$ such that $0\leq p\leq\dim(X)$ we let $$\Omega_{X}^{[p]}:=j_{*}\Omega^{p}_{X_{reg}}=\big(\wedge^{p}\Omega_{X}\big)^{**},$$whose global sections are called \textit{reflexive $p-$forms} on $X$. A reflexive $p-$form on $X$ is then a holomorphic $p-$form on $X_{reg}$. 

If $f:Y\longrightarrow X$ is a finite, dominant morphism between two irreducible normal varieties, then there is a morphism $f^{*}\Omega_{X}^{[p]}\longrightarrow\Omega_{Y}^{[p]}$ induced by the usual pull-back morphism of forms on the smooth loci, giving a morphism $f^{[*]}:H^{0}(X,\Omega_{X}^{[p]})\longrightarrow H^{0}(Y,\Omega_{Y}^{[p]})$, called \textit{reflexive pull-back morphism}.

We first recall from \cite{B2} the definitions of symplectic form and symplectic variety we will use to define the singular analogue of irreducible symplectic manifolds.

\begin{defn}
{\rm Let $X$ be a normal complex algebraic variety.
\begin{enumerate}
 \item A \textit{symplectic form} on $X$ is a closed reflexive $2-$form $\sigma$ on $X$ which is non-degenerate at each point of $X_{reg}$.
 \item If $\sigma$ is a symplectic form on $X$, the pair $(X,\sigma)$ is a \textit{symplectic variety} if for every resolution $f:\widetilde{X}\longrightarrow X$ of the singularities of $X$, the holomorphic symplectic form $\sigma_{reg}:=\sigma_{|X_{reg}}$ extends to a holomorphic $2-$form on $\widetilde{X}$.
 \item If $(X,\sigma)$ is a symplectic variety and $f:\widetilde{X}\longrightarrow X$ is a resolution of the singularities over which $\sigma_{reg}$ extends to a holomorphic symplectic form on $\widetilde{X}$, we say that $f$ is a \textit{symplectic resolution}.
\end{enumerate}}
\end{defn}

A symplectic variety has trivial canonical bundle and canonical (and hence rational) singularities. Conversely, by Theorem 6 of \cite{N3} a normal variety having rational Gorenstein singularities and whose regular locus carries a holomorphic symplectic form is a symplectic variety. 

Moreover, a normal variety having a symplectic form and whose singular locus has codimension at least 4 is a symplectic variety (see \cite{Fle}), and a symplectic variety has terminal singularities if and only if its singular locus has codimension at least 4 (Corollary 1 of \cite{N1}).

We now define irreducible Calabi-Yau and irreducible symplectic varieties following \cite{GKP}. If $X$ and $Y$ are two irreducible normal projective varieties, a \textit{finite quasi-\'etale morphism} $f:Y\longrightarrow X$ is a finite morphism which is \'etale in codimension one.

\begin{defn}
\label{defn:irrvar}
{\rm Let $X$ be an irreducible normal projective variety of dimension $d\geq 2$ with trivial canonical divisor and canonical singularities.
\begin{enumerate}
 \item The variety $X$ is \textit{irreducible Calabi-Yau} if for every $0<p<d$ and for every finite quasi-\'etale morphism $Y\longrightarrow X$, we have $H^{0}(Y,\Omega_{Y}^{[p]})=0$.
 \item The variety $X$ is \textit{irreducible symplectic} if it has a symplectic form $\sigma$, and for every finite quasi-\'etale morphism $f:Y\longrightarrow X$ the exterior algebra of reflexive forms on $Y$ is spanned by $f^{[*]}\sigma$.
\end{enumerate}}
\end{defn}

The definition of irreducible symplectic variety is motivated by the description of the algebra of holomorphic forms of an irreducible symplectic manifold, which is spanned by a holomorphic symplectic form.

By Proposition A.1 of \cite{HNW}, a smooth irreducible symplectic variety is an irreducible symplectic manifold. Moreover, by Corollary 13.3 of \cite{GGK} an irreducible symplectic variety $X$ is simply connected. In particular, the $\mathbb{Z}-$module $H^{2}(X,\mathbb{Z})$ is free.

Definition \ref{defn:irrvar} proves to be the good one in view of the Bogomolov Decomposition Theorem in the singular projective setting. H\"oring and Peternell (see Theorem 1.5 of \cite{HP}) show that if $X$ is an irreducible normal projective variety with klt singularities and has numerically trivial canonical bundle, then it admits a finite quasi-\'etale cover $Y$ such that $$Y=T\times\prod_{i=1}^{n}X_{i}\times\prod_{j=1}^{m}Y_{j},$$where $T$ is a complex torus, the $X_{i}$'s are irreducible Calabi-Yau varieties and the $Y_{j}$'s are irreducible symplectic varieties (see Proposition 5.20 of \cite{KM2}). 

Moreover, by Proposition 1.10 of \cite{PR3} if $X$ is an irreducible symplectic variety, then $X$ is a \textit{Namikawa symplectic variety}, i. e. a normal projective variety such that $h^{1}(X,\mathcal{O}_{X})=0$ and $h^{0}(X,\Omega_{X}^{[2]})=1$ (the converse is false, see Examples 1.12 and 1.13 of \cite{PR3}). Namikawa symplectic varieties, and hence irreducible symplectic varieties, share many features with irreducible symplectic manifolds. 

A first example of this is given by a symplectic variety $X$ admitting a symplectic resolution $f:\widetilde{X}\longrightarrow X$ which is an irreducible symplectic manifold: by Proposition 1.9 of \cite{PR3} the variety $X$ is a Namikawa symplectic variety (even if not necessarily an irreducible symplectic variety, see again Examples 1.12 and 1.13 of \cite{PR3}). 

Under these assumptions $X$ has rational singularities, and as a consequence $f^{*}:H^{2}(X,\mathbb{Z})\longrightarrow H^{2}(\widetilde{X},\mathbb{Z})$ is an inclusion of mixed Hodge structures: it follows that the $\mathbb{Z}-$module $H^{2}(X,\mathbb{Z})$ is free, has a pure weight-two Hodge structure and a nondegenerate integral quadratic form of signature $(3,b_{2}(X)-3)$. A Local and a Global Torelli Theorem for these varieties are proved in \cite{BL}.

A second example, more interesting for our purposes, is given by a $\mathbb{Q}-$factorial Namikawa symplectic variety $X$ having terminal singularities. The variety $X$ has no symplectic resolutions, and the following results hold.
\begin{enumerate}
 \item The (free part of the) $\mathbb{Z}-$module $H^{2}(X,\mathbb{Z})$ has a pure weight-two Hodge structure and a nondegenerate quadratic form $b_{X}$ of signature $(3,b_{2}(X)-3)$, called the \textit{Beauville form} of $X$ (see \cite{N3} and \cite{N4}).
 \item The deformations of $X$ are unobstructed and locally trivial (see \cite{N1} and \cite{N2}).
 \item There is a positive rational number $C_{X}$, called \textit{Fujiki constant} of $X$, verifying the same equality as the Fujiki constant of an irreducible symplectic manifold (see \cite{S}). As a consequence, the Beauville form and the Fujiki constant are deformation invariant.
 \item A Local and a Global Torelli Theorem hold (see \cite{N3}, Theorem 8, and \cite{BL2}).
\end{enumerate}

The case of Namikawa symplectic varieties having terminal singularities but which are not $\mathbb{Q}-$factorial has been studied in \cite{Kir}. 

Examples of irreducible symplectic varieties are obtained by several authors (see \cite{Pe} for an overview). Among them we cite the partial resolution of the quotient of a Hilbert scheme of two points on a K3 surface (resp. of a generalized Kummer variety of dimension 4) by the action of a symplectic involution, whose second integral cohomology is described in \cite{Men1} (resp. in \cite{KaMe}), and the quotients of Hilbert schemes of two points by the action of a symplectic automorphism of order 3, 5, 7 or 11, whose second integral cohomology is studied in \cite{Men2} and \cite{Men3}.

A family of examples coming from moduli spaces of semistable sheaves over K3 surfaces or Abelian surfaces is described in \cite{PR3}: the aim of this work is to describe their second integral cohomology.

\subsection{Notation and main results of the paper} 

In what follows $S$ will be a projective K3 surface or an Abelian surface, and we let $$\epsilon(S):=\left\{\begin{array}{ll} 1, & S\,\,{\rm is}\,\,{\rm K3}\\ 0, & S\,\,{\rm is}\,\,{\rm Abelian}\end{array}\right.$$We let $\rho(S)$ be the rank of the N\'eron-Severi group $NS(S)$ of $S$.

We let $\widetilde{H}(S,\mathbb{Z}):=H^{2*}(S,\mathbb{Z})$, and recall that $\widetilde{H}(S,\mathbb{Z})$ has a pure weight-two Hodge structure and a lattice structure with respect to the Mukai pairing $(.,.)$ (see \cite{HL}, Definitions 6.1.5 and 6.1.11). We let $v^{2}:=(v,v)$ for every $v\in\widetilde{H}(S,\mathbb{Z})$, and we call $\widetilde{H}(S,\mathbb{Z})$ the \textit{Mukai lattice} of $S$.

An element $v\in\widetilde{H}(S,\mathbb{Z})$ will be written $v=(v_{0},v_{1},v_{2})$, where $v_{i}\in H^{2i}(S,\mathbb{Z})$, and $v_{0},v_{2}\in\mathbb{Z}$. It will be called \textit{Mukai vector} if $v_{0}\geq 0$, $v_{1}\in NS(S)$ and if $v_{0}=0$, then either $v_{1}$ is the first Chern class of an effective divisor, or $v_{1}=0$ and $v_{2}>0$. 

If $\mathscr{F}$ is a coherent sheaf on $S$, we define its \textit{Mukai vector} as $$v(\mathscr{F}):=ch(\mathscr{F})\sqrt{td(S)}=(rk(\mathscr{F}),c_{1}(\mathscr{F}),ch_{2}(\mathscr{F})+\epsilon(S)rk(\mathscr{F})).$$

Let now $v$ be a Mukai vector on $S$ and suppose that $H$ is a $v-$generic polarization (see section 2.1 of \cite{PR3} for the definition). We write $M_{v}(S,H)$ (resp. $M_{v}^{s}(S,H)$) for the moduli space of Gieseker $H-$semistable (resp. $H-$stable) sheaves on $S$ with Mukai vector $v$. If $S$ is Abelian and $v^{2}>0$, we have a dominant isotrivial fibration $a_{v}:M_{v}(S,H)\longrightarrow S\times\widehat{S}$ (see section 4.1 of \cite{Y2}), where $\widehat{S}$ is the dual of $S$. We let $K_{v}(S,H):=a_{v}^{-1}(0_{S},\mathscr{O}_{S})$, and $K_{v}^{s}(S,H):=K_{v}(S,H)\cap M^{s}_{v}(S,H)$. If no confusion on $S$ and $H$ is possible, we drop them from the notation.

If $M_{v}^{s}\neq\emptyset$, then it is a holomorphically symplectic quasi-projective manifold of dimension $v^{2}+2$ (see \cite{M1}). If $S$ is Abelian, $v^{2}>0$ and $K_{v}^{s}\neq\emptyset$, then $K_{v}^{s}$ is a holomorphically symplectic quasi-projective manifold of dimension $v^{2}-2$ (see \cite{Y2}).

We write $v=mw$, where $m\in\mathbb{N}$ and $w$ is a primitive Mukai vector on $S$. If $S$ is K3, then $M_{v}^{s}\neq\emptyset$ if and only if $w^{2}\geq-2$ (see Theorem 0.1 of \cite{Y1}), while if $S$ is Abelian, then $M_{v}^{s}\neq\emptyset$ if and only if $w^{2}\geq 0$ (see Theorem 0.1 of \cite{Y2}, and compare with section 2.4 of \cite{KLS}). If $w^{2}>0$, then $M_{v}$ and $K_{v}$ are normal, irreducible projective varieties (see Theorem 4.4 of \cite{KLS} and Remark A.1 of \cite{PR2}).

We recall the following definition (where we let $\mathbb{N}^{*}:=\mathbb{N}\setminus\{0\}$), which was introduced in \cite{PR3}:

\begin{defn}
{\rm Let $S$ be a projective K3 or Abelian surface, $v$ a Mukai vector, $H$ an ample divisor on $S$ and $m,k\in\mathbb{N}^{*}$. We say that $(S,v,H)$ is an $(m,k)-$\textit{triple} if the following conditions are verified:
\begin{enumerate}
 \item the polarization $H$ is primitive and $v-$generic;
 \item we have $v=mw$, where $w$ is primitive and $w^{2}=2k$;
 \item if $w=(0,w_{1},w_{2})$ and $\rho(S)>1$, then $w_{2}\neq 0$.
\end{enumerate}}
\end{defn}

If $(S,v,H)$ is an $(m,k)-$triple, then $M_{v}$ is a nonempty, irreducible, normal projective variety of dimension $2m^{2}k+2$ (see Theorem 4.4 of \cite{KLS}), which is symplectic and whose regular locus is $M^{s}_{v}$. If $S$ is Abelian and $(m,k)\neq(1,1)$, then $K_{v}$ is a nonempty, irreducible, normal projective variety of dimension $2m^{2}k-2$, which is symplectic and whose regular locus is $K^{s}_{v}$. If $(m,k)=(1,1)$, then $M_{v}$ is isomorphic to $S\times\widehat{S}$ and $K_{v}$ is just a point.

The starting point of the present paper is the following classification of the moduli spaces $M_{v}$ and $K_{v}$ whose proof can be found in \cite{PR3} (and see references therein). Suppose first that $S$ is a K3 surface.
\begin{enumerate}
 \item If $k<0$, then $M_{v}$ is either empty (if $k<-1$) or a point (if $k=-1$).
 \item If $k=0$, then $M_{w}$ is a K3 surface and $M_{v}\simeq Sym^{m}(M_{w})$. If $m\geq 2$, then $M_{v}$ is a Namikawa symplectic variety which is not irreducible symplectic. 
 \item If $k>0$ and $m=1$, then $M_{v}$ is an irreducible symplectic manifold which is deformation equivalent to $Hilb^{k+1}(S)$.
 \item If $k=1$ and $m=2$, then $M_{v}$ is an irreducible symplectic variety which has a symplectic resolution $\widetilde{M}_{v}$ of the singularities, and $\widetilde{M}_{v}$ is an irreducible symplectic manifold in the deformation class $OG_{10}$.
 \item In all other cases $M_{v}$ is a locally factorial irreducible symplectic variety with terminal singularities.
\end{enumerate}
If $S$ is K3 and $M_{v}$ is neither empty nor a point, then it is a Namikawa symplectic variety, so $H^{2}(M_{v},\mathbb{Z})$ carries a nondegenerate integral quadratic form of signature $(3,b_{2}-3)$ and a compatible pure weight-two Hodge structure.

Suppose now that $S$ is an Abelian surface.
\begin{enumerate}
 \item If $k<0$, then $M_{v}=K_{v}=\emptyset$. 
 \item If $k=0$, then $M_{w}$ is an Abelian surface and $M_{v}\simeq Sym^{m}(M_{w})$. The fiber $K_{v}$ of the sum morphism $Sym^{m}(M_{w})\longrightarrow M_{w}$ is a point (if $m=1$) or a Namikawa symplectic variety which is not irreducible symplectic (if $m>1$).
 \item If $k>0$ and $m=1$, then $K_{v}$ is either a point (if $k=1$), a K3 surface (if $k=2$) or an irreducible symplectic manifold which is deformation equivalent to $Kum^{k-1}(S)$ (if $k\geq 3$).
 \item If $k=1$ and $m=2$, then $K_{v}$ is an irreducible symplectic variety which has a symplectic resolution $\widetilde{K}_{v}$ of the singularities, and $\widetilde{K}_{v}$ is an irreducible symplectic manifold in the deformation class $OG_{6}$.
 \item In all other cases $K_{v}$ is a locally factorial irreducible symplectic variety with terminal singularities.
\end{enumerate}
If $S$ is Abelian and $K_{v}$ is neither empty nor a point, then it is a Namikawa symplectic variety, so $H^{2}(K_{v},\mathbb{Z})$ carries a nondegenerate integral quadratic form of signature $(3,b_{2}-3)$ and a compatible pure weight-two Hodge structure.

The aim of this paper is to study in detail the lattice and Hodge structures on the second integral cohomology of the moduli spaces $M_{v}(S,H)$ associated with $(m,k)-$triples $(S,v,H)$ where $S$ is a K3 surface, and of the moduli spaces $K_{v}(S,H)$ associated with $(m,k)-$triples $(S,v,H)$ where $S$ is an Abelian surface. 

The first result we show is purely topological, and provides us the rank of $H^{2}(M_{v},\mathbb{Z})$ and of $H^{2}(K_{v},\mathbb{Z})$ for every $(m,k)-$triple (see Theorem \ref{thm:b22}):

\begin{thm}
\label{thm:b2}Let $(S,v,H)$ be an $(m,k)-$triple. 
\begin{enumerate}
 \item If $S$ is K3, then $b_{2}(M_{v})=23$.
 \item If $S$ is Abelian and $(m,k)\neq(1,1),(1,2)$, then $b_{2}(K_{v})=7$.
\end{enumerate}
\end{thm}

This result is well-known for $(1,k)-$triples: in this case if $S$ is K3 then $M_{v}$ is deformation equivalent to $Hilb^{k+1}(S)$, for which $b_{2}=23$, and if $S$ is Abelian then $K_{v}$ is deformation equivalent to $Kum^{k-1}(S)$, for which $b_{2}=7$. For $(2,1)-$triples this is showed in \cite{PR}. 

The strategy of the proof of this result is already used in \cite{R1}: by \cite{PR3} we just need to prove the statement for a particular $(m,k)-$triple, so we use a surface $S$ such that $NS(S)=\mathbb{Z}\cdot h$ where $h$ is the first Chern class of an ample line bundle $H$ with $H^{2}=2k$, and $v=m(0,h,0)$. We then relate the fibrations $\phi:M_{v}(S,H)\longrightarrow|mH|$ and $\psi:M_{u}(S,H)\longrightarrow|mH|$ where $u=(0,mh,1-m^{2}k)$. 

\begin{oss}
\label{oss:symb2}
{\rm If $S$ is K3 and $v=mw$ for $w^{2}\leq 0$, then $M_{v}$ is either empty, a point, or a symmetric product of a K3 surface $X$. In this last case, we have an isomorphism $H^{2}(M_{v},\mathbb{Z})\simeq H^{2}(X,\mathbb{Z})$, hence $b_{2}(M_{v})=22$. If $S$ is Abelian and $v=mw$ for $w^{2}\leq 0$, then $M_{v}$ is either empty or a symmetric product of an Abelian surface $X$. In this last case, any fiber $K$ of the sum morphism $M_{v}\longrightarrow X$ is such that $H^{2}(K,\mathbb{Z})\simeq H^{2}(X,\mathbb{Z})$, so that $b_{2}(K)=6$. Moreover, if $(m,k)=(1,1)$ then $K_{v}$ is a point, while if $(m,k)=(1,2)$ then $K_{v}$ is a K3 surface, hence $b_{2}(K_{v})=22$.} 
\end{oss}

The second step in the description of the Beauville-Namikawa form and of the Hodge structure on $H^{2}(M_{v},\mathbb{Z})$ (resp. $H^{2}(K_{v},\mathbb{Z})$) is to construct a morphism $\lambda_{v}:v^{\perp}\longrightarrow H^{2}(M_{v},\mathbb{Z})$ (resp. $\lambda^{0}_{v}:v^{\perp}\longrightarrow H^{2}(K_{v},\mathbb{Z})$). 

The morphism $\lambda_{v}$ was first constructed by Donaldson for rank 2 vector bundles with trivial determinant. Mukai gave a more general construction for moduli spaces of stable sheaves of rank 2 and any first Chern class. O'Grady \cite{OG1} and Yoshioka \cite{Y1} constructed $\lambda_{v}$ for primitive Mukai vectors, showing that it is always an isomorphism which is a Hodge isometry. In \cite{PR} we gave a construction of $\lambda_{v}$ for $(2,1)-$triples, showing that even in this case it is an isomorphism which is a Hodge isometry.

The general construction we present in this paper is a cohomological (and $K-$theoretical) version of the Le Potier morphism (see Chapter 8 of \cite{HL}): if $R_{v}^{s}$ is the open subset of the Quot-scheme parameterizing stable sheaves, and whose quotient (under the action of a linear group) is $M^{s}_{v}$, we first make use of a universal family $\mathcal{Q}^{s}_{v}$ on $S\times R_{v}^{s}$ to produce a morphism $$\mu_{\mathcal{Q}_{v}^{s}}:v^{\perp}\longrightarrow H^{2}(R_{v}^{s},\mathbb{Z}).$$We show that the class $\mu_{\mathcal{Q}^{s}_{v}}(\alpha)$ descends to a class in $H^{2}(M_{v}^{s},\mathbb{Z})$, which extends uniquely to a class $\lambda_{v}(\alpha)\in H^{2}(M_{v},\mathbb{Z})$. When $S$ is Abelian, the morphism $\lambda^{0}_{v}$ is obtained by the same construction, but composing furthermore with the restriction from $M_{v}$ to $K_{v}$.

Theorem \ref{thm:b2} is a key step in order to prove the main result of the present paper, namely that the morphism $\lambda_{v}$ (resp. $\lambda_{v}^{0}$) is an isomorphism between $v^{\perp}$ and $H^{2}(M_{v},\mathbb{Z})$ if $S$ is K3 (resp. between $v^{\perp}$ and $H^{2}(K_{v},\mathbb{Z})$ if $S$ is Abelian) which is a Hodge isometry. Here on $v^{\perp}$ we consider the Hodge structure induced from the one on the Mukai lattice and the lattice structure given by the Mukai pairing, and on $H^{2}(M_{v},\mathbb{Z})$ (resp. on $H^{2}(K_{v},\mathbb{Z})$) we consider the natural Hodge structure and the lattice structure given by the Beauville-Namikawa form. This result is resumed in the following:

\begin{thm}
\label{thm:main}Let $(S,v,H)$ be an $(m,k)-$triple. 
\begin{enumerate}
 \item If $S$ is K3, the morphism $\lambda_{v}:v^{\perp}\longrightarrow H^{2}(M_{v},\mathbb{Z})$ is an isomorphism of $\mathbb{Z}-$modules and a Hodge isometry.
 \item If $S$ is Abelian and $(m,k)\neq(1,1),(1,2)$, the morphism $\lambda^{0}_{v}:v^{\perp}\longrightarrow H^{2}(K_{v},\mathbb{Z})$ is an isomorphism of $\mathbb{Z}-$modules and a Hodge isometry.
\end{enumerate}
\end{thm}

The proof of this result will be discussed in sections 5 and 6. If $v=mw$, the main idea of the proof is to relate the morphism $\lambda_{v}$ to the morphism $\lambda_{pw}$ for $1\leq p<m$, and to conclude the statement by induction on $m$ since for $(1,k)-$ and $(2,1)-$triples the result is known to hold (see \cite{Y1}, \cite{Y2}, \cite{PR2}).

The final result of this paper is the calculation of the Fujiki constants $C_{M_{v}}$ of $M_{v}$ and $C_{K_{v}}$ of $K_{v}$.

\begin{thm}
\label{thm:fujiki}Let $(S,v,H)$ be an $(m,k)-$triple. 
\begin{enumerate}
 \item If $S$ is K3, then the Fujiki constant of $M_{v}(S,H)$ is $$C_{M_{v}}=\frac{(2m^{2}k+2)!}{(m^{2}k+1)!2^{m^{2}k+1}}.$$
 \item If $S$ is Abelian and $(m,k)\neq(1,1),(1,2)$, the Fujiki constant of $K_{v}(S,H)$ is $$C_{K_{v}}=\frac{(2m^{2}k-2)!m^{2}k}{(m^{2}k-1)!2^{m^{2}k-1}}.$$
\end{enumerate}
\end{thm}

For $(1,k)-$triples the Fujiki constants of $M_{v}$ and $K_{v}$ are computed in \cite{B} (see even \cite{R1}), for $(2,1)-$triples they are computed in \cite{R}, \cite{R1} and \cite{PR}. 

The idea of the proof is as follows. First, recall that the Fujiki constant is invariant under deformation, and that for two different $(m,k)-$triples, the corresponding moduli spaces are analytically locally trivially deformation equivalent. This allows us to reduce to the case of an $(m,k)-$triple $(S,v,H)$ where $S$ is a projective K3 or Abelian surface $S$ such that $NS(S)=\mathbb{Z}\cdot h$ where $h$ is the first Chern class of the polarization $H$ of degree $2k$, and where $v=m(0,h,0)$. We then show that if $S$ is K3 (resp. if $S$ is Abelian) then the Fujiki constant $C_{M_{v}}$ of $M_{v}(S,H)$ (resp. $C_{K_{v}}$ of $K_{v}(S,H)$) equals the Fujiki constant of the moduli space $M_{u}(S,H)$ (resp. $K_{u}(S,H)$), where $u=(0,mh,1-m^{2}k)$. As $u$ is primitive $C_{M_{u}}$ and $C_{K_{u}}$ are known and the statement follows. 

\begin{oss}
\label{oss:fujnew} 
{\rm Theorem \ref{thm:main} implies that if $(S,v,H)$ is an $(m,k)-$triple where $S$ is K3 and $v=mw$, then $H^{2}(M_{v},\mathbb{Z})$ is isomorphic and Hodge isometric to $w^{\perp}$, and if $S$ is Abelian then $H^{2}(K_{v},\mathbb{Z})$ is isomorphic and Hodge isometric to $w^{\perp}$. Theorem \ref{thm:fujiki} shows that the if $S$ is K3 then the Fujiki constant of $M_{v}$ only depends on the dimension of $M_{v}$, and if $S$ is Abelian then the Fujiki constant on $K_{v}$ only depends on the dimension of $K_{v}$.} 
\end{oss}

\section{Moduli spaces associated with $(m,k)-$triples}

In this section we recall basic results on moduli spaces associated with $(m,k)-$triples. We first recall the notion of $v-$generic polarizations and the main results on the moduli spaces we obtained in \cite{PR}. Then we recall the Hodge structure and the lattice structure on $H^{2}(M_{v},\mathbb{Z})$ (and on $H^{2}(K_{v},\mathbb{Z})$ for Abelian surfaces).

\subsection{Notation about moduli spaces}

First we recall some basic elements about the construction of moduli spaces of semistable sheaves that will be used all along the paper: we refer the reader to \cite{HL} for further details, in particular to Chapter 4 therein. 

If $S$ is a projective K3 surface or an Abelian surface, $v$ is a Mukai vector on $S$, and $H$ is a polarization on $S$, we let $P_{v,H}$ be the Hilbert polynomial with respect to $H$ of any coherent sheaf $\mathcal{F}$ with Mukai vector $v$, i. e. for every $n\in\mathbb{N}$ we have $$P_{v,H}(n):=\chi(\mathcal{F}\otimes\mathscr{O}_{S}(nH)).$$This polynomial only depends on $v$ and $H$ by the Hirzebruch-Riemann-Roch Theorem.

There is an integer $l\gg 0$ (depending on $v$) such that every $H-$semistable sheaf $\mathcal{F}$ of Mukai vector $v$ is globally generated and has $h^{0}(\mathcal{F}\otimes\mathscr{O}_{S}(lH))=P_{v,H}(l)$. In what follows, we will let $N_{v}:=P_{v,H}(l)$ and $V_{v}$ a complex vector space of dimension $N_{v}$. Moreover, we will let $\mathcal{H}_{v}:=V_{v}\otimes\mathscr{O}_{S}(-lH)$.

For every $H-$semistable sheaf $\mathcal{F}$ of Mukai vector $v$ there is then a surjective morphism $\rho:\mathcal{H}_{v}\longrightarrow\mathcal{F}$, defining a point $[\rho:\mathcal{H}_{v}\longrightarrow\mathcal{F}]$ in the Grothendieck Quot-scheme $Quot(\mathcal{H}_{v},P_{v,H})$ parameterizing quotients of $\mathcal{H}_{v}$ with Hilbert polynomial with respect to $H$ equal to $P_{v,H}$.

We will let $R^{ss}_{v}$ be the open subset of $Quot(\mathcal{H}_{v},P_{v,H})$ parameterizing the quotients $[\rho:\mathcal{H}_{v}\longrightarrow\mathcal{F}]$ where $\mathcal{F}$ is $H-$semistable of Mukai vector $v$, and $R^{s}_{v}$ the open subset of $Quot(\mathcal{H}_{v},P_{v,H})$ of quotients where $\mathcal{F}$ is $H-$stable of Mukai vector $v$. Notice that $R_{v}^{s}\subseteq R_{v}^{ss}$.

The group $PGL(N_{v})=PGL(V_{v})=Aut(\mathcal{H}_{v})/\mathbb{C}^{*}$ acts naturally on $Quot(\mathcal{H}_{v},P_{v,H})$, and $R^{ss}_{v}$ and $R^{s}_{v}$ are both invariant under this action. The moduli space $M_{v}(S,H)$ is the universal good quotient of $R^{ss}_{v}$ under the action of $PGL(N_{v})$, and we let $q_{v}:R^{ss}_{v}\longrightarrow M_{v}$ be the quotient morphism. The moduli space $M^{s}_{v}(S,H)$ is a geometric quotient of $R^{s}_{v}$ under the action of $PGL(N_{v})$, and we let $q^{s}_{v}:R^{s}_{v}\longrightarrow M^{s}_{v}$ be the restriction of $q_{v}$ to $R^{s}_{v}$.

We have two open embeddings that will be used all along the paper: the first one is $j_{v}:R^{s}_{v}\longrightarrow R^{ss}_{v}$, and the second one is $i_{v}:M^{s}_{v}\longrightarrow M_{v}$, so that we have the following commutative diagram
\begin{equation}
\label{eq:commquot}
\begin{CD}
R^{s}_{v} @>{j_{v}}>> R^{ss}_{v}\\
@V{q^{s}_{v}}VV          @VV{q_{v}}V\\
M^{s}_{v} @>{i_{v}}>> M_{v}
\end{CD}
\end{equation}

If $S$ is an Abelian surface and $v^{2}>0$ we moreover have the morphism $a_{v}:M_{v}\longrightarrow S\times\widehat{S}$ (see \cite{Y2} or section 2.2 in \cite{PR3}). As already recalled in the introduction, we let $K_{v}:=a_{v}^{-1}(0_{S},\mathscr{O}_{S})$ and $K_{v}^{s}:=K_{v}\cap M_{v}^{s}$. 

We will make use of the following notation: we let $i_{v}^{0}:K_{v}^{s}\longrightarrow K_{v}$ be the open embedding, $R^{ss,0}_{v}:=q_{v}^{-1}(K_{v})$, $R^{s,0}_{v}:=R^{ss,0}_{v}\cap R^{s}_{v}$ and $j^{0}_{v}:R^{s,0}_{v}\longrightarrow R^{ss,0}_{v}$ the open embedding. Finally, we let $q^{0}_{v}:R^{ss,0}_{v}\longrightarrow K_{v}$ and $q^{s,0}_{v}:R^{s,0}_{v}\longrightarrow K^{s}_{v}$ for the restrictions of $q_{v}$ to $R^{ss,0}_{v}$ and of $q^{s}_{v}$ to $R^{s,0}_{v}$, respectively. We have a commutative diagram
\begin{equation}
\label{eq:commquotab}
\begin{CD}
R^{s,0}_{v} @>{j^{0}_{v}}>> R^{ss,0}_{v}\\
@V{q^{s,0}_{v}}VV          @VV{q^{0}_{v}}V\\
K^{s}_{v} @>{i^{0}_{v}}>> K_{v}
\end{CD}
\end{equation}

\subsection{Generic polarizations}

We now recall the definition of $v-$generic polarization for a Mukai vector $v$ on a projective K3 or Abelian surface $S$. We let $v=(v_{0},v_{1},v_{2})$, and if $v_{0}>0$ we let $|v|:=\frac{v_{0}^{2}}{4}(v,v)+\frac{v_{0}^{2+2\epsilon(S)}}{2}$ and $$W_{v}:=\{D\in NS(S)\,|\,-|v|\leq D^{2}<0\}.$$

If $v_{0}=0$ and $\rho(S)>1$ (so that $v_{2}\neq 0$), with each pair $(\mathscr{E},\mathscr{F})$ of a pure coherent sheaf $\mathscr{E}$ with Mukai vector $v$ and coherent subsheaf $\mathscr{F}$ of $\mathscr{E}$, we associate a divisor $D=u_{2}v_{1}-v_{2}u_{1}$, where $(0,u_{1},u_{2})$ is the Mukai vector of $\mathscr{F}$. We let $W_{v}$ be the set of nontrivial divisors obtained in this way. 

\begin{defn}
The polarization $H$ is $v-$generic if and only if $H\cdot D\neq 0$ for every $D\in W_{v}$.  
\end{defn}

We refer the reader to section 2.1 of \cite{PR3} for the main properties of $v-$generic polarizations. For further purposes, the following result relating generic polarizations for different multiples of the same Mukai vector will be useful.

\begin{lem}
\label{lem:genvw}
Let $v$ be a Mukai vector on a projective K3 or Abelian surface $S$, and write $v=mw$ where $m\in\mathbb{N}$ and $w$ is a primitive Mukai vector. If $H$ is a $v-$generic polarization, then it is $pw-$generic for every $1\leq p\leq m$. In particular, if $(S,v,H)$ is an $(m,k)-$triple, then $(S,pw,H)$ is a $(p,k)-$triple.
\end{lem}

\proof We write $v=(v_{0},v_{1},v_{2})$ and $w=(w_{0},w_{1},w_{2})$. If $v_{0},w_{0}\neq 0$, notice that $$|v|=\frac{m^{2}}{p^{2}}\bigg(\frac{m^{2}}{p^{2}}\cdot\frac{(pw_{0})^{2}}{4}(pw,pw)+\bigg(\frac{m}{p}\bigg)^{2\epsilon(S)}\frac{(pw_{0})^{2+2\epsilon(S)}}{2}\bigg).$$As $m\geq p$ we get that $|v|\geq|pw|$, so $W_{pw}\subseteq W_{v}$ and the statement follows.

If $v_{0}=w_{0}=0$ and $v_{2}=mw_{2}\neq 0$, let $\mathscr{E}$ be a pure sheaf of Mukai vector $pw$ and $\mathscr{F}\subseteq\mathscr{E}$ a subsheaf with Mukai vector $u=(0,u_{1},u_{2})$. Letting $D':=u_{2}w_{1}-w_{2}u_{1}$, the divisor associated with $(\mathscr{E},\mathscr{F})$ is then $pD'$. Consider now a pure sheaf $\mathscr{G}$ of Mukai vector $(m-p)w$, so $\mathscr{E}\oplus\mathscr{G}$ is a pure sheaf of Mukai vector $v$ and $\mathscr{F}\oplus\mathscr{G}$ is a subsheaf of $\mathscr{E}\oplus\mathscr{G}$ of Mukai vector $u+(m-p)w$. The divisor associated with the pair $(\mathscr{E}\oplus\mathscr{G},\mathscr{F}\oplus\mathscr{G})$ is $mD'$. 

This shows us that we have a natural injection $W_{pw}\longrightarrow W_{v}$ mapping $D\in W_{v}$ to $\frac{m}{p}D$. If now $H$ is $v-$generic, then $H\cdot\Delta\neq 0$ for every $\Delta\in W_{v}$. If $D\in W_{pw}$ then $\frac{m}{p}D\in W_{v}$, so $H\cdot\frac{m}{p}D\neq 0$, and hence $H\cdot D\neq 0$, showing that $H$ is $pw-$generic.\endproof

\subsection{Results on moduli spaces}

If $(S,v,H)$ is an $(m,k)-$triple, in \cite{PR3} we proved the following result, which is the starting point of our investigation:

\begin{thm}
\label{thm:mio}Let $(S,v,H)$ be a $(m,k)-$triple.
\begin{enumerate}
 \item The analytically locally trivial deformation class of $M_{v}$ only depends on $m$ and $k$. If $S$ is Abelian, the analytically locally trivial deformation class of $K_{v}$ only depends on $m$ and $k$. 
 \item If $S$ is K3, the moduli spaces $M_{v}$ and $M^{s}_{v}$ are simply connected. If $S$ is Abelian then $K_{v}$ is simply connected, and if $(m,k)\neq(2,1)$ then $K^{s}_{v}$ is simply connected, while if $(m,k)=(2,1)$ then $\pi_{1}(K_{v})=\mathbb{Z}/2\mathbb{Z}$.
 \item If $S$ is K3, then the moduli space $M_{v}$ is an irreducible symplectic variety. If $S$ is Abelian and $(m,k)\neq(1,1)$, then $K_{v}$ is an irreducible symplectic variety.
\end{enumerate}
\end{thm}

Let us then resume the main informations we have on the moduli spaces associated with $(m,k)-$triples.
\begin{enumerate}
 \item If $m=1$, the proof of Theorem \ref{thm:mio} is contained in \cite{OG1} and \cite{Y1}.
 \begin{itemize} 
  \item If $S$ is K3 then $M_{v}$ is an irreducible symplectic manifold which is analytically locally trivially deformation equivalent to $Hilb^{k+1}(S)$. It follows that $H^{2}(M_{v},\mathbb{Z})$ is a free $\mathbb{Z}-$module of rank 22 if $k=1$ and 23 if $k\geq 2$.
	\item If $S$ is Abelian and $k\neq 1$, then $K_{v}$ is a $K3$ surface if $k=2$, and an irreducible symplectic manifold which is analytically locally trivially deformation equivalent to $Kum^{k-1}(S)$ if $k\geq 3$. It follows that $H^{2}(K_{v},\mathbb{Z})$ is a free $\mathbb{Z}-$module of rank 22 if $k=2$ and 7 if $k\geq 3$.
 \end{itemize}
We notice that in this case we find one deformation class in dimension 2 (namely the deformation class of K3 surfaces) and two different deformation classes in dimension at least 4.
	\item If $(m,k)=(2,1)$, the proof of Theorem \ref{thm:mio} is in \cite{PR} and \cite{PR3}. 
 \begin{itemize}
	\item If $S$ is K3 then $M_{v}$ is a $10-$dimensional irreducible symplectic variety whose singular locus has codimension 2, it is either locally factorial or $2-$factorial (by Theorem 1.1 of \cite{PR2}) and has a symplectic resolution $\widetilde{M}_{v}$ which is an irreducible symplectic manifold whose analytically locally trivial deformation class is $OG_{10}$. Finally $H^{2}(M_{v},\mathbb{Z})$ is a free $\mathbb{Z}-$module of rank 23.
  \item	If $S$ is Abelian then $K_{v}$ is a $6-$dimensional irreducible symplectic variety whose singular locus has codimension 2, it is $2-$factorial (by Theorem 1.2 of \cite{PR2}) and has a symplectic resolution $\widetilde{K}_{v}$ which is an irreducible symplectic manifold whose analytically locally trivial deformation class is $OG_{6}$. Finally $H^{2}(K_{v},\mathbb{Z})$ is a free $\mathbb{Z}-$module of rank 7.
 \end{itemize}
 \item For all other cases, the proof of Theorem \ref{thm:mio} is contained in \cite{PR3}. More precisely, point (1) was first shown in \cite{Y4} (see even Theorem 1.3 of \cite{PR3}), point (2) is Lemma 3.1 and Theorem 3.2 of \cite{PR3} and point (3) is Theorem 1.5 of \cite{PR3}. 
 \begin{itemize}
  \item If $S$ is K3 then $M_{v}$ is a $(2km^{2}+2)-$dimensional irreducible symplectic variety whose singular locus has codimension at least 4, and its analytically locally trivial deformation class only depends on $m$ and $k$. The moduli space $M_{v}$ is locally factorial and has no symplectic resolutions of the singularities (see Theorems A and B of \cite{KLS}), and none of its smooth birational models is an irreducible symplectic manifold (see \cite{GLR}, Proposition 6.4). The $\mathbb{Z}-$module $H^{2}(M_{v},\mathbb{Z})$ is free (by Theorem I of \cite{GGK}).
	\item If $S$ is Abelian then $K_{v}$ is a $(2km^{2}-2)-$dimensional irreducible symplectic variety whose singular locus has codimension at least 4, and the analytically locally trivial deformation class only depends on $m$ and $k$. The moduli space $K_{v}$ is locally factorial (see Proposition A.2 of \cite{PR2}), it has no symplectic resolutions of the singularities (see Theorem 6.2 of \cite{KLS}) and none of its smooth birational models is an irreducible symplectic manifold (see \cite{GLR}, Proposition 6.4). The $\mathbb{Z}-$module $H^{2}(K_{v},\mathbb{Z})$ is free (by Theorem I of \cite{GGK}).
 \end{itemize}
\end{enumerate}

We notice that besides the cases of $(1,k)-$ and $(2,1)-$triples we do not know the rank of the free $\mathbb{Z}-$modules $H^{2}(M_{v},\mathbb{Z})$ and $H^{2}(K_{v},\mathbb{Z})$. As explained in the introduction, the calculation of these ranks will be one of the main results of the present paper (see Theorems \ref{thm:b2} and \ref{thm:b22}).

\subsection{The structure of $H^{2}(M_{v},\mathbb{Z})$}

Let $(S,v,H)$ be an $(m,k)-$triple. By point (2) of Theorem \ref{thm:mio}, if $S$ is K3 the moduli space $M_{v}$ is simply connected, and if $S$ is Abelian then $K_{v}$ is simply connected. In any case, it follows that and if $S$ is K3 then $H^{2}(M_{v},\mathbb{Z})$ is a free $\mathbb{Z}-$module, and if $S$ is Abelian then $H^{2}(K_{v},\mathbb{Z})$ is a free $\mathbb{Z}-$module. 

Moreover, we know that $M_{v}$ and $K_{v}$ are irreducible symplectic varieties (with the exception of $K_{v}$ when $(m,k)=(1,1)$): they are smooth if $m=1$, and they have canonical non-terminal singularities if $(m,k)=(2,1)$. In all other cases, they have terminal singularities by Corollary 1 of \cite{N3} (as the singular locus of $M_{v}$ has codimension at least 4), and they are locally factorial by \cite{KLS} and \cite{PR2}.

In any case, the free $\mathbb{Z}-$modules $H^{2}(M_{v},\mathbb{Z})$ and $H^{2}(K_{v},\mathbb{Z})$ carry a pure weight-two Hodge structure and a lattice structure, that we describe here (we refer the reader to \cite{PR3} for more details).

\subsubsection{The Hodge structure.}

If the triple $(S,v,H)$ is a $(1,k)-$triple, then $M_{v}$ and $K_{v}$ are smooth projective manifolds. Hence $H^{2}(M_{v},\mathbb{Z})$ and $H^{2}(K_{v},\mathbb{Z})$ carry a natural pure weight-two Hodge structure. For $M_{v}$ this is $$H^{2}(M_{v},\mathbb{C})=H^{2,0}(M_{v})\oplus H^{1,1}(M_{v})\oplus H^{0,2}(M_{v}),$$where $H^{2,0}=H^{0}(M_{v},\Omega_{M_{v}})$, $H^{1,1}(M_{v})=H^{1}(M_{v},\Omega_{M_{v}})$ and $H^{0,2}(M_{v})=H^{2}(M_{v},\mathscr{O}_{M_{v}})$. If $S$ is Abelian, the same formula holds for $K_{v}$. 

If $m\geq 2$, the description of the Hodge structure on the second integral cohomology on $H^{2}(M_{v},\mathbb{Z})$ and $H^{2}(K_{v},\mathbb{Z})$ is similar, and we present here only the one for $H^{2}(M_{v},\mathbb{Z})$.

Let $\pi:\widetilde{M}_{v}\longrightarrow M_{v}$ be a resolution of the singularities: as the singularities are canonical, they are all rational singularities (see \cite{E}), and hence $$\pi^{*}:H^{2}(M_{v},\mathbb{Z})\longrightarrow H^{2}(\widetilde{M}_{v},\mathbb{Z})$$is an injective morphism of mixed Hodge structures. 

As $\widetilde{M}_{v}$ is smooth and projective, the mixed Hodge structure on $H^{2}(\widetilde{M}_{v},\mathbb{Z})$ is pure of weight two, hence the one on $H^{2}(M_{v},\mathbb{Z})$ is pure of weight two, and we have $$H^{2,0}(M_{v}):=(\pi^{*})^{-1}(H^{2,0}(\widetilde{M}_{v})\cap\pi^{*}(H^{2}(M_{v},\mathbb{C}))),$$ $$H^{1,1}(M_{v}):=(\pi^{*})^{-1}(H^{1,1}(\widetilde{M}_{v})\cap\pi^{*}(H^{2}(M_{v},\mathbb{C}))),$$ $$H^{0,2}(M_{v}):=(\pi^{*})^{-1}(H^{0,2}(\widetilde{M}_{v})\cap\pi^{*}(H^{2}(M_{v},\mathbb{C}))).$$

\subsubsection{The lattice structure}

The second integral cohomology of an irreducible symplectic variety may be endowed with a lattice structure with respect to an integral nondegenerate quadratic form, at least in some cases that we resume in this section. As the moduli space $M_{v}$ (resp. $K_{v}$) associated with an $(m,k)-$triple $(S,v,H)$ where $S$ is K3 (resp. $S$ is Abelian) will fall in one of these cases, this will allow us to endow its second integral cohomology with the lattice structure we will describe. 

The first case is when $X$ is a smooth irreducible symplectic variety, i. e. an irreducible symplectic manifold. On $H^{2}(X,\mathbb{Z})$ we then have the Beauville form $b_{X}$ of $X$, whose signature is $(3,b_{2}(X)-3)$. With respect to the bilinear form associated with $b_{X}$ we have that $H^{2,0}(X)\oplus H^{0,2}(X)$ is orthogonal to $H^{1,1}(X)$. By \cite{F} there is a positive rational number $C_{X}$, the Fujiki constant of $X$, such that if $2n$ is the complex dimension of $X$ we have $$\int_{X}\alpha^{2n}=C_{X}b_{X}(\alpha)^{n}$$for every $\alpha\in H^{2}(X,\mathbb{Z})$. Both the Beauville form and the Fujiki constant of $X$ are deformation invariant.

The second case is when $X$ is a singular irreducible symplectic variety having a symplectic resolution of the singularities $\pi:\widetilde{X}\longrightarrow X$. Remark 1.6 of \cite{PR3} implies that $\widetilde{X}$ is an irreducible symplectic manifold, so on $H^{2}(\widetilde{X},\mathbb{Z})$ we have the Beauville form $b_{\widetilde{X}}$. 

Since the pull-back morphism $\pi^{*}:H^{2}(X,\mathbb{Z})\longrightarrow H^{2}(\widetilde{X},\mathbb{Z})$ is injective, we may define an integral nondegenerate quadratic form $b_{X}$ on $H^{2}(X,\mathbb{Z})$ by using $b_{\widetilde{X}}$: for every $\alpha\in H^{2}(X,\mathbb{Z})$ we let $$b_{X}(\alpha):=b_{\widetilde{X}}(\pi^{*}(\alpha)),$$so that $b_{X}$ is a nondegenerate quadratic form. The signature of $b_{X}$ is again $(3,b_{2}(X)-3)$, and we have that $H^{2,0}(X)\oplus H^{0,2}(X)$ is orthogonal to $H^{1,1}(X)$ with respect to the bilinear form associated with $b_{X}$. If $2n$ is the complex dimension of $X$, then we have $$\int_{X}\alpha^{2n}=C_{\widetilde{X}}b_{X}(\alpha)^{n}$$for every $\alpha\in H^{2}(X,\mathbb{Z})$.

The last case is when $X$ is an irreducible symplectic variety having $\mathbb{Q}-$factorial terminal singularities. Then $X$ is a Namikawa symplectic variety by Proposition 1.10 of \cite{PR3}, and by Corollary 1 of \cite{N3} its singular locus has codimension at least 4. By Theorem 8 of \cite{N1} we conclude that on $H^{2}(X,\mathbb{R})$ there is a nondegenerate quadratic form. 

Up to renormalization, it can be shown exactly as in the smooth case treated by \cite{B} that it defines an integral nondegenerate quadratic form $b_{X}$ on $H^{2}(X,\mathbb{Z})$, that we call the \textit{Beauville-Namikawa form} of $X$. Its signature is again $(3,b_{2}(X)-3)$, and we have that $H^{2,0}(X)\oplus H^{0,2}(X)$ is orthogonal to $H^{1,1}(X)$ with respect to the bilinear form associated with $b_{X}$ (see Corollary 8 of \cite{N4}). By \cite{Ma3}, \cite{Ma4} and \cite{S} if $2n$ is the complex dimension of $X$, then there is a positive rational number $C_{X}$, the Fujiki constant of $X$, such that $$\int_{X}\alpha^{2n}=C_{X}b_{X}(\alpha)^{n}$$for every $\alpha\in H^{2}(X,\mathbb{Z})$.

Now, let us suppose that $(S,v,H)$ is an $(m,k)-$triple. As a consequence of Theorem \ref{thm:mio} and of the previous discussion, we may resume the lattice structure on the second integral cohomology of $M_{v}$ or $K_{v}$ as follows.

Suppose first that $S$ is a K3 surface.
\begin{enumerate}
 \item If $m=1$, then $M_{v}$ is an irreducible symplectic manifold, so on $H^{2}(M_{v},\mathbb{Z})$ we have the Beauville form $b_{v}$ of $M_{v}$. Its signature is $(3,19)$ if $k=1$ and $(3,20)$ if $k>1$, and the value of the Fujiki constant $C_{v}$ of $M_{v}$ is $\frac{(2k+2)!}{(k+1)!2^{k+1}}$ (see \cite{B}).
 \item If $(m,k)=(2,1)$, then by \cite{OG2}, \cite{LS}, \cite{PR} and \cite{PR3} we know that $M_{v}$ is an irreducible symplectic variety having a symplectic resolution $\widetilde{M}_{v}$ of the singularities, so on $H^{2}(M_{v},\mathbb{Z})$ we have the integral nondegenerate quadratic form $b_{v}$ induced by the Beauville form of $\widetilde{M}_{v}$. The signature of $b_{v}$ is again $(3,20)$ (by \cite{PR}), and the Fujiki constant $C_{v}$ of $M_{v}$ in this case is 945, as shown in \cite{R1}.
 \item In all other cases, by point (3) of Theorem \ref{thm:mio} we know that $M_{v}$ is an irreducible symplectic variety, which is moreover locally factorial by \cite{KLS}. Its singular locus has codimension at least 4, so on $H^{2}(M_{v},\mathbb{Z})$ we have the Beauville-Namikawa form $b_{v}$. 
\end{enumerate}

Suppose now that $S$ is an Abelian surface.
\begin{enumerate}
 \item If $m=1$ and $k>2$, then $K_{v}$ is an irreducible symplectic manifold, so on $H^{2}(K_{v},\mathbb{Z})$ we have the Beauville form $b_{v}$ of $K_{v}$. Its signature is $(3,4)$, and the Fujiki constant $D_{v}$ of $K_{v}$ is $\frac{(2k-2)!k}{(k-1)!2^{k-1}}$ (see \cite{B}).
 \item If $(m,k)=(2,1)$, then by \cite{OG3}, \cite{LS}, \cite{PR} and \cite{PR3} we know that $K_{v}$ is an irreducible symplectic variety having a symplectic resolution $\widetilde{K}_{v}$ of the singularities, so on $H^{2}(K_{v},\mathbb{Z})$ we have the integral nondegenerate quadratic form $b_{v}$ induced by the Beauville form of $\widetilde{K}_{v}$. The signature of $b_{v}$ is again $(3,4)$ (by \cite{PR}), and the Fujiki constant $D_{v}$ of $K_{v}$ in this case is 60, as shown in \cite{R}.
 \item In all other cases, by point (3) of Theorem \ref{thm:mio} we know that $K_{v}$ is an irreducible symplectic variety which is locally factorial by \cite{PR2}. Its singular locus has codimension at least 4, so on $H^{2}(K_{v},\mathbb{Z})$ we have the Beauville-Namikawa form $b_{v}$. 
\end{enumerate}

We notice that besides the cases of $(1,k)-$ and $(2,1)-$triples we do not know the value of the Fujiki constants of $M_{v}$ and $K_{v}$. As explained in the introduction, the calculation of these Fujiki constants will be one of the main results of the present paper (see Theorem \ref{thm:fujiki}).

\section{Cohomological properties of the moduli spaces}

This section is devoted to prove some preliminary results on the second integral cohomology of the moduli spaces associated with $(m,k)-$triples that will be used all along the second part of the paper. 

The content of the section is the following: we first summarize some basic results about the first cohomology groups of $M_{v}$, $M_{v}^{s}$, $K_{v}$ and $K_{v}^{s}$ which are a consequence of Theorem \ref{thm:mio}. A second step will be to show that the restriction morphism from $H^{2}(M_{v},\mathbb{C})$ to $H^{2}(M^{s}_{v},\mathbb{C})$ is an isomorphism (of pure weight-two Hodge structures). As final step we will show that the restriction morphism from $H^{2}(M_{v},\mathbb{Z})\longrightarrow H^{2}(M^{s}_{v},\mathbb{Z})$ is an isomorphism (of $\mathbb{Z}-$modules) if $(m,k)\neq(2,1)$. Similar results will be shown for $K_{v}$ and $K^{s}_{v}$ in the case of Abelian surfaces.

\subsection{First cohomology groups}

We first summarize some useful properties of the first cohomology groups of the various moduli spaces we are interested in, which are a consequence of Theorem \ref{thm:mio}.

\begin{lem}
\label{lem:firstcohom}
Let $(S,v,H)$ be an $(m,k)-$triple.
\begin{enumerate}
 \item If $S$ is K3 then $H^{1}(M_{v},\mathbb{Z})=H^{1}(M^{s}_{v},\mathbb{Z})=0$.
 \item If $S$ is Abelian then $H^{1}(K_{v},\mathbb{Z})=H^{1}(K^{s}_{v},\mathbb{Z})=0$. Moreover, we have that $H^{1}(M_{v},\mathbb{Z})$ and $H^{1}(M^{s}_{v},\mathbb{Z})$ are free and the restriction morphism $i_{v}^{*}:H^{1}(M_{v},\mathbb{Z})\longrightarrow H^{1}(M_{v}^{s},\mathbb{Z})$ is an isomorphism.
\end{enumerate}
\end{lem}

\proof Point (1) is an immediate consequence of the fact that both $M_{v}$ and $M_{v}^{s}$ are simply connected (see point (2) of Theorem \ref{thm:mio}). If $S$ is Abelian, by point (2) of Theorem \ref{thm:mio} we have that $K_{v}$ is simply connected and $K^{s}_{v}$ is either simply connected or $\pi_{1}(K^{s}_{v})\simeq\mathbb{Z}/2\mathbb{Z}$. In both cases we have $$H^{1}(K_{v},\mathbb{Z})=H^{1}(K^{s}_{v},\mathbb{Z})=0.$$

Consider now the Leray spectral sequence associated with the morphisms $a_{v}:M_{v}\longrightarrow S\times\widehat{S}$ and its restriction $a_{v}^{s}:M^{s}_{v}\longrightarrow S\times\widehat{S}$: we get a commutative diagram 
$$
\begin{CD}
0 @>{}>> H^{1}(S\times \hat{S},\mathbb{Z})    @>{a^{*}_{v}}>>     H^1(M_v,\mathbb{Z})    @>{}>>  H^0(R^1a_{v*}(\mathbb{Z}))\\
@V{}VV     @V{{\rm id}}VV                                            @VV{i_{v}^{*}}V                 @V{}VV\\
0 @>{}>> H^{1}(S\times \hat{S},\mathbb{Z})    @>>{a^{s*}_{v}}>    H^1(M^s_v,\mathbb{Z})  @>{}>>  H^0(R^1 a^s_{v*}(\mathbb{Z}))
\end{CD}
$$
As $H^{1}(K_{v},\mathbb{Z})=H^{1}(K^{s}_{v},\mathbb{Z})=0$ we get $H^{0}(R^{1}a_{v*}\mathbb{Z})=H^{1}(R^{1}a^{s}_{v*}\mathbb{Z})=0$ and hence $a_{v}^{*}$ and $a^{s*}_{v}$ are isomorphisms. Then $i_{v}^{*}:H^{1}(M_{v},\mathbb{Z})\longrightarrow H^{1}(M^{s}_{v},\mathbb{Z})$ is an isomorphism, and $H^{2}(M_{v},\mathbb{Z})$ and $H^{2}(M^{s}_{v},\mathbb{Z})$ are both free.\endproof

\subsection{Relation between $H^{2}(M_{v},\mathbb{C})$ and $H^{2}(M^{s}_{v},\mathbb{C})$}

We first relate $H^{2}(M_{v},\mathbb{C})$ and $H^{2}(M_{v}^{s},\mathbb{C})$ via the pull-back under the inclusion morphism $i_{v}:M^{s}_{v}\longrightarrow M_{v}$. This is based on the following general result, due to Namikawa:

\begin{lem}
\label{lem:inclqfatt}
Let $X$ be a $\mathbb{Q}-$factorial Namikawa symplectic variety having terminal singularities, let $X^{s}$ be its smooth locus and $i:X^{s}\longrightarrow X$ the inclusion. Then the mixed Hodge structure on $H^{2}(X^{s},\mathbb{C})$ is pure of weight two, and the morphism $$i^{*}:H^{2}(X,\mathbb{C})\longrightarrow H^{2}(X^{s},\mathbb{C})$$is an isomorphism of pure Hodge structures.
\end{lem}

\proof Theorem 1 of \cite{Oh} proves that the mixed Hodge structure on $H^{2}(X^{s},\mathbb{C})$ is pure of weight two (for an algebraic proof of this see Lemma 2.6 of \cite{N2}). The proof that $i^{*}$ is an isomorphism of pure Hodge structures can be found in the proof of Proposition 9 of \cite{N1} (see part (b) of the proof, page 143).\endproof

As a consequence of Lemma \ref{lem:inclqfatt} we get the following, which will be used in the proof of Theorem \ref{thm:b2} and which shows that the restriction morphism from $M_{v}$ to $M^{s}_{v}$ (respectively from $K_{v}$ to $K_{v}^{s}$) induces an isomorphism of pure Hodge structures in cohomology for all $(m,k)-$triples.

\begin{cor}
\label{cor:restriction}Let $(S,v,H)$ be an $(m,k)-$triple. 
\begin{enumerate}
 \item If $i_{v}:M^{s}_{v}\longrightarrow M_{v}$ is the inclusion, the restriction morphism $$i_{v}^{*}:H^{2}(M_{v},\mathbb{C})\longrightarrow H^{2}(M_{v}^{s},\mathbb{C})$$is an isomorphism of pure weight two Hodge structures.
 \item If $S$ is Abelian and $i_{v}^{0}:K^{s}_{v}\longrightarrow K_{v}$ is the inclusion, the restriction morphism $$(i_{v}^{0})^{*}:H^{2}(K_{v},\mathbb{C})\longrightarrow H^{2}(K_{v}^{s},\mathbb{C})$$is an isomorphism of pure weight two Hodge structures.
\end{enumerate}
\end{cor}

\proof The structure of the proof is the following: we first prove point (1) for K3 surfaces; we then turn to Abelian surfaces and we first prove point (2), and deduce point (1) from it.

\textbf{Proof of point (1) for K3 surfaces}. Suppose first that $S$ is a K3 surface. The statement is trivial for $(1,k)-$triples. 

For $(m,k)-$triples with $m>1$ and $(m,k)\neq(2,1)$, then $M_{v}$ is an irreducible symplectic variety by Theorem 1.19 of \cite{PR3}, and hence a Namikawa symplectic variety by Proposition 1.10 of \cite{PR3}. Moreover $M_{v}$ is locally factorial by Theorem A of \cite{KLS} and has terminal singularities (by Corollary 1 of \cite{N3}, as $M_{v}$ is a symplectic variety whose singular locus has codimension at least 4). The statement then follows from Lemma \ref{lem:inclqfatt}. 

We are left with $(2,1)-$triples, so suppose that $(S,v,H)$ is a $(2,1)-$triple where $S$ is a K3 surface, and we write $v=2w$ where $w$ is primitive and $w^{2}=2$. By Lemma 3.7 of \cite{PR} we have that $i_{v}^{*}$ is injective, so if we show that $b_{2}(M_{v})=b_{2}(M^{s}_{v})$ we will be done: the injectivity of $i_{v}^{*}$ implies then that $i_{v}^{*}$ is an isomorphism, and since it is a morphism of mixed Hodge structures and the Hodge structure on $H^{2}(M_{v},\mathbb{C})$ is pure, the Hodge structure on $H^{2}(M^{s}_{v},\mathbb{C})$ will be pure as well. 

We then just need to prove that $b_{2}(M_{v})=b_{2}(M_{v}^{s})$, hence by Theorem 1.7 of \cite{PR} we just need to prove that $b_{2}(M^{s}_{v})=23$. To do so, let $\Sigma$ be the singular locus of $M_{v}$, $\Omega$ the singular locus of $\Sigma$, $\Sigma^{0}:=\Sigma\setminus\Omega$ and $\pi:\widetilde{M}\longrightarrow M_{v}$ be the blow-up of $M_{v}$ along $\Sigma$ with reduced structure. We let $\widetilde{\Sigma}$ be the exceptional divisor of $\pi$, $\widetilde{\Omega}$ the singular locus of $\widetilde{\Sigma}$ and $\widetilde{\Sigma}^{0}:=\widetilde{\Sigma}\setminus\widetilde{\Omega}$. Notice that $\pi^{-1}(M^{s}_{v})\subseteq\widetilde{M}\setminus\widetilde{\Omega}$ and that $\pi^{-1}(M^{s}_{v})\simeq M^{s}_{v}$.

Consider the exact sequence $$H^{2}(\widetilde{M}\setminus\widetilde{\Omega},\pi^{-1}(M^{s}_{v}))\rightarrow H^{2}(\widetilde{M}\setminus\widetilde{\Omega})\rightarrow H^{2}(\pi^{-1}(M^{s}_{v}))\rightarrow H^{3}(\widetilde{M}\setminus\widetilde{\Omega},\pi^{-1}(M^{s}_{v}))$$associated with the pair $(\widetilde{M}\setminus\widetilde{\Omega},\pi^{-1}(M^{s}_{v}))$. As $\widetilde{\Sigma}^{0}=(\widetilde{M}\setminus\widetilde{\Omega})\setminus\pi^{-1}(M^{s}_{v})$ is smooth of real codimension 2, by the Thom isomorphism we have $$H^{2}(\widetilde{M}\setminus\widetilde{\Omega},\pi^{-1}(M^{s}_{v}))\simeq H^{0}(\widetilde{\Sigma}^{0},\mathbb{C}),\,\,\,\,\,\,\,H^{3}(\widetilde{M}\setminus\widetilde{\Omega},\pi^{-1}(M^{s}_{v}))\simeq H^{1}(\widetilde{\Sigma}^{0},\mathbb{C}).$$Since $\pi^{-1}(M^{s}_{v})\simeq M^{s}_{v}$ the previous exact sequence may then be written as 
\begin{equation}
\label{eq:exactsigma0}
H^{0}(\widetilde{\Sigma}^{0})\stackrel{\gamma}\longrightarrow H^{2}(\widetilde{M}\setminus\widetilde{\Omega})\longrightarrow H^{2}(M^{s}_{v})\longrightarrow H^{1}(\widetilde{\Sigma}^{0}).
\end{equation}

The codimension of $\widetilde{\Omega}$ in $\widetilde{M}$ is at least 4, so $H^{2}(\widetilde{M}\setminus\widetilde{\Omega})\simeq H^{2}(\widetilde{M})$. Moreover, by \cite{OG2} we know that $\widetilde{\Sigma}$ is irreducible, so $\widetilde{\Sigma}^{0}$ is irreducible as well and $H^{0}(\widetilde{\Sigma}^{0})\simeq\mathbb{C}$. As $\gamma$ maps the generator of $H^{0}(\widetilde{\Sigma}^{0})$ to the first Chern class of $\widetilde{\Sigma}^{0}$, it follows that $\gamma$ is injective.

Now, recall that $\Sigma\simeq Sym^{2}(M_{w})$ where $M_{w}$ is an irreducible symplectic manifold of dimension 4 (see for example \cite{PR}). It follows that $\Sigma^{0}$ is isomorphic to the smooth locus of $Sym^{2}(M_{w})$. Its universal covering is then isomorphic to $(M_{w}\times M_{w})\setminus\Delta$ (where $\Delta$ is the diagonal): as this covering is $2:1$, it follows that $\pi_{1}(\Sigma^{0})=\mathbb{Z}/2\mathbb{Z}$. 

Now, recall by \cite{OG2} that $\pi:\pi^{-1}(\Sigma^{0})\longrightarrow\Sigma^{0}$ is a $\mathbb{P}^{1}-$bundle, hence $\pi_{1}(\pi^{-1}(\Sigma^{0}))\simeq\mathbb{Z}/2\mathbb{Z}$. But $\pi^{-1}(\Sigma^{0})$ is open in $\widetilde{\Sigma}^{0}$, so $\pi_{1}(\pi^{-1}(\Sigma^{0}))$ surjects onto $\pi_{1}(\widetilde{\Sigma}^{0})$, so $\pi_{1}(\widetilde{\Sigma}^{0})$ is torsion and $H^{1}(\widetilde{\Sigma}^{0},\mathbb{C})=0$. The exact sequence (\ref{eq:exactsigma0}) then gives an exact sequence $$0\longrightarrow\mathbb{C}\longrightarrow H^{2}(\widetilde{M},\mathbb{C})\longrightarrow H^{2}(M^{s}_{v},\mathbb{C})\longrightarrow 0.$$Since by \cite{R1} and \cite{PR} we have $b_{2}(\widetilde{M})=24$, we finally get that $b_{2}(M^{s}_{v})=23$, completing the proof of the statement for K3 surfaces.

\textbf{Proof of point (2) for Abelian surfaces}. Suppose now that $S$ is an Abelian surface. The case of $(1,k)-$triples is trivial again, while the case of $(m,k)-$triples where $m>1$ and $(m,k)\neq(2,1)$ follows from Lemma \ref{lem:inclqfatt} as in the case of K3 surfaces: $K_{v}$ is again a Namikawa symplectic variety by Theorem 1.19 and Proposition 1.10 of \cite{PR3}; it is locally factorial by Proposition A.2 of \cite{PR2}, and it has terminal singularities since it is a symplectic variety whose singular locus has codimension at least 4. 

Suppose now that $(S,v,H)$ is a $(2,1)-$triple where $S$ is an Abelian surface, and write $v=2w$ where $w$ is primitive and $w^{2}=2$. Lemma 3.7 of \cite{PR} provides again the injectivity of $(i^{0}_{v})^{*}$, hence we only need to prove that $b_{2}(K_{v})=b_{2}(K^{s}_{v})$. Since by Theorem 1.7 of \cite{PR} we have $b_{2}(K_{v})=7$, we then just need to prove that $b_{2}(K^{s}_{v})=7$. 

Theorem 1.6 of \cite{PR3} shows that if $(S,v,H)$ and $(S',v',H')$ are two $(2,1)-$triples, then $K_{v}^{s}(S,H)$ and $K_{v'}^{s}(S',H')$ are deformation equivalent, hence it is enough to prove that $b_{2}(K^{s}_{v})=7$ for one particular $(2,1)-$triple. We choose it as follows: we let $J$ be the Jacobian of a smooth projective curve $C$ of genus 2, $\Theta$ the theta divisor on $J$ and $v=2(0,\Theta,1)$. It is easy to see that $(J,v,\Theta)$ is a $(2,1)-$triple, and we show that $b_{2}(K^{s}_{v})=7$ in this case.

The argument is similar to the one we used for K3 surfaces, and we use the same notation: let $\Sigma$ be the singular locus of $K_{v}$, $\Omega$ the singular locus of $\Sigma$, $\Sigma^{0}:=\Sigma\setminus\Omega$ and $\pi:\widetilde{K}\longrightarrow K_{v}$ be the blow-up of $K_{v}$ along $\Sigma$ with reduced structure. We let $\widetilde{\Sigma}$ be the exceptional divisor of $\pi$, $\widetilde{\Omega}$ the singular locus of $\widetilde{\Sigma}$ and $\widetilde{\Sigma}^{0}:=\widetilde{\Sigma}\setminus\widetilde{\Omega}$. Notice that $\pi^{-1}(K^{s}_{v})\subseteq\widetilde{K}\setminus\widetilde{\Omega}$ and that $\pi^{-1}(K^{s}_{v})\simeq K^{s}_{v}$.

As before, since by \cite{OG3} and \cite{PR} we know that $b_{2}(\widetilde{K})=8$, we just need to deduce an exact sequence 
\begin{equation}
\label{eq:exkk}
0\longrightarrow\mathbb{C}\longrightarrow H^{2}(\widetilde{K},\mathbb{C})\longrightarrow H^{2}(K^{s}_{v},\mathbb{C})\longrightarrow 0
\end{equation}
from the exact sequence 
\begin{equation}
\label{eq:exactsigmak0}
H^{0}(\widetilde{\Sigma}^{0})\stackrel{\gamma}\longrightarrow H^{2}(\widetilde{K}\setminus\widetilde{\Omega})\longrightarrow H^{2}(K^{s}_{v})\longrightarrow H^{1}(\widetilde{\Sigma}^{0})
\end{equation}
associated with the pair $(\widetilde{K}\setminus\widetilde{\Omega},\pi^{-1}(K^{s}_{v}))$.

To do so, notice that the fact that $H^{0}(\widetilde{\Sigma}^{0})\simeq\mathbb{C}$ and that $\gamma$ is injective follow again from the fact that $\widetilde{\Sigma}$ is irreducible (see \cite{OG3}). We then just need to prove that $H^{1}(\widetilde{\Sigma}^{0},\mathbb{C})=0$.

To do so, consider the map $f:J\times\widehat{J}\longrightarrow\Sigma$ mapping $(x,L)$ to $t_{x*}L\oplus t_{-x*}L$ (where if $p\in J$ then $t_{p}:J\longrightarrow J$ is the translation by $p$). This is a $2:1$ covering which is \'etale outside the 256 $2-$torsion points of $J\times\widehat{J}$: if $\mathbb{Z}/2\mathbb{Z}$ acts on $J\times\widehat{J}$ as $-1$, the morphism $f$ is invariant with respect to this action. 

We consider the following commutative diagram
\begin{equation}
\label{eq:sigmab}
\begin{CD}
T @>{\pi}>> J\times\widehat{J}\\
@V{g}VV                   @VV{f}V\\
A @>>{\phi}> \Sigma
\end{CD}
\end{equation}
where $\pi:T\longrightarrow J\times\widehat{J}$ is the blow-up of $J\times\widehat{J}$ along its 256 $2-$torsion points and $\phi:A\longrightarrow\Sigma$ is the resolution of the singularities of $\Sigma$ obtained by blowing up its 256 singular points.

The action of $\mathbb{Z}/2\mathbb{Z}$ over $J\times\widehat{J}$ extends naturally to an action on $T$ so that $g:T\longrightarrow A$ is invariant with respect to this action. It follows that $H^{1}(A,\mathbb{C})\simeq H^{1}(T,\mathbb{C})^{\mathbb{Z}/2\mathbb{Z}}=0$. Using now the exact sequence $$H^{1}(A)\longrightarrow H^{1}(\phi^{-1}(\Sigma^{0}))\longrightarrow H^{2}(A,\phi^{-1}(\Sigma^{0}))\stackrel{\beta}\longrightarrow H^{2}(A)$$coming from the pair $(A,\phi^{-1}(\Sigma^{0}))$, since $H^{1}(A)=0$ and the map $\beta$ is injective, it follows that $H^{1}(\phi^{-1}(\Sigma^{0}),\mathbb{C})=0$. But as $\phi:\phi^{-1}(\Sigma^{0})\longrightarrow\Sigma^{0}$ is an isomorphism, it follows that $H^{1}(\Sigma^{0},\mathbb{C})=0$.

Now, recall by \cite{OG2} that $\pi:\pi^{-1}(\Sigma^{0})\longrightarrow\Sigma^{0}$ is a $\mathbb{P}^{1}-$bundle: it follows that $H^{1}(\pi^{-1}(\Sigma^{0}),\mathbb{C})=0$. As $\pi^{-1}(\Sigma^{0})$ is open in $\widetilde{\Sigma}^{0}$, we get $H^{1}(\widetilde{\Sigma}^{0},\mathbb{C})=0$, concluding the proof.

\textbf{Proof of point (1) for Abelian surfaces}. The inclusion morphism $i_v: M_v^s\longrightarrow M_v$ induces a morphism between the Leray spectral sequences of the fibrations $a_v:M_v\longrightarrow S\times\widehat{S}$ and $a^{s}_v:M^{s}_v\longrightarrow S\times\widehat{S}$, where $a^{s}_{v}=a_{v}\circ i_{v}$.

By point (2) of the statement that we just proved for Abelian surfaces, the morphism $i_v$ induces an isomorphism $$H^0(R^2 a_{v*}(\mathbb{C}))\simeq H^0(R^2 a^{s}_{v*}(\mathbb{C})).$$Since, by Theorem 3.4 and Theorem 3.7 of \cite{PR3} we have
$$H^{1}(K_v,\mathbb{C})\simeq H^{1}(K_v^s,\mathbb{C})=0,\,\,\,\,\,\,\,H^{0}(K_v,\mathbb{C})\simeq H^{0}(K_v^s,\mathbb{C})\simeq \mathbb{C},$$the morphism $i_v$ also induces isomorphisms $$H^1(R^1 a_{v*}(\mathbb{C}))\simeq H^1(R^1 a^{s}_{v*}(\mathbb{C})),\,\,\,\,\,\,\,\,\,H^{i}(R^0 a_{v*}(\mathbb{C}))\simeq H^{i}(R^0 a^{s}_{v*}(\mathbb{C}))$$for $i=2,3$. By analyzing the morphism of spectral sequences we obtain that $i_{v}^{*}:H^{2}(M_{v},\mathbb{C})\longrightarrow H^{2}(M_{v}^{s},\mathbb{C})$ is an isomorphism. As it is a morphism of mixed Hodge structures, and as the mixed Hodge structure on $H^{2}(M_{v},\mathbb{C})$ is pure of weight two (see section 2.4.1), it follows that the mixed Hodge structure on $H^{2}(M^{s}_{v},\mathbb{C})$ is pure as well.\endproof
 
\subsection{Relation between $H^{2}(M_{v},\mathbb{Z})$ and $H^{2}(M^{s}_{v},\mathbb{Z})$}

The aim of this section is to show that if $(S,v,H)$ is an $(m,k)-$triple for $(m,k)\neq(2,1)$, then the restriction morphism $H^{2}(M_{v},\mathbb{Z})\longrightarrow H^{2}(M^{s}_{v},\mathbb{Z})$ is an isomorphism. 

We use the following notation: for every $n\in\mathbb{N}$, $n>1$, we let $\mu_{n}$ be the multiplicative group of the $n-$th roots of 1. Moreover, for every variety $X$ we let $Br(X)$ be the Brauer group of $X$, i. e. the group of equivalence classes of Azumaya algebras on $X$ (the interested reader may refer to \cite{Gr} or \cite{Mi2} for the general properties of the Brauer group).

The key result is the following, relating the second cohomology with coefficients in $\mu_{n}$ of a normal projective variety $X$ with that of its smooth locus, and showing that the restriction morphism to a smooth open subset is injective in cohomology if it is injective at the level of Brauer groups. 

\begin{lem}
\label{lem:munrest}
Let $X$ be a locally factorial, quasi-projective variety, $U\subseteq X$ a smooth open subset such that $X\setminus U$ has codimension at least 2 in $X$, and $i:U\longrightarrow X$ the inclusion. Suppose that $i^{*}:Br(X)\longrightarrow Br(U)$ is injective.
\begin{enumerate}
 \item For every $n\in\mathbb{N}$, $n>1$, the morphism $i^{*}:H^{2}(X,\mu_{n})\longrightarrow H^{2}(U,\mu_{n})$ is injective.
 \item The restriction morphism $i^{*}:H^{2}(X,\mathbb{Z})\longrightarrow H^{2}(U,\mathbb{Z})$ is injective. 
 \item If $H^{2}(U,\mathbb{Z})$ is free, then the image of $i^{*}$ is saturated in $H^{2}(U,\mathbb{Z})$. If moreover $H^{2}(X,\mathbb{Z})$ and $H^{2}(U,\mathbb{Z})$ have the same rank, the morphism $i^{*}:H^{2}(X,\mathbb{Z})\longrightarrow H^{2}(U,\mathbb{Z})$ is an isomorphism.
\end{enumerate} 
\end{lem}

\proof By point (ii) of Theorem 4.1 in \cite{SGA4}, Tome 3, Chapitre XVI we know that the \'etale and the classical cohomology of constructible sheaves of torsion Abelian groups are isomorphic. To prove (1) we then just need to show that the restriction morphism $i^{*}:H^{2}_{\acute et}(X,\mu_{n})\longrightarrow H^{2}_{\acute et}(U,\mu_{n})$ is injective.

To do so, for every $n>1$ we consider the following exact sequence of sheaves (with respect to the \'etale topology) $$0\longrightarrow\mu_{n}\longrightarrow\mathcal{O}_{X}^{*}\stackrel{x^{n}}\longrightarrow\mathcal{O}_{X}^{*}\longrightarrow 0,$$which induces the following exact sequence in cohomology $$0\longrightarrow Pic(X)[n]\longrightarrow H^{2}_{\acute et}(X,\mu_{n})\longrightarrow Tors_{n}(H^{2}_{\acute et}(X,\mathcal{O}_{X}^{*}))\longrightarrow 0$$(recall that $Pic(X)\simeq H^{1}_{\acute et}(X,\mathcal{O}_{X}^{*})$), where $Pic(X)[n]$ is the quotient of $Pic(X)$ by the image of the morphism $Pic(X)\longrightarrow Pic(X)$ mapping $L\in Pic(X)$ to $L^{\otimes n}$, and $Tors_{n}(H^{2}_{\acute et}(X,\mathcal{O}_{X}^{*}))$ is the subgroup of $H^{2}_{\acute et}(X,\mathcal{O}_{X}^{*})$ given by $n-$torsion elements.

The same exact sequence is available on $U$, and we have a commutative diagram
$$\begin{array}{ccccccccc}
0 & \longrightarrow & Pic(X)[n] & \longrightarrow & H^{2}_{\acute et}(X,\mu_{n}) & \longrightarrow & Tors_{n}(H^{2}_{\acute et}(X,\mathcal{O}_{X}^{*})) & \longrightarrow & 0\\
  &                 & \downarrow &                & \downarrow              &                 & \downarrow                       
	 &                 &  \\
0 & \longrightarrow & Pic(U)[n] & \longrightarrow & H^{2}_{\acute et}(U,\mu_{n}) & \longrightarrow & Tors_{n}(H^{2}_{\acute et}(U,\mathcal{O}_{U}^{*})) & \longrightarrow & 0
\end{array}$$
whose horizontal rows are exact sequences, and whose vertical arrows are induced by the restriction morphism. 

As $X$ and $U$ are quasi-projective varieties, by Gabber Theorem (see Theorem 1.1 of \cite{dJ}) we have $Tors(H^{2}_{\acute et}(X,\mathcal{O}_{X}^{*}))\simeq Br(X)$ and $Tors(H^{2}_{\acute et}(U,\mathcal{O}_{U}^{*}))\simeq Br(U)$. Hence the previous commutative diagram reads as
\begin{equation}
\label{eq:picbr}
\begin{array}{ccccccccc}
0 & \longrightarrow & Pic(X)[n] & \longrightarrow & H^{2}_{\acute et}(X,\mu_{n}) & \longrightarrow & Br_{n}(X) & \longrightarrow & 0\\
  &                 & \downarrow &                & \downarrow              &                 & \downarrow                       
	 &                 &  \\
0 & \longrightarrow & Pic(U)[n] & \longrightarrow & H^{2}_{\acute et}(U,\mu_{n}) & \longrightarrow & Br_{n}(U) & \longrightarrow & 0
\end{array}
\end{equation}
where $Br_{n}(X)$ (resp. $Br_{n}(U)$) is the $n-$torsion subgroup of $Br(X)$ (resp. of $Br(U)$). 

As $X$ is locally factorial, $U$ is smooth and $X\setminus U$ has codimension at least 2 in $X$, it follows that the left-hand vertical arrow in diagram (\ref{eq:picbr}) is an isomorphism. The right-hand vertical arrow in diagram (\ref{eq:picbr}) is injective by hypothesis: the Snake Lemma then implies that $i^{*}:H^{2}_{\acute et}(X,\mu_{n})\longrightarrow H^{2}_{\acute et}(U,\mu_{n})$ is injective, completing the proof of point (1) of the statement.

We now show that the injectivity of $i^{*}:H^{2}(X,\mu_{n})\longrightarrow H^{2}(U,\mu_{n})$ implies the injectivity of $i^{*}:H^{2}(X,\mathbb{Z})\longrightarrow H^{2}(U,\mathbb{Z})$. By the Universal Coefficient Theorem for cohomology we have the following commutative diagram
\begin{equation}
\label{eq:univcoeff}
\resizebox{.93\hsize}{!}{$\begin{array}{ccccccccc}
0 & \longrightarrow & H^{2}(X,\mathbb{Z})\otimes\mu_{n} & \longrightarrow & H^{2}(X,\mu_{n}) & \longrightarrow & Tor_{1}(H^{3}(X,\mathbb{Z}),\mu_{n}) & \longrightarrow & 0\\
  &                 & \downarrow &                & \downarrow              &                 & \downarrow                       
	 &                 &  \\
0 & \longrightarrow & H^{2}(U,\mathbb{Z})\otimes\mu_{n} & \longrightarrow & H^{2}(U,\mu_{n}) & \longrightarrow & Tor_{1}(H^{3}(U,\mathbb{Z}),\mu_{n}) & \longrightarrow & 0
\end{array}$}
\end{equation}
in which the vertical arrows are restriction morphisms. As we just proved, the central vertical arrow is injective for every $n>1$: it follows that the first vertical arrow is injective for every $n>1$, i. e. the restriction morphism induces an injection $i^{*}_{n}:H^{2}(X,\mathbb{Z})\otimes\mu_{n}\longrightarrow H^{2}(U,\mathbb{Z})\otimes\mu_{n}$ for every $n>1$.

As a consequence of this, the morphism $i^{*}:H^{2}(X,\mathbb{Z})\longrightarrow H^{2}(U,\mathbb{Z})$ is injective. To show this, suppose that $\alpha\in H^{2}(X,\mathbb{Z})$ is such that $i^{*}\alpha=0$ in $H^{2}(U,\mathbb{Z})$. Then for every $n>1$ we have $i^{*}\alpha=0$ in $H^{2}(U,\mathbb{Z})\otimes\mu_{n}$, i. e. $i^{*}_{n}\alpha=0$ for every $n>1$. Since $i^{*}_{n}$ is injective for every $n>1$, it follows that $\alpha=0$ in $H^{2}(X,\mathbb{Z})\otimes\mu_{n}$ for every $n>1$, i. e. $n\alpha=0$ in $H^{2}(X,\mathbb{Z})$ for every $n>1$. But this implies $\alpha=0$ in $H^{2}(X,\mathbb{Z})$, and the proof of (2) is complete.

We are left with the proof of (3). So, suppose that $H^{2}(U,\mathbb{Z})$ is free, we show that the image of $i^{*}$ in $H^{2}(U,\mathbb{Z})$ is saturated. To do so, let $T$ be the quotient of $H^{2}(U,\mathbb{Z})$ by $i^{*}(H^{2}(X,\mathbb{Z}))$, so that we have an exact sequence 
\begin{equation}
\label{eq:tquot}
0\longrightarrow H^{2}(X,\mathbb{Z})\stackrel{i^{*}}\longrightarrow H^{2}(U,\mathbb{Z})\longrightarrow T\longrightarrow 0.
\end{equation}
Now, as $H^{2}(U,\mathbb{Z})$ is free we have $Tor_{1}(H^{2}(U,\mathbb{Z}),\mu_{n})=0$. It follows that tensoring with $\mu_{n}$ the exact sequence (\ref{eq:tquot}) we get an exact sequence $$0\longrightarrow Tor_{1}(T,\mu_{n})\longrightarrow H^{2}(X,\mathbb{Z})\otimes\mu_{n}\stackrel{i^{*}_{n}}\longrightarrow H^{1}(U,\mathbb{Z})\otimes\mu_{n}\longrightarrow T\otimes\mu_{n}\longrightarrow 0$$and as $i^{*}_{n}$ is injective, it follows that $Tor_{1}(T,\mu_{n})=0$. As a consequence $T$ has no $n-$torsion elements for every $n\in\mathbb{N}$, i. e. $T$ is free. But this implies that $i^{*}(H^{2}(X,\mathbb{Z}))$ is saturated in $H^{2}(U,\mathbb{Z})$. The last part of point (3) is an immediate consequence of this.\endproof

We now will use Lemma \ref{lem:munrest} and Proposition \ref{lem:codimrsv} (in the Appendix) to relate the second integral cohomology of the moduli space $M_{v}$ (resp. $K_{v}$) associated with an $(m,k)-$triple $(S,v,H)$ with the second integral cohomology of its smooth locus: this is the content of the following, which is the main result of this section (see section 1.1 for the notation):

\begin{prop}
\label{prop:extmr}Let $(S,v,H)$ be an $(m,k)-$triple where $m\geq 2$ and $(m,k)\neq(2,1)$.
\begin{enumerate}
 \item The morphism $j_{v}^{*}:H^{2}(R^{ss}_{v},\mathbb{Z})\longrightarrow H^{2}(R^{s}_{v},\mathbb{Z})$ is injective, where $j_{v}:R^{s}_{v}\longrightarrow R^{ss}_{v}$ is the inclusion.
 \item The morphism $i_{v}^{*}:H^{2}(M_{v},\mathbb{Z})\longrightarrow H^{2}(M^{s}_{v},\mathbb{Z})$ is an isomorphism, where $i_{v}:M_{v}^{s}\longrightarrow M_{v}$ is the inclusion. 
 \item If $S$ is Abelian the morphism $(i_{v}^{0})^{*}:H^{2}(K_{v},\mathbb{Z})\longrightarrow H^{2}(K_{v}^{s},\mathbb{Z})$ is an isomorphism, where $i_{v}^{0}:K^{s}_{v}\longrightarrow K_{v}$ is the inclusion.
\end{enumerate}
\end{prop}

\proof We prove the three points of the statement separately.

\textbf{Proof of (1)}. By Proposition \ref{lem:codimrsv} the injectivity of $j_{v}^{*}$ is a consequence of point (2) of Lemma \ref{lem:munrest} as soon as we know that $j_{v}^{*}:Br(R^{ss}_{v})\longrightarrow Br(R^{s}_{v})$ is injective. To prove this we show that $j_{v}^{*}:Br_{n}(R^{ss}_{v})\longrightarrow Br_{n}(R^{s}_{v})$ is injective for every $n>1$, so consider $\alpha\in Br_{n}(R^{ss}_{v})$ and suppose that $j_{v}^{*}\alpha=1$. 

As $\alpha\in Br_{n}(R^{ss}_{v})$, there is a projective bundle $\pi:P\longrightarrow R^{ss}_{v}$ of rank $n-1$ whose Brauer class is $\alpha$. Let $\pi_{s}:P^{s}\longrightarrow R^{s}_{v}$ be the restriction of $P$ to $R^{s}_{v}$, so we have a commutative diagram
$$\begin{CD}
\label{CD:commPR}
P^{s} @>{\iota}>> P\\
@V{\pi_{s}}VV          @VV{\pi}V\\
R^{s}_{v} @>{j_{v}}>> R^{ss}_{v}
\end{CD}$$
where $\iota$ is the inclusion. The Brauer class of $P^{s}$ is $j_{v}^{*}\alpha=1$: there is then a holomorphic vector bundle $V$ of rank $n$ over $R^{s}_{v}$ such that $P^{s}\simeq\mathbb{P}(V)$ (as projective bundles).

By Proposition 3.10 of \cite{KLS} we know that $R^{ss}_{v}$ is locally complete intersection and its singular locus has codimension at least 3: as $P$ is a projective bundle on $R^{ss}_{v}$, it follows that $P$ is locally complete intersection and its singular locus has codimension at least 3. By Corollary 3.14 in Exp. XI of \cite{SGA2} we then have that $P$ is locally factorial, and since $P^{s}$ is a smooth open subset of $P$ the restriction morphism gives an isomorphism $Pic(P)\simeq Pic(P^{s})$.

The tautological line bundle on $P^{s}$ (which exists as $P^{s}\simeq\mathbb{P}(V)$) extends then to a line bundle $\mathcal{O}(1)$ on $P$, which fiberwise is the ample line bundle of degree 1. Hence $P$ has a tautological line bundle, so by Lemma 1.4 of \cite{Y5} there is a holomorphic vector bundle $\overline{V}$ on $R^{ss}_{v}$ of rank $n$ such that $P\simeq\mathbb{P}(\overline{V})$ (as projective bundles). But this implies that $\alpha=1$, completing the proof of point (1) of the statement.

\textbf{Proof of (2)}. By Theorem \ref{thm:mio} we know that if $S$ is K3, then $M_{v}$ and $M^{s}_{v}$ are simply connected, while if $S$ is Abelian, then $a_v:M_v\longrightarrow S\times\widehat{S}$ and its restriction to $M^{s}_{v}$ are isotrivial fibrations with simply connected fibers. In both cases $\pi_{1}(M_{v})$ and $\pi_{1}(M^{s}_{v})$ are free Abelian groups, so $H^{2}(M_{v},\mathbb{Z})$ and $H^{2}(M^{s}_{v},\mathbb{Z})$ are free $\mathbb{Z}-$modules of the same rank (by Corollary \ref{cor:restriction}).

Since $M_{v}$ is locally factorial by \cite{KLS}, the statement follows from point (3) of Lemma \ref{lem:munrest} as soon as we know that $i_{v}^{*}:Br(M_{v})\longrightarrow Br(M^{s}_{v})$ is injective, i. e. that $i_{v}^{*}:Br_{n}(M_{v})\longrightarrow Br_{n}(M^{s}_{v})$ is injective for every $n>1$. 

To show this, let $\alpha\in Br_{n}(M_{v})$ be such that $i_{v}^{*}\alpha=1$. Let $\pi:P\longrightarrow M_{v}$ be a projective bundle of rank $n-1$ on $M_{v}$ whose Brauer class is $\alpha$, and let $\pi^{s}:P^{s}\longrightarrow M^{s}_{v}$ be its restriction to $M^{s}_{v}$. The Brauer class of $P^{s}$ is $i_{v}^{*}\alpha=1$: it follows that there is a holomorphic vector bundle $V$ of rank $n$ on $M^{s}_{v}$ such that $P^{s}\simeq\mathbb{P}(V)$ (as projective bundles), so there is a line bundle $L$ on $P^{s}$ whose restriction to every fiber has degree 1.

Now, let $q_{v}:R^{ss}_{v}\longrightarrow M_{v}$ and $q_{v}^{s}:R^{s}_{v}\longrightarrow M^{s}_{v}$ be the quotients with respect to the action of $PGL(N_{v})$. Let $P_{R}:=P\times_{M_{v}}R^{ss}_{v}$, which is a projective bundle of rank $n-1$ over $R^{ss}_{v}$, and let $\pi_{R}:P_{R}\longrightarrow R^{ss}_{v}$ and $q_{R}:P_{R}\longrightarrow P$ be the two projections. Let moreover $P^{s}_{R}:=\pi_{R}^{-1}(R^{s}_{v})$, which is a smooth open subset of $P_{R}$, and $\pi^{s}_{R}:P^{s}_{R}\longrightarrow R^{s}_{v}$ be the restriction of $\pi_{R}$ to $P^{s}_{R}$, which is a projective bundle of rank $n-1$ over $R^{s}_{v}$. Letting $q^{s}_{R}:P^{s}_{R}\longrightarrow P^{s}$ be the restriction of $q_{R}$ to $P^{s}_{R}$, we get a commutative diagram
$$\begin{tikzcd}[row sep=tiny, column sep=tiny]
P^{s}_{R} \arrow[rr,"\pi^{s}_{R}"] \arrow[dr] \arrow[dd,"q^{s}_{R}"] & & R^{s}_{v} \arrow[dd,"q_{v}^{s}" near start] \arrow[dr] & \\
& P_{R} \arrow[rr,crossing over,"\pi_{R}" near start] & & R^{ss}_{v} \arrow[dd,"q_{v}"] \\
P^{s} \arrow[dr] \arrow[rr,"\pi^{s}" near start] & & M^{s}_{v} \arrow[dr] & \\
& P \arrow[from=uu,crossing over,"q_{R}" near start] \arrow[rr,"\pi"] & & M_{v}\\
\end{tikzcd}$$                        
The morphism $q_{R}:P_{R}\longrightarrow P$ is the map induced by $q_{v}:R^{ss}_{v}\longrightarrow M_{v}$ under base change with $\pi:P\longrightarrow M_{v}$. As $q_{v}:R^{ss}_{v}\longrightarrow M_{v}$ is a universal good quotient with respect to the action of $PGL(N_{v})$ (see Theorems 4.2.10, 4.3.3 and 4.3.4 of \cite{HL}) by definition (see Definition 4.2.2 of \cite{HL}) the action of $PGL(N_{v})$ on $R^{ss}_{v}$ extends to a compatible action on $P_{R}$, whose quotient is $q_{R}:P_{R}\longrightarrow P$.

As in the proof of point (1) we see that $P_{R}$ is locally complete intersection and its singular locus has codimesion at least 3, so it is locally factorial. As $P^{s}_{R}\subseteq P_{R}$ is smooth, the line bundle $(q_{R}^{s})^{*}L$ on $P^{s}_{R}$ then extends to a line bundle $\mathcal{O}_{P_{R}}(1)$ on $P_{R}$. As the class of $P_{R}$ in $Br(R^{ss}_{v})$ is $q_{v}^{*}\alpha$, we get $q_{v}^{*}\alpha=1$.

Notice that $\mathcal{O}_{P_{R}}(1)_{|P^{s}_{R}}=(q_{R}^{s})^{*}L$ is $PGL(N_{v})-$invariant. Moreover, since by Proposition \ref{lem:codimrsv} we have that $R^{ss}_{v}\setminus R^{s}_{v}$ has codimension at least 3 in $R^{ss}_{v}$, it follows that $P_{R}\setminus P^{s}_{R}$ has codimension at least 3 in $P_{R}$: as $P_{R}$ is locally factorial, it follows that $\mathcal{O}_{P_{R}}(1)$ is $PGL(N_{v})-$invariant. 

As $\mathcal{O}_{P_{R}}(1)_{|P^{s}_{R}}$ descends to the tautological line bundle $L$ of $P^{s}$, it follows that $\mathcal{O}_{P_{R}}(1)_{|P^{s}_{R}}$ admits a $PGL(N_{v})-$linearization. As $P_{R}\setminus P^{s}_{R}$ has codimension at least 3 in $P_{R}$ and as the conditions defining a $PGL(N_{v})-$linearization are closed, it follows that $\mathcal{O}_{P_{R}}(1)$ has a $PGL(N_{v})-$linearization as well.

We claim that $\mathcal{O}_{P_{R}}(1)$ descends to a line bundle $\mathcal{O}_{P}(1)$ on $P$. If so, then $$(q^{s}_{R})^{*}L=\mathscr{O}_{P_{R}}(1)_{|P^{s}_{R}}=(q_{R}^{*}\mathscr{O}_{P}(1))_{|P^{s}_{R}}=(q^{s}_{R})^{*}(\mathscr{O}_{P}(1)_{|P^{s}}).$$Hence $L$ and $\mathscr{O}_{P}(1)_{|P^{s}}$ are both descent bundles of the same line bundle. As the descent bundle is unique, this implies that $L\simeq\mathscr{O}_{P}(1)_{|P^{s}}$, and as $L$ is the tautological line bundle on $P^{s}$ we finally get that $\mathscr{O}_{P}(1)$ is the tautological line bundle on $P$. By Lemma 1.4 of \cite{Y5} it follows that the Brauer class $\alpha$ of $P$ is 1, completing the proof of point (2) of the statement.

We are then left with the proof of the previous claim. To do so we use Kempf's criterion (see Th\'eor\`eme 2.3 of \cite{DN}): the line bundle $\mathcal{O}_{P_{R}}(1)$ descends to a line bundle $\mathcal{O}_{P}(1)$ if and only if for every closed point $p\in P_{R}$ whose $PGL(N_{v})-$orbit is closed, the stabilizer of $p$ in $PGL(N_{v})$ acts trivially on $\mathcal{O}_{P_{R}}(1)_{p}$.

We first describe the orbits and the stabilizers of points in $P_{R}$. Consider $p\in P_{R}$, so $p=(t,x)$ where $x\in R^{ss}_{v}$, $t\in P$ and $q_{v}(x)=\pi(t)$. For every $g\in PGL(N_{v})$ we have that $g\cdot p=(t,g\cdot x)$, hence the the orbit of $p$ in $P_{R}$ is isomorphic, via $\pi_{R}$, to the orbit of $\pi_{R}(p)$ in $R^{ss}_{v}$, and the stabilizer of $p$ in $P_{R}$ is isomorphic to the stabilizer of $\pi_{R}(p)$ in $R^{ss}_{v}$.

The morphism $\pi_{R}:P_{R}\longrightarrow R^{ss}_{v}$ being projective, and hence closed, we see that if $p\in P_{R}$ is a closed point whose orbit is closed in $P_{R}$, then $\pi_{R}(p)$ is a closed point of $R^{ss}_{v}$ whose orbit is closed in $R^{ss}_{v}$. In particular, this implies that $\pi_{R}(p)=[\gamma:\mathcal{H}_{v}\longrightarrow E]\in R^{ss}_{v}$, where $E$ is an $H-$polystable sheaf (see Theorem 4.3.3 of \cite{HL}).

Let us now show that the stabilizer of a closed point $p$ with closed orbit in $P_{R}$ acts trivially on $\mathcal{O}_{P_{R}}(1)_{p}$. To do so, we mainly follow the strategy used in the proof of Theorem 5.3 of \cite{KLS}. 

Since $\pi_{R}(p)=[\gamma:\mathcal{H}_{v}\longrightarrow E]$ where $E$ is an $H-$polystable sheaf of Mukai vector $v$, and since $H$ is $v-$generic, we have $$E=\oplus_{i=1}^{s}E_{i}^{\oplus n_{i}}$$where $E_{i}$ is $H-$stable with Mukai vector $v_{i}=m_{i}w$, so that $\sum_{i=1}^{s}m_{i}n_{i}=m$. 

For every $1\leq i\leq s$ let $P_{i}$ be the Hilbert polynomial associated with $v_{i}$, $N_{i}:=P_{i}(k)$ for some $k\gg 0$, $\mathcal{H}_{i}:=\mathcal{O}_{S}(-kH)^{\oplus N_{i}}$ so that $E_{i}$ is quotient of $\mathcal{H}_{i}$ for every $i$. Moreover, up to choosing $k\gg 0$ we have $\mathcal{H}_{v}=\oplus_{i=1}^{n}\mathcal{H}^{\oplus n_{i}}_{i}$, and we let $N:=\sum_{i=1}^{s}n_{i}N_{i}$. 

Finally, let $R^{ss}_{i}$ be the open subset of $Quot(\mathcal{H}_{i},P_{H,v_{i}})$ of semistable quotients, $q_{i}:R^{ss}_{i}\longrightarrow M_{v_{i}}$ the quotient morphism, and let $$\phi:\prod_{i=1}^{s}R^{ss}_{i}\longrightarrow R^{ss}_{v},\,\,\,\,\,\,\,\,\phi([\mathcal{H}_{i}\longrightarrow F_{i}]):=[\mathcal{H}_{v}\longrightarrow\oplus_{i=1}^{s}F^{\oplus n_{i}}_{i}],$$whose image is denoted $Z$.

As $M_{v_{i}}$ is irreducible for every $i$ (by Theorem 4.4 of \cite{KLS}), the $R^{ss}_{i}$'s and $Z$ are irreducible. As $\pi_{R}:P_{R}\longrightarrow R^{ss}_{v}$ is a projective bundle, it follows that $\pi_{R}^{-1}(Z)$ is irreducible.  Moreover $Z$ contains the point $[\gamma]$ as well as a point $[\gamma':\mathcal{H}_{v}\longrightarrow E_{0}^{\oplus m}]$ for some $E_{0}\in M^{s}_{w}$. This means that we may choose $p'\in\pi_{R}^{-1}(Z)$ so that $\pi_{R}(p')=[\gamma']$. Finally, we let $$G:=\bigg(\prod_{i=1}^{s}GL(n_{i})\bigg)/\mathbb{C}^{*}\subseteq PGL(N_{v}).$$This group acts on $R^{ss}_{v}$ and stabilizes $Z$ pointwise, so it stabilizes $\pi_{R}^{-1}(Z)$ pointwise. It is exactly the stabilizer of $[\gamma]$, and hence the stabilizer of $p$, and it is contained in the stabilizer of $[\gamma']$, and hence in the stabilizer of $p'$.

Now, recall that $\mathcal{O}_{P_{R}}(1)$ has a $PGL(N_{v})-$linearization. The group $G$ acts on $\mathcal{O}_{P_{R}}(1)_{|\pi_{R}^{-1}(Z)}$ with a locally constant character: as $\pi_{R}^{-1}(Z)$ is irreducible, this character is constant. Now, the stabilizer of $p'$ is isomorphic to $PGL(m)$ (as it is isomorphic to the stabilizer of $[\gamma']$), so it has no nontrivial characters. As $G$ is contained in the stabilizer of $p'$, it acts with trivial character on $\mathcal{O}_{P_{R}}(1)_{|\pi_{R}^{-1}(Z)}$. But as $G$ is the stabilizer of $p$, and $p\in\pi_{R}^{-1}(Z)$, it follows that $G$ acts trivially on $\mathcal{O}_{P_{R}}(1)_{p}$, concluding the proof of point (2).

\textbf{Proof of (3)}. The proof is essentially the same of that of point (2): $K_{v}$ is locally factorial by Proposition A.2 of \cite{PR2}, and as $(m,k)\neq(2,1)$ we know that $K_{v}$ and $K^{s}_{v}$ are simply connected by Theorem \ref{thm:mio}, so $H^{2}(K_{v},\mathbb{Z})$ and $H^{2}(K^{s}_{v},\mathbb{Z})$ are free, and have the same rank by point (2) of Corollary \ref{cor:restriction}. The statement then follows from point (3) of Lemma \ref{lem:munrest} as soon as we know that the restriction morphism $Br(K_{v})\longrightarrow Br(K^{s}_{v})$ is injective. 

The proof of this is as in the proof of point (2), but replacing $R^{ss}_{v}$ with $R^{ss,0}_{v}$ and $R^{s}_{v}$ with $R^{s,0}_{v}$ (see the notation in section 2.1) and using that by Corollary \ref{cor:codimrsv} the variety $R^{ss,0}_{v}$ is locally factorial and the codimension of $R^{ss,0}_{v}\setminus R^{s,0}_{v}$ in $R^{ss,0}_{v}$ is at least $3$.\endproof

We conclude this section with some considerations about the relation between the cohomology groups of $M_{v}$ and $R^{ss}_{v}$, and of $M^{s}_{v}$ and $R^{s}_{v}$.

\begin{cor}
\label{cor:iniet}
Let $(S,v,H)$ be an $(m,k)-$triple. 
\begin{enumerate}
 \item The morphism $q_{v}^{s*}:H^i(M^{s}_{v},\mathbb{Z})\rightarrow H^i(R_{v}^{s},\mathbb{Z})$ is injective for $i=1,2$.
 \item The morphism $q_{v}^{*}:H^i(M_{v},\mathbb{Z})\rightarrow H^i(R_{v}^{ss},\mathbb{Z})$ is injective for $i=1,2$.
\end{enumerate}
\end{cor}

\proof For point (1), use the exact sequence induced by the Leray spectral sequence of $q_{v}^{s}$, namely 
\begin{multline}
\label{eq:leray1}
0\longrightarrow H^{1}(M_{v}^{s},\mathbb{Z})\stackrel{q_{v}^{s*}}\longrightarrow H^{1}(R^{s}_{v},\mathbb{Z})\stackrel{\psi}\longrightarrow H^{0}(M^{s}_{v},R^{1}q^{s}_{v*}\mathbb{Z})\longrightarrow \\ \longrightarrow H^{2}(M_{v}^{s},\mathbb{Z}) \stackrel{q_{v}^{s*}}\longrightarrow H^{2}(R^{s}_{v},\mathbb{Z})\stackrel{\psi}\longrightarrow H^{0}(M^{s}_{v},R^{2}q^{s}_{v*}\mathbb{Z}).
\end{multline}
As $q_{v}^{s}:R_{v}^{s}\longrightarrow M^{s}_{v}$ is a $PGL(N_{v})-$bundle (see for example Corollary 4.3.5 of \cite{HL}), the fibers of $q_{v}^{s}$ are $PGL(N_{v})$, so $H^{0}(M_{v}^{s},R^{1}q^{s}_{v*}\mathbb{Z})=0$: point (1) follows immediately.

For point (2), we use the following commutative diagram where $i=1,2$, coming from diagram (\ref{eq:commquot}): 
\begin{equation}
\label{eq:commquotnew}
\begin{CD}
H^{i}(M_{v},\mathbb{Z}) @>{i^{*}_{v}}>> H^{i}(M^{s}_{v},\mathbb{Z})\\
@V{q^{*}_{v}}VV          @VV{q^{s*}_{v}}V\\
H^{i}(R^{ss}_{v},\mathbb{Z}) @>{j^{*}_{v}}>> H^{i}(R^{s}_{v},\mathbb{Z})
\end{CD}
\end{equation}
The morphism $q_{v}^{s*}$ is injective by point (1). The morphism $i_v^*$ is an isomorphism if $i=2$ and $(m,k)\neq(2,1)$ by point (2) of Proposition \ref{prop:extmr}; it is injective if $i=2$ and $(m,k)=(2,1)$ by Lemma 3.7 of \cite{PR}, and if $i=1$ by Lemma \ref{lem:firstcohom}. The commutativity of diagram (\ref{eq:commquotnew}) implies then that $q^{*}_{v}$ is injective.\endproof

\section{The morphism $\lambda_{v}$}

The aim of this section is the construction of a morphism of $\mathbb{Z}-$modules $$\lambda_{v}:v^{\perp}\longrightarrow H^{2}(M_{v},\mathbb{Z}).$$

Let us first fix some notation that will be used all along the section. If $X$ and $Y$ are two schemes, then $p_{X}:X\times Y\longrightarrow X$ and $p_{Y}:X\times Y\longrightarrow Y$ will denote the two projections. Moreover, if $S$ is a compact complex surface and $\alpha\in H^{2*}(S,\mathbb{Q})$ is such that $\alpha=(\alpha_{0},\alpha_{1},\alpha_{2})$, we let $\alpha^{\vee}:=(\alpha_{0},-\alpha_{1},\alpha_{2})$. 

The construction of $\lambda_{v}$ is done in two steps: first, we provide a general construction for a morphism $\lambda^{s}_{v}:v^{\perp}\longrightarrow H^{2}(M^{s}_{v},\mathbb{Z})$, which is a cohomological version of the construction of line bundles on moduli spaces of sheaves by Le Potier. The second step is to extend the classes of $H^{2}(M^{s}_{v},\mathbb{Z})$ in the image of $\lambda_{v}^{s}$ to classes in $H^{2}(M_{v},\mathbb{Z})$.

The section is divided into two main parts: first we provide the actual construction of $\lambda_{v}$; finally, we will relate the morphism $\lambda_{v}$ with the morphism $\lambda_{v'}$ for $v=mw$ and $v'=pw$, where $1\leq p\leq m$.

\subsection{The construction of the morphism $\lambda_{v}$.}

Let $S$ be a projective K3 surface or an Abelian surface and $Y$ any complex algebraic variety. We will associate with every $Y-$flat family $\mathscr{F}$ of coherent sheaves on $S\times Y$ a morphism $\mu_{\mathscr{F}}: H^{2*}(S,\mathbb{Z})\rightarrow H^{2}(Y,\mathbb{Z})$ that factors through the topological $K-$theory of $Y$.

If $Z$ is a complex algebraic variety, we let $K^{0}_{top}(Z)$ and $K^{0}_{alg}(Z)$ be the topological and the algebraic Grothendieck groups of topological and algebraic complex vector bundles on $Z$. If $X$ and $Y$ are complex, quasi-projective algebraic varieties and $f:X\rightarrow Y$ is a proper, locally complete intersection morphism, Baum, Fulton and MacPherson constructed in \cite{BFMa} a push-forward in topological $K-$theory $$f_{!,top}:K^{0}_{top}(X)\rightarrow K^{0}_{top}(Y)$$(see the Definition on page 141 therein). 

In \cite{BFM} they also proved (see section 4.2 therein) that $f_{!,top}$ extends the push forward in algebraic $K-$theory $f_{!,alg}:K^{0}_{alg}(X)\rightarrow K^{0}_{alg}(Y)$ (whose definition is given at page 170 of \cite{BFM} or at page 36 of \cite{HL}): more precisely if we let $\alpha_{X}:K^{0}_{alg}(X)\longrightarrow K^{0}_{top}(X)$ be the natural morphism mapping the class of an algebraic complex vector bundle on $X$ to the class of the underlying topological vector bundle, then the following diagram commutes
$$\begin{CD}
K^{0}_{alg}(X) @>{f_{!,alg}}>> K^{0}_{alg}(Y)\\
@V{\alpha_{X}}VV          @VV{\alpha_{Y}}V\\
K^{0}_{top}(X) @>{f_{!,top}}>> K^{0}_{top}(Y)
\end{CD}$$
 
\begin{oss}
\label{oss:commfshrik}
{\rm As a consequence of the definition of the push-forward morphism in topological $K-$theory one may show that if $X$ is a complex projective variety and $f:Y'\longrightarrow Y$ is a morphism of complex algebraic varieties, letting $p:X\times Y\longrightarrow Y$ and $p':X\times Y'\longrightarrow Y'$ be the two projections, then the following diagram is commutative:
\begin{equation}
\label{eq:fppf}
\begin{tikzcd}
K^{0}_{top}(X\times Y) \arrow[d,swap,"p_{!,top}"] \arrow[r,"(id_{X}\times f)^{*}_{top}"] & K^{0}_{top}(X\times Y') \arrow[d,"p'_{!,top}"] \\
K^{0}_{top}(Y) \arrow[r,swap,"f^{*}_{top}"] & K^{0}_{top}(Y')
\end{tikzcd}
\end{equation}}
\end{oss}

We will define the morphism $\mu_{\mathscr{F}}$ by means of the push-forward in topological $K-$theory. Since the surface $S$ is smooth and $\mathscr{F}$ is a coherent sheaf on $S\times Y$ which is flat over $Y$, by Proposition 2.1.10 of \cite{HL} the sheaf $\mathscr{F}$ admits a finite locally free resolution. Hence it defines a class $[\mathscr{F}]_{alg}$ in $K^{0}_{alg}(S\times Y)$ and a class $[\mathscr{F}]_{top}$ in $K^{0}_{top}(S\times Y)$, so that $\alpha_{S\times Y}([\mathscr{F}]_{alg})=[\mathscr{F}]_{top}$.

Since for K3 surfaces and Abelian surfaces the Mukai vector gives an isomorphism between $K^{0}_{top}(S)$ and $H^{2*}(S,\mathbb{Z})$, for every $\alpha\in H^{2*}(S,\mathbb{Z})$ there is a class $E_{\alpha}\in K^{0}_{top}(S)$ such that $v(E_{\alpha})=\alpha$. We moreover let $E_{\alpha}^{\vee}$ be the dual of $E_{\alpha}$, so that $v(E_{\alpha}^{\vee})=\alpha^{\vee}$.

This allows us to give the following definition, for which we need the following notation: write $\otimes_{top}$ for the tensor product in $K^{0}_{top}(S\times Y)$, and if $\gamma\in K^{0}_{top}(Y)$ we let $c_1(\gamma)\in H^{2}(Y,\mathbb{Z})$ be the first Chern class of $\gamma$.

\begin{defn}
If $\mathscr{F}$ is a $Y-$flat family of coherent sheaves on $S\times Y$, we define the group morphism $\mu_{\mathscr{F}}:H^{2*}(S,\mathbb{Z})\rightarrow H^{2}(Y,\mathbb{Z})$ as follows: for every $\alpha\in H^{2*}(S,\mathbb{Z})$ we let $$\mu_{\mathscr{F}}(\alpha):=c_{1}(p_{Y!,top}(p_{S,top}^{*}(E^{\vee}_{\alpha})\otimes_{top} [\mathscr{F}])).$$We write $\mu^{\mathbb{Q}}_{\mathscr{F}}:H^{2*}(S,\mathbb{Q})\longrightarrow H^{2}(Y,\mathbb{Q})$ for the $\mathbb{Q}-$linear morphism induced by $\mu_{\mathscr{F}}$ by tensorization with $\mathbb{Q}$.
\end{defn}

The morphism $\mu_{\mathscr{F}}$ we just defined coincides, when calculated on algebraic classes, with the one defined in Chapter 8 of \cite{HL} (see Definition 8.1.1 therein). The following collects the general properties of $\mu_{\mathscr{F}}$, and are similar to those verified by the Le Potier morphism (see Lemma 8.1.2 of \cite{HL}).

\begin{lem}
\label{lem:propmu}Let $Y$ be a complex algebraic variety and $\mathscr{F}$ a $Y-$flat family on $S\times Y$ of coherent sheaves on $S$ of Mukai vector $v$.
\begin{enumerate}
 \item If $0\longrightarrow\mathscr{F}'\longrightarrow\mathscr{F}\longrightarrow\mathscr{F}''\longrightarrow 0$ is an exact sequence of $Y-$flat coherent sheaves on $S\times Y$, then $\mu_{\mathscr{F}}=\mu_{\mathscr{F}'}+\mu_{\mathscr{F}''}$.
 \item If $f:Y'\longrightarrow Y$ is a morphism of complex algebraic varieties and $\mathscr{F}$ is a $Y-$flat coherent sheaf on $S\times Y$, then $\mu_{(id_{S}\times f)^{*}\mathscr{F}}=f^{*}\circ\mu_{\mathscr{F}}$.
 \item If  $V$ is a locally free sheaf of rank $r$ on $Y$, then for every class $\alpha\in H^{2*}(S,\mathbb{Z})$ we have $$\mu_{\mathscr{F}\otimes p_{Y}^{*}V}(\alpha)=r\mu_{\mathscr{F}}(\alpha)+(\alpha,v)c_{1}(V)$$where $(\alpha,v)$ is the Mukai pairing between $\alpha$ and $v$.
\item For every class $\alpha\in H^{2*}(S,\mathbb{Q})$ we have $$\mu^{\mathbb{Q}}_{\mathscr{F}}(\alpha)=[p_{Y*}(p_{X}^{*}(\alpha^{\vee}\cdot\sqrt{td(S)})\cdot ch(\mathscr{F}))]_{2},$$ where for a class $\beta\in H^{*}(Y,\mathbb{Q})$ we let $[\beta]_{2}$ be its component in $H^{2}(Y,\mathbb{Q})$.
\end{enumerate}
\end{lem}

\proof For point (1), since $p_{Y!,top}: K^{0}_{top}(S\times Y)\longrightarrow K^{0}_{top}(Y)$ is well defined and preserves exact sequences, the additivity of the first Chern class on exact sequences of complex vector bundles implies the result. 

For point (2), let $p'_{S}:S\times Y'\longrightarrow Y'$ be the projection, so that $p'_{S}=p_{S}\circ (id_{S}\times f)$. We then have $$\mu_{(id_{S}\times f)^{*}\mathscr{F}}(\alpha)=c_{1}(p_{Y'!top}((p'_{S})_{top}^{*}[E^{\vee}_{\alpha}]\otimes_{top}[(id_{S}\times f)^{*}\mathscr{F}]))=$$ $$=c_{1}(p_{Y'!top}((id_{S}\times f)_{top}^{*}(p_{S,top}^{*}[E^{\vee}_{\alpha}])\otimes_{top}(id_{S}\times f)_{top}^{*}[\mathscr{F}]))=$$ $$=c_{1}(p_{Y'!top}((id_{S}\times f)_{top}^{*}(p_{S,top}^{*}[E^{\vee}_{\alpha}]\otimes_{top}[\mathscr{F}]))).$$By Remark \ref{oss:commfshrik}, and as $c_{1}\circ f^{*}_{top}=f^{*}\circ c_{1}$, this last is equal to $$c_{1}(f_{top}^{*}(p_{Y!top}(q_{S,top}^{*}[E^{\vee}_{\alpha}]\otimes_{top}[\mathscr{F}])))=f^{*}(\mu_{\mathscr{F}}(\alpha)).$$

For point (3) we just need to use the projection formula in topological $K-$theory (see Proposition 4.4 (b) of \cite{BFMa}), giving $$p_{Y!,top}(p_{S,top}^{*}(E^{\vee}_{\alpha})\otimes_{top} [\mathscr{F}\otimes p_{Y}^{*}V])=p_{Y!,top}(p_{S,top}^{*}(E^{\vee}_{\alpha})\otimes_{top} [\mathscr{F}])\otimes_{top}^Y[V]$$where $\otimes_{top}^Y$ denotes the tensor product on $K^{0}_{top}(Y)$.

As the rank of $p_{Y!,top}(p_{S,top}^{*}(E^{\vee}_{\alpha})\otimes_{top} [\mathscr{F}])$ equals the Mukai pairing $(\alpha,v)$,
using the formula for the first Chern class of a tensor product we obtain $$\mu_{\mathscr{F}\otimes p_{Y}^{*}V}(\alpha)=c_1(p_{Y!,top}(p_{S,top}^{*}(E^{\vee}_{\alpha})\otimes_{top} [\mathscr{F}])\otimes_{top}^Y[V])=$$ $$=rc_1(p_{Y!,top}(p_{S,top}^{*}(E^{\vee}_{\alpha})\otimes_{top} [\mathscr{F}])+(\alpha,v)c_1(V)=$$ $$=r\mu_{\mathscr{F}}(\alpha)+(\alpha,v)c_1(V).$$ 

Point (5) is a direct consequence of Grothendieck-Riemann-Roch for topological $K-$theory, which is the commutativity of diagram (1) in Remark (2) of Section 5 in \cite{BFM}.\endproof

We now use the morphism $\mu_{\mathscr{F}}$ in two different situations. One is when $Y=M^{s}_{v}$ and $\mathscr{F}$ is a quasi-universal family on $S\times M^{s}_{v}$ (see Definition 4.6.1 of \cite{HL}). The other is when $Y=R^{s}_{v}$ and $\mathscr{F}=\mathcal{Q}^{s}_{v}$ is the restriction to $S\times R^{s}_{v}$ of the tautological quotient $Q_{v}$ on $S\times Quot(\mathcal{H}_{v},P_{v,H})$ (so the quotient map $q_{v}^{s}:R_{v}^{s}\longrightarrow M^{s}_{v}$ is the modular map associated with $\mathcal{Q}_{v}^{s}$). 

The following shows that these two constructions are strictly related:

\begin{prop}
\label{prop:descent}
Let $S$ be a projective K3 surface or an Abelian surface, $v$ a Mukai vector and $H$ a $v-$generic polarization on it.
\begin{enumerate}
 \item There is a unique morphism $$\lambda^{s}_{v}: v^{\perp}\rightarrow H^2(M^{s}_{v},\mathbb{Z})$$such that $q_{v}^{s*}\circ \lambda^{s}_{v}=\mu_{\mathcal{Q}_{v}^{s}}$. 
 \item If $\mathscr{F}$ is a quasi-universal family on $S\times M^{s}_{v}$ of similitude $\rho$, then $\rho\lambda^{s}_{v}(\alpha)=\mu_{\mathscr{F}}(\alpha)$ for every $\alpha\in v^{\perp}$.
 \item Let $\lambda_{v}^{s,\mathbb{Q}}:v^{\perp}\otimes\mathbb{Q}\longrightarrow H^{2}(M^{s}_{v},\mathbb{Q})$ be the morphism induced by $\lambda_{v}^{s}$ by tensorization with $\mathbb{Q}$. If $\mathscr{G}$ is a flat family of stable sheaves on $S$ with Mukai vector $v$ parametrized by $Y$ and $f_\mathscr{G}:Y\rightarrow M^{s}_{v}$ is the associated modular map, then for all $\alpha\in v^{\perp}\otimes\mathbb{Q}$ we have $$f_\mathscr{G}^{*}(\lambda^{s,\mathbb{Q}}_v(\alpha))=\mu^{\mathbb{Q}}_{\mathscr{G}}(\alpha)$$as elements of $H^{2}(Y,\mathbb{Q})$.
\end{enumerate}
\end{prop}

\proof By point (1) of Corollary \ref{cor:iniet}, point (1) follows if we show that $\mu_{\mathcal{Q}_{v}^{s}}(\alpha)$ descends to a class $\lambda_{v}^{s}(\alpha)\in H^{2}(M_{v}^{s},\mathbb{Z})$. 

To prove this, consider the exact sequence induced by the Leray spectral sequence of $q_{v}^{s}$, which is $$0\longrightarrow H^{2}(M^{s}_{v},\mathbb{Z})\longrightarrow H^{2}(R^{s}_{v},\mathbb{Z})\longrightarrow H^{0}(M^{s}_{v},R^{2}q^{s}_{v*}\mathbb{Z})$$(see Corollary \ref{cor:iniet} and the exact sequence (\ref{eq:leray1}) in its proof). 

If for every $p\in M^{s}_{v}$ we let $i_{p}:(q_{v}^{s})^{-1}(p)\longrightarrow R_{v}^{s}$ be the closed embedding, to prove that $\mu_{\mathcal{Q}_{v}^{s}}(\alpha)$ descends to a class $\lambda_{v}^{s}(\alpha)\in H^{2}(M_{v}^{s},\mathbb{Z})$ it is enough to show that $i_{p}^{*}(\mu_{\mathcal{Q}_{v}^{s}}(\alpha))=0$ for every $p\in M^{s}_{v}$.

Since $(q_{v}^{s})^{-1}(p)\simeq PGL(N_{v})$, if $p$ corresponds to the isomorphism class of an $H-$stable sheaf $F$ then the restriction $(id_{S}\times i_{p})^{*}\mathcal{Q}_{v}^{s}$ of the tautological quotient family $\mathcal{Q}_{v}^{s}$ to $S\times(q_{v}^{s})^{-1}(p)$ is of the form $p_{1}^{*}F\otimes p_{2}^{*}L$, where we let $p_{1}:S\times PGL(N_{v})\longrightarrow S$ and $p_{2}:S\times PGL(N_{v})\longrightarrow PGL(N_{v})$ be the two projections and $L$ a holomorphic line bundle on $PGL(N_{v})$.

By points (2) and (3) of Lemma \ref{lem:propmu}, and since $\alpha \in v^{\perp}$, we have$$(i_{p}^{*}\circ \mu_{\mathcal{Q}_{v}^{s}})(\alpha)=\mu_{(id_{S}\times i_{p})^{*}\mathcal{Q}_{v}^{s}}(\alpha)=\mu_{p_{1}^{*}F\otimes p_{2}^{*}L}(\alpha)=$$ $$=\mu_{p_{1}^{*}F}(\alpha)+(\alpha,v)c_1(L)=\mu_{p_{1}^{*}F}(\alpha).$$

As $PGL(N_{v})=(q^{s})^{-1}(p)$, we may view $p_{1}:S\times PGL(N_{v})\longrightarrow S$ as the morphism $id_{S}\times q_{p}^{s}:S\times PGL(N_{v})\longrightarrow S\times\{p\}$ where $q_{p}^{s}:PGL(N_{v})\longrightarrow\{p\}$ is the projection. By point (2) of Lemma \ref{lem:propmu} we get $$\mu_{p_{1}^{*}F}(\alpha)=\mu_{(id_{S}\times q_{p}^{s})^{*}F}(\alpha)=q_{p}^{s*}(\mu_{F}(\alpha))=0.$$As a consequence $i_{p}^{*}(\mu_{\mathcal{Q}_{v}^{s}}(\alpha))=0$, completing the proof of point (1).

For point (2) we need to compare the tautological quotient $\mathcal{Q}_{v}^{s}$ and the pull back $(id_S\times q_{v}^{s})^{*}(\mathscr{F})$ on $S\times R_{v}^{s}$. 

By definition, if $F$ is an $H-$stable sheaf giving a point $a\in M_{v}^{s}$ and if $b\in(q_{v}^{s})^{-1}(a)$, then the restriction of $\mathcal{Q}_{v}^{s}$ to the fiber $pr_{R}^{-1}(b)$ of the projection $pr_{R}:S\times R^{s}_{v}\longrightarrow R^{s}_{v}$ is isomorphic to $F$. Moreover, the restriction to $pr_{R}^{-1}(b)$ of $(id_S\times q_{v}^{s})^{*}(\mathscr{F})$ is isomorphic to $F^{\oplus\rho}$. 

As $R_{v}^{s}$ is reduced, the sheaf $pr_{R*}\mathcal{H}om(\mathcal{Q}_{v}^{s},(id_S\times q^{s}_{v})^{*}(\mathscr{F}))$ is isomorphic to a rank $\rho$ vector bundle $\mathcal{V}$, and the evaluation map gives an isomorphism $\mathcal{Q}_{v}^{s}\otimes pr_R^{*}\mathcal{V}\simeq (id_S\times q_{v}^{s})^{*}(\mathscr{F})$. By point (2) of Lemma \ref{lem:propmu} we then have 
\begin{equation}
\label{eq:muqv1}
\mu_{\mathcal{Q}_{v}^{s}\otimes pr_R^{*}\mathcal{V}}(\alpha)=\mu_{(id_S\times q_{v}^{s})^{*}(\mathscr{F})}(\alpha)=(q_{v}^{s*}\circ \mu_{\mathscr{F}})(\alpha)
\end{equation}
for every $\alpha \in v^{\perp}$. 

On the other hand, by point (3) of Lemma \ref{lem:propmu} and point (1) of this Proposition we just proved, for every $\alpha\in v^{\perp}$ we have 
\begin{equation}
\label{eq:muqv2}
\mu_{\mathcal{Q}_{v}^{s}\otimes pr_R^{*}\mathcal{V}}(\alpha)=\rho\mu_{\mathcal{Q}_{v}^{s}}(\alpha)+(\alpha,v)c_1(\mathcal{V})=\rho\mu_{\mathcal{Q}_{v}^{s}}(\alpha)=\rho (q_{v}^{s*}\circ\lambda^{s}_{v})(\alpha).
\end{equation}
Comparing the last terms of equations (\ref{eq:muqv1}) and (\ref{eq:muqv2}) and using the injectivity of $q_{v}^{s*}$ we get $\mu_{\mathscr{F}}(\alpha)=\rho\lambda^{s}_{v}(\alpha)$, so point (2) is proved.
 
We are now left with the proof of point (3): we have 
\begin{equation}
\label{eq:lambdamu1}
\mu^{\mathbb{Q}}_{(id_{S}\times f_{\mathscr{G}})^{*}\mathscr{F}}(\alpha)=f_{\mathscr{G}}^{*}(\mu^{\mathbb{Q}}_{\mathscr{F}}(\alpha))=\rho f^{*}_{\mathscr{G}}(\lambda^{s,\mathbb{Q}}_{v}(\alpha))
\end{equation}
where the first equality comes from the fact that the formula in point (4) of Lemma \ref{lem:propmu} is compatible with base change, and the second equality is the previous part of the proof.

But as $\mathscr{F}$ is semiuniversal of similitude $\rho$ there exists a rank $\rho$ vector bundle $V$ on $Y$ such that 
$(id_S\times f_\mathscr{G})^* \mathscr{F}=\mathscr{G}\otimes p_Y^*(V)$. By point (2) of Lemma \ref{lem:propmu} we then get 
\begin{equation}
\label{eq:lambdamu2}
\mu^{\mathbb{Q}}_{(id_S\times f_\mathscr{G})^* \mathscr{F}}(\alpha)=\mu^{\mathbb{Q}}_{\mathscr{G}\otimes p_Y^*(V)}(\alpha)=\rho \mu^{\mathbb{Q}}_{\mathscr{G}}(\alpha).
\end{equation}
Comparing equations (\ref{eq:lambdamu1}) and (\ref{eq:lambdamu2}) we are done.\endproof

\begin{oss}
\label{oss:muvlv}
{\rm If $(S,v,H)$ is a $(1,k)-$ or $(2,1)-$triple, in \cite{OG1} and \cite{Y1} (for $(1,k)-$triples) and in \cite{PR} (for $(2,1)-$triples) it is defined a morphism $\mu^{s}_{v}:v^{\perp}\longrightarrow H^{2}(M^{s}_{v},\mathbb{Q})$ as follows: if $\mathscr{F}$ is a quasi-universal family of similitude $\rho$ on $S\times M_{v}^{s}$, for every $\alpha\in v^{\perp}$ we let $$\mu^{s}_{v}(\alpha)=\frac{1}{\rho}[p_{M^{s}_{v}*}(p_{S}^{*}(\alpha^{\vee}\cdot\sqrt{td(S)})\cdot ch(\mathscr{F}))]_{2}.$$By point (4) of Lemma \ref{lem:propmu} and point (2) of Proposition \ref{prop:descent} we then get $$\rho\mu^{s}_{v}(\alpha)=\mu_{\mathscr{F}}^{\mathbb{Q}}(\alpha)=\rho \lambda^{s,\mathbb{Q}}_{v}(\alpha),$$so that $\mu^{s}_{v}=\lambda^{s,\mathbb{Q}}_{v}$.}
\end{oss}

We may now give the main definition of this section:

\begin{defn} 
Let $(S,v,H)$ be an $(m,k)-$triple.
\begin{enumerate}
 \item We let $\lambda_{v}:v^{\perp}\longrightarrow H^{2}(M_{v},\mathbb{Z})$ be the unique group morphism such that
$i_v^{*}\circ \lambda_{v}=\lambda^{s}_{v}$.
 \item If $S$ is an Abelian surface, let $\lambda_{v}^{s,0}:v^{\perp}\longrightarrow H^{2}(K_{v}^{s},\mathbb{Z})$ denote the composition of $\lambda^{s}_{v}$ with the pull-back by the inclusion of $K_{v}^{s}$ in $M_{v}^{s}$. We let $\lambda^{0}_{v}:v^{\perp}\longrightarrow H^{2}(K_{v},\mathbb{Z})$ be the unique group morphism such that $i_v^{0*}\circ \lambda^{0}_{v}=\lambda^{s,0}_{v}$.
\end{enumerate}
\end{defn}

\begin{oss}
\label{oss:bendef}
{\rm We notice first of all that the uniqueness of $\lambda_{v}$ is a consequence of the injectivity of $i_{v}^{*}$ (which is point (2) of Proposition \ref{prop:extmr} for $(m,k)\neq(2,1)$, and it is Lemma 3.7 of \cite{PR} for $(m,k)=(2,1)$. Similarly, the uniqueness of $\lambda^{0}_{v}$ is a consequence of the injectivity of $i^{0*}_{v}$ (which is point (3) of Proposition \ref{prop:extmr} for $(m,k)\neq(2,1)$, and it is Lemma 3.7 of \cite{PR} for $(m,k)=(2,1)$).}

{\rm The existence of the morphism $\lambda_{v}$ (resp. $\lambda^{0}_{v}$) is as follows.
\begin{itemize}
 \item If $m=1$, this is immediate as $M_{v}=M^{s}_{v}$ (resp. $K_{v}=K^{s}_{v}$).
 \item If $(m,k)=(2,1)$, by Remark \ref{oss:muvlv} the class $\mu^{s}_{v}(\alpha)\in H^{2}(M^{s}_{v},\mathbb{Q})$ is the rational class induced by $\lambda^{s}_{v}(\alpha)$. The extension $\mu_{v}(\alpha)\in H^{2}(M_{v},\mathbb{Z})$ of $\mu_{v}^{s}(\alpha)$ given by Proposition 3.9 of \cite{PR} is then the unique extension of $\lambda^{s}_{v}(\alpha)$ to $M_{v}$.
 \item If $(m,k)\neq(2,1)$ and $m\neq 1$, this is a consequence of the fact that $i_{v}^{*}:H^{2}(M_{v},\mathbb{Z})\longrightarrow H^{2}(M^{s}_{v},\mathbb{Z})$ is an isomorphism by point (2) of Proposition \ref{prop:extmr} (resp. of the fact that $i_{v}^{0*}:H^{2}(K_{v},\mathbb{Z})\longrightarrow H^{2}(K^{s}_{v},\mathbb{Z})$ is an isomorphism by point (3) of Proposition \ref{prop:extmr}).
\end{itemize}}
\end{oss}

As a formal consequence of the definition we have the following:

\begin{prop}
\label{prop:lambdav}
Let $(S,v,H)$ be an $(m,k)-$triple and equip $v^{\perp}$ with the weight-two Hodge structure induced by the Hodge structure of the Mukai lattice of $S$.  
\begin{enumerate}
 \item The morphism $\lambda_{v}:v^{\perp}\longrightarrow H^{2}(M_{v},\mathbb{Z})$ is a morphism of pure weight-two Hodge structures.
 \item If $S$ is Abelian, the morphism $\lambda^{0}_{v}:v^{\perp}\longrightarrow H^{2}(K_{v},\mathbb{Z})$ is a morphism of pure weight-two Hodge structures. 
\end{enumerate}
\end{prop}

\proof If $m=1$ this is due to Mukai, while if $(m,k)=(2,1)$, this is proved in \cite{PR} (see Remark \ref{oss:muvlv}). We are left with all other cases. 

First, recall by Corollary \ref{cor:restriction} that the weight two Hodge structures of $H^{2}(M_{v},\mathbb{Z})\simeq H^{2}(M^{s}_{v},\mathbb{Z})$ and $H^{2}(K_{v},\mathbb{Z})\simeq H^{2}(K^{s}_{v},\mathbb{Z})$ are pure. Moreover by point (3) of Lemma \ref{lem:propmu} the morphism $\lambda_{v}^{s}$ equals the restriction of $\mu_{\mathscr{F}}/\rho$ to $v^{\perp}$ for a quasi-universal family $\mathscr{F}$. 

Since by point (3) of Proposition \ref{prop:descent} the morphism $\mu_{\mathscr{F}}$ may be defined by means of the algebraic cycle $ch(\mathscr{F})$, it follows that $\lambda_{v}^{s}:v^{\perp}\rightarrow H^{2}(M^{s}_{v},\mathbb{Z})$ is a Hodge morphism. Since $i_v^{*}$ is a Hodge isomorphism, point (1) holds.

Point (2) follows from point (1) since the closed embedding $K_{v}^{s}\subset M_{v}^{s}$ induces an Hodge morphism and $i_v^{0*}$ is a Hodge isomorphism.\endproof

We conclude this section with two useful characterizations of the morphism $\lambda_{v}$: the first holds if $(m,k)\neq(2,1)$, and may be viewed as parallel to the definition of $\lambda_{v}^{s}$ (see point (1) of Proposition \ref{prop:descent}).

\begin{oss}\bf{Characterization of $\lambda_{v }$ for $(m,k)\ne (2,1)$}.
\label{oss:carattlambda}
{\rm Let $(S,v,H)$ be an $(m,k)$-triple with $(m,k)\ne (2,1)$ and $\mathcal{Q}_{v}$ the restriction to $S\times R^{ss}_{v}$ of the universal family of $S\times Quot(\mathcal{H},v)$. If $\alpha\in v^{\perp}$, the class $\lambda_{v}(\alpha)$ may be characterized as the unique class in $H^{2}(M_{v},\mathbb{Z})$ such that $$q_{v}^{*}\lambda_{v}(\alpha)=\mu_{\mathcal{Q}_{v}}(\alpha).$$Uniqueness is a consequence of the injectivity of $q_{v}^{*}$ (see point (2) of Corollary \ref{cor:iniet}). For $q_{v}^{*}(\lambda_{v}(\alpha))=\mu_{\mathcal{Q}_{v}}(\alpha)$, the commutativity of diagram (\ref{eq:commquotnew}) and the definitions of $\lambda_{v}$ and $\lambda_{v}^{s}$ give $$j_{v}^{*}(q_{v}^{*}(\lambda_{v}(\alpha)))=q_{v}^{s*}(i_{v}^{*}(\lambda_{v}(\alpha)))=q_{v}^{s*}(\lambda_{v}^{s}(\alpha))=\mu_{\mathcal{Q}_{v}^{s}}(\alpha),$$where $\mathcal{Q}_{v}^{s}=(id_{S}\times j_{v})^{*}\mathcal{Q}_{v}$ is the universal quotient on $R^{s}_{v}$. But since $$j_{v}^{*}(\mu_{\mathcal{Q}_{v}}(\alpha))=\mu_{(id_{S}\times j_{v})^{*}\mathcal{Q}_{v}}(\alpha)=\mu_{\mathcal{Q}_{v}^{s}}(\alpha)$$by point (2) of Lemma \ref{lem:propmu}, we see that $j_{v}^{*}q_{v}^{*}(\lambda_{v}(\alpha))=j_{v}^{*}(\mu_{\mathcal{Q}_{v}}(\alpha))$. Now, as $(m,k)\neq(2,1)$, by point (1) of Proposition \ref{prop:extmr} we have that $j_{v}^{*}$ is injective, hence $q_{v}^{*}(\lambda_{v}(\alpha))=\mu_{\mathcal{Q}_{v}}(\alpha)$.}
\end{oss}

The characterization of $\lambda_{v}(\alpha)$ in terms of its image in $H^{2}(R^{ss}_{v},\mathbb{Z})$ presented in Remark \ref{oss:carattlambda} works for $(m,k)\neq(2,1)$ because of the injectivity of $j_{v}^{*}$. As this does not hold for $(m,k)=(2,1)$, the proof of the same characterization does not work for $(2,1)-$triples. However there is another characterization of $\lambda_{v}(\alpha)$ in terms of its image in $H^2(R^s_v,\mathbb{Q})$ holding for all $(m,k)-$triples that we present in the following:

\begin{oss}\bf{Characterization of $\lambda_{v}$ for any $(m,k)$}.
\label{oss:carattlambdaII} 
{\rm Let us denote $j_{v,\mathbb{Q}}^{*}:H^{2}(R^{ss}_{v},\mathbb{Z})\longrightarrow H^{2}(R^{s}_{v},\mathbb{Q})$ and $q_{v,\mathbb{Q}}^{s*}:H^{2}(M^{s}_{v},\mathbb{Z})\longrightarrow H^{2}(R^{s}_{v},\mathbb{Q})$ the compositions of $j_{v}^{*}$ and $q_{v}^{s*}$, respectively, with the natural morphism $H^{2}(R^{s}_{v},\mathbb{Z})\longrightarrow H^{2}(R^{s}_{v},\mathbb{Q})$ induced by tensorization with $\mathbb{Q}$. If $\alpha\in v^{\perp}$, then $\lambda_{v}(\alpha)\in H^{2}(M_{v},\mathbb{Z})$ is the unique class such that 
\begin{equation}
\label{eq:carlambda} 
j_{v,\mathbb{Q}}^{*}(q_{v}^{*}(\lambda_{v}(\alpha)))=q_{v,\mathbb{Q}}^{s*}(i_v^{*}(\lambda_{v}(\alpha)))=\mu_{\mathcal{Q}_{v}^{s}}^{\mathbb{Q}}(\alpha).
\end{equation}
The first equality follows from commutativity of diagram (\ref{eq:commquotnew}). For the second equality recall that by definition of $\lambda_{v}$ we have $i_v^{*}(\lambda_{v}(\alpha))=\lambda_{v}^{s}(\alpha)$, so $q^{s*}_{v,\mathbb{Q}}(i_{v}^{*}(\lambda_{v}(\alpha)))=q^{s*}_{v,\mathbb{Q}}(\lambda^{s}_{v}(\alpha))$, and by definition of $\lambda_{v}^{s}$ (see point (1) of Proposition \ref{prop:descent}) we have $q^{s*}_{v,\mathbb{Q}}(\lambda^{s}_{v}(\alpha))=\mu_{\mathcal{Q}_{v}^{s}}^{\mathbb{Q}}(\alpha)$.}

{\rm The fact that $\lambda_{v}(\alpha)$ is the unique class verifying equation (\ref{eq:carlambda}) follows if we know that $q^{s*}_{v,\mathbb{Q}}\circ i_{v}^{*}$ is injective. But this follows since $i_{v}^{*}$ is injective (by point (2) of Proposition \ref{prop:extmr} for $(m,k)\neq(2,1)$, and by Lemma 3.7 of \cite{PR} for $(m,k)=(2,1)$), and $q^{s*}_{v,\mathbb{Q}}$ is injective (since $q_{v}^{s*}$ is injective by point (1) of Corollary \ref{cor:iniet}).}
\end{oss}

\subsection{Relation between $\lambda_{v}$ and $\lambda_{w}$.}

If $(S,v,H)$ is an $(m,k)-$triple where $S$ is a K3 surface or an Abelian surface, write $v=mw$ where $w$ is a primitive Mukai vector. We will discuss the relation between $\lambda_{v}$ and $\lambda_{pw}$ for $1\leq p<m$ by constructing a modular morphism $$f_{p}:M_{pw}\times M_{(m-p)w}\longrightarrow M_{v}$$mapping a pair of S-equivalence classes $([F],[G])\in M_{pw}\times M_{(m-p)w}$ to the S-equivalence class $[F\oplus G]\in M_{v}$ and, for any $H-$semistable sheaf $G$ with Mukai vector $(m-p)w$, a modular morphism $$f_{p,[G]}:M_{pw}\longrightarrow M_{v}$$ obtained by restriction of $f_{p}$ to $M_{pw}\times\{[G]\}$.

To construct $f_{p}$, recall that for a Mukai vector $u$ we have a quotient morphism $q_{u}:R^{ss}_{u}\longrightarrow M_{u}$, which exhibits $M_{u}$ as a universal good quotient under the action of $PGL(N_{u})$. It follows that the morphism $$q_{pw}\times q_{(m-p)w}:R^{ss}_{pw}\times R^{ss}_{(m-p)w}\longrightarrow M_{pw}\times M_{(m-p)w}$$gives a universal good quotient for the action of $PGL(N_{pw})\times PGL(N_{(m-p)w})$ on $R^{ss}_{pw}\times R^{ss}_{(m-p)w}$.

We first show the existence of morphism $f_{p}^{R}:R^{ss}_{pw}\times R^{ss}_{(m-p)w}\rightarrow R^{ss}_{v}$ by using the existence of a universal family for the Quot functor. 

Let $\mathcal{Q}_{pw}$ and $\mathcal{Q}_{(m-p)w}$ be the universal quotients families on $S\times R^{ss}_{pw}$ and $S\times R^{ss}_{(m-p)w}$ respectively. We moreover let $\pi_1 ^R:R^{ss}_{pw}\times R^{ss}_{(m-p)w}\rightarrow R^{ss}_{pw}$ and $\pi_2 ^R:R^{ss}_{pw}\times R^{ss}_{(m-p)w}\rightarrow R^{ss}_{(m-p)w}$ be the two projections, and set $$\overline{\mathcal{Q}}_{pw}:=(id_S\times \pi_1 ^R)^*\mathcal{Q}_{pw},\,\,\,\,\,\,\,\overline{\mathcal{Q}}_{(m-p)w}:=(id_S\times \pi_2^R)^*\mathcal{Q}_{(m-p)w}.$$

Since $P_{pw,H}+P_{(m-p)w,H}=P_{v,H}$, the coherent sheaf $\overline{\mathcal{Q}}_{pw}\oplus\overline{\mathcal{Q}}_{(m-p)w}$ is a family on $S\times R^{ss}_{pw}\times R^{ss}_{(m-p)w}$ of $H$-semistable quotients of $\mathcal{H}_{v}$ with Mukai vector $v$. If $\mathcal{Q}_{v}$ is a universal family on $S\times R^{ss}_{v}$, there exists then a modular morphism $f_{p}^{R}:R^{ss}_{pw}\times R^{ss}_{(m-p)w}\rightarrow R^{ss}_{v}$ such that 
\begin{equation}
\label{eq:mancava}
(id_S\times f_{p}^{R})^{*}\mathcal{Q}_v\simeq\overline{\mathcal{Q}}_{pw}\oplus\overline{\mathcal{Q}}_{(m-p)w}.
\end{equation}
    
Now, remark that quotients in $R^{ss}_{pw}\times R^{ss}_{(m-p)w}$ belonging to the same $PGL(N_{pw})\times PGL(N_{(m-p)w})$-orbit
give isomorphic sheaves. It follows that the composition $q_v \circ f_{p}^{R}$ is $PGL(N_{pw})\times PGL(N_{(m-p)w})-$invariant 
and descends to the universal good quotient $M_{pw}\times M_{(m-p)w}$, i. e. there is a morphism $f_{p}:M_{pw}\times M_{(m-p)w}\longrightarrow M_{v}$ and a commutative diagram
\begin{equation}
\label{eq:diagmpw}
\begin{CD}
R^{ss}_{pw}\times R^{ss}_{(m-p)w} @>{f_{p}^{R}}>> R^{ss}_{v}\\
@V{q_{pw}\times q_{(m-p)w}}VV                   @VV{q_{v}}V\\
M_{pw}\times  M_{(m-p)w}@>{f_{p}}>> M_{v}
\end{CD}
\end{equation}
and since the upper map is induced by the family $\overline{\mathcal{Q}}_{pw}\oplus\overline{\mathcal{Q}}_{(m-p)w}$ of quotients of $\mathcal{H}_{v}$, by construction we have $f_{p}([F],[G])=[F\oplus G]$ for every $[F]\in M_{pw}$ and $[G]\in M_{(m-p)w}$.

If we now let $G$ be an $H-$semistable sheaf whose Mukai vector is $(m-p)w$ and $Q_G\in R^{ss}_{(m-p)w}$ be such that $q_{(m-p)w}(Q_G)=[G]$, 
the commutative diagram (\ref{eq:diagmpw}) induces the following commutative diagram 
\begin{equation}
\label{eq:diagmpw2}
\begin{CD}
R^{ss}_{pw} @>{f_{p,Q_G}^{R}}>> R^{ss}_{v}\\
@V{q_{pw}}VV                   @VV{q_{v}}V\\
M_{pw}@>{f_{p,[G]}}>> M_{v}
\end{CD}
\end{equation}
where $f_{p,[G]}$ is the restriction of $f_{p}$ to $M_{pw}\times\{[G]\}\simeq M_{pw}$, and $f_{p,Q_G}^{R}$ is the restriction of $f_{p}^{R}$ to $R^{ss}_{pw}\times\{Q_G\}\simeq R^{ss}_{pw}$. Equivalently, if we let $p_S:S\times R^{ss}_{pw}\rightarrow S$ be the projection, the morphism $f_{p,Q_G}^{R}$ is the modular map to $R^{ss}_{v}$ associated with the family $\mathcal{Q}_{pw}\oplus p^*_S(Q_G)$. 

Since $v^{\perp}=(pw)^{\perp}=((m-p)w)^{\perp}$, the morphisms $f_{p,[G]}$ and $f_{p}$ allow us to compare $\lambda_v$, $\lambda_{pw}$ and $\lambda_{(m-p)w}$. This is the content of the following:

\begin{prop}
\label{prop:diag1} Let $(S,v,H)$ be an $(m,k)-$triple where $(m,k)\neq(2,1)$ and where $v=mw$ for a primitive Mukai vector $w$.
\begin{enumerate}
 \item For every $1\leq p<m$ and $G\in M_{(m-p)w}$, the diagram 
\begin{equation}
\label{eq:diagramma2gen}
\begin{CD}
v^{\perp} @>{\lambda_{v}}>> H^{2}(M_{v},\mathbb{Z})\\
@V{\rm id}VV                   @VV{f_{p,[G]}^{*}}V\\
(pw)^{\perp} @>>{\lambda_{pw}}> H^{2}(M_{pw},\mathbb{Z})
\end{CD}
\end{equation} 
is commutative, i. e. we have $\lambda_{pw}=f_{p,[G]}^{*}\circ\lambda_{v}$.
 \item Let $\pi_{i}$ for $i=1,2$ be the projection of $M_{pw}\times M_{(m-p)w}$ on the $i$-th factor. For every $1\leq p<m$ the diagram 
\begin{equation}
\label{eq:diagramma5gen}
\begin{CD}
v^{\perp}                                                        @>{\lambda_{v}}>>             H^{2}(M_{v},\mathbb{Z})\\
@V{(\lambda_{pw},\lambda_{(m-p)w})}VV                                                 @VV{f_{p}^{*}}V\\
H^{2}(M_{pw},\mathbb{Z})\oplus H^{2}(M_{(m-p)w},\mathbb{Z})  @>>{\pi_{1}^{*}+\pi_{2}^{*}}>  H^{2}(M_{pw}\times M_{(m-p)w},\mathbb{Z})
\end{CD}
\end{equation}
is commutative, i.e. $f_{p}^{*}\circ\lambda_{v}=\pi_{1}^{*}\circ\lambda_{pw}+\pi_{2}^{*}\circ\lambda_{(m-p)w}$.
\end{enumerate} 
\end{prop}

\proof To prove point (1), let $\alpha\in v^{\perp}$ and notice that by the characterization of $\lambda_{pw}(\alpha)$ given in Remark \ref{oss:carattlambdaII} to show that $f_{p,[G]}^{*}(\lambda_{v}(\alpha))=\lambda_{pw}(\alpha)$ we just need to show that $$j_{pw,\mathbb{Q}}^*(q_{pw}^*(f_{p,[G]}^{*}(\lambda_{v}(\alpha)))=j_{pw,\mathbb{Q}}^*(q_{pw}^*(\lambda_{pw}(\alpha))).$$

Using the commutativity of diagram (\ref{eq:diagmpw2}) we get $$j_{pw,\mathbb{Q}}^*(q_{pw}^*(f_{p,[G]}^{*}(\lambda_{v}(\alpha))))=j_{pw,\mathbb{Q}}^*(f_{p,Q_G}^{R*}(q_v^*(\lambda_v(\alpha)))).$$By Remark \ref{oss:carattlambda} (that we may apply since $(m,k)\neq(2,1)$) we get $$j_{pw,\mathbb{Q}}^*(f_{p,Q_G}^{R*}(q_v^*(\lambda_v(\alpha))))=j_{pw,\mathbb{Q}}^*(f_{p,Q_G}^{R*}(\mu_{\mathcal{Q}_v}(\alpha))).$$As $(id_S\times f_{p,Q_G}^{R})^{*}(\mathcal{Q}_v)\simeq\mathcal{Q}_{pw}\oplus p^*_S(Q_{G})$, by point (2) of Lemma \ref{lem:propmu} we get $$j_{pw,\mathbb{Q}}^*(f_{p,Q_G}^{R*}(\mu_{\mathcal{Q}_v}(\alpha)))=j_{pw,\mathbb{Q}}^*(\mu_{\mathcal{Q}_{pw}\oplus p^*_S(G)}(\alpha)).$$By point (1) of Lemma \ref{lem:propmu} we have $$j_{pw,\mathbb{Q}}^*(\mu_{\mathcal{Q}_{pw}\oplus p^*_S(Q_{G})}(\alpha))=j_{pw,\mathbb{Q}}^*(\mu_{\mathcal{Q}_{pw}}(\alpha))+j_{pw,\mathbb{Q}}^*(\mu_{p^*_S(Q_{G})}(\alpha)).$$By (2) of Lemma \ref{lem:propmu} and Remark \ref{oss:carattlambdaII} we have $$j_{pw,\mathbb{Q}}^*(\mu_{\mathcal{Q}_{pw}}(\alpha))=\mu_{\mathcal{Q}_{pw}^s}(\alpha)=j_{pw,\mathbb{Q}}^*(q_{pw}^*(\lambda_{pw}(\alpha))).$$Since by point (3) of Lemma \ref{lem:propmu} we have that $\mu_{p^*_S(Q_{G})}$ factors through the integral cohomology of a point, we get that $\mu_{p^*_S(Q_{G})}(\alpha)=0$, completing the proof of point (1). 

We now prove point (2). By Lemma \ref{lem:firstcohom} we know that either $$H^{1}(M_{pw},\mathbb{Z})\simeq H^{1}(M_{(m-p)w},\mathbb{Z})=0$$(if $S$ is K3), or $H^{1}(M_{pw},\mathbb{Z})$ and $H^{1}(M_{(m-p)w},\mathbb{Z})$ are free. In both cases $Tor(H^{j}(M_{pw},\mathbb{Z}),H^{3-j}(M_{(m-p)w},\mathbb{Z}))=0$ for $j\le3$, so the K\"unneth formula for integral cohomology implies that the $\mathbb{Z}-$module $H^{2}(M_{pw}\times M_{(m-p)w},\mathbb{Z})$ is isomorphic to $$H^{2}(M_{pw},\mathbb{Z})\oplus H^{2}(M_{(m-p)w},\mathbb{Z})\oplus(H^{1}(M_{pw},\mathbb{Z})\otimes H^{1}(M_{(m-p)w},\mathbb{Z}))$$and the restriction of this isomorphism to $H^{2}(M_{pw},\mathbb{Z})\oplus H^{2}(M_{(m-p)w},\mathbb{Z})$ is $\pi_{1}^{*}+\pi_{2}^{*}$.

For every $G'\in M_{pw}$ and $G''\in M_{(m-p)w}$ and for every $\alpha\in v^{\perp}$, by point (1) the restrictions of $f^{*}_{p}(\lambda_v(\alpha))$ to $M_{pw}\times [G'']$ and $[G']\times M_{(m-p)w}$ are $\lambda_{pw}(\alpha)$ and $\lambda_{(m-p)w}(\alpha)$ respectively: as a consequence, the components of $f^{*}_{p}(\lambda_v(\alpha))$ in $H^{2}(M_{pw},\mathbb{Z})$ and $H^{2}(M_{(m-p)w},\mathbb{Z})$ are $\lambda_{pw}(\alpha)$ and $\lambda_{(m-p)w}(\alpha)$ respectively. It remains to show that$$f^{*}_{p}(\lambda_v(\alpha))\in \pi_1^{*}(H^{2}(M_{pw},\mathbb{Z}))\oplus\pi_2^{*}(H^{2}(M_{(m-p)w},\mathbb{Z})).$$ 

To prove this, recall that by Corollary \ref{cor:iniet} the morphisms $$q_{pw}^{*}:H^{i}(M_{pw},\mathbb{Z})\rightarrow H^{i}(R^{ss}_{pw},\mathbb{Z}),$$ $$q_{(m-p)w}^{*}:H^{i}(M_{(m-p)w},\mathbb{Z})\rightarrow H^{i}(R^{ss}_{(m-p)w},\mathbb{Z})$$are injective for $i=1,2$. It follows that $$(q_{pw}\times q_{(m-p)w})^{*}:H^{2}(M_{pw}\times M_{(m-p)w},\mathbb{Z})\longrightarrow H^{2}(R^{ss}_{pw}\times R^{ss}_{(m-p)w},\mathbb{Z})$$is injective and we just need to show that $$(q_{pw}\times q_{(m-p)w})^{*}(f^{*}_{p}(\lambda_v(\alpha)))\in \pi_1^{R*}(H^{2}(R_{pw}^{ss},\mathbb{Z}))\oplus\pi_2^{R*}(H^{2}(R_{(m-p)w}^{ss},\mathbb{Z})).$$

To do so, notice that the commutativity of diagram (\ref{eq:diagmpw}) implies $$(q_{pw}\times q_{(m-p)w})^{*}(f^{*}_{p}(\lambda_v(\alpha)))=f^{R*}_{p}(q_{v}^{*}(\lambda_v(\alpha))).$$By Remark \ref{oss:carattlambda} and point (2) of Lemma \ref{lem:propmu} we have $$f^{R*}_{p}(q_{v}^{*}(\lambda_v(\alpha)))=f^{R*}_{p}(\mu_{\mathcal{Q}_{v}}(\alpha))=\mu_{(id_S\times f^{R}_{p})^{*}\mathcal{Q}_{v}}(\alpha).$$

Finally, by equation (\ref{eq:mancava}) and point (2) of Lemma \ref{lem:propmu} we obtain $$\mu_{(id_S\times f^{R}_{p})^{*}\mathcal{Q}_{v}}(\alpha)=\mu_{\overline{\mathcal{Q}}_{pw}\oplus\overline{\mathcal{Q}}_{(m-p)w}}(\alpha)=\pi_1^{R*}(\mu_{\mathcal{Q}_{pw}}(\alpha))+\pi_2^{R*}(\mu_{\mathcal{Q}_{(m-p)w}}(\alpha)).$$As this last is in $\pi_1^{R*}(H^{2}(R_{pw}^{ss},\mathbb{Z}))\oplus\pi_2^{R*}(H^{2}(R_{(m-p)w}^{ss},\mathbb{Z}))$, we are done.\endproof

If $S$ is an Abelian surface and $[G]\in K_{(m-p)w}$, then the two morphisms $f_{p,[G]}:M_{pw}\longrightarrow M_v$ and $f_{p}:M_{pw}\times M_{(m-p)w}\longrightarrow M_v$ induce by restrictions the morphisms $$f_{p,[G]}^{0}:K_{pw}\longrightarrow K_v,\,\,\,\,\,\,\,\,\,\,\,\,\,f_{p}^{0}:K_{pw}\times K_{(m-p)w}\longrightarrow K_v$$that allow to state the analogue of Proposition \ref{prop:diag1} in the Abelian case.

\begin{prop}
\label{prop:diag1ab}
Let $(S,v,H)$ be an $(m,k)-$triple where $S$ is an Abelian surface, $(m,k)\neq(2,1)$ and $v=mw$ for a primitive Mukai vector $w$.
\begin{enumerate}
 \item For every $1\leq p<m$ and $G\in K_{(m-p)w}$, the diagram 
\begin{equation}
\label{eq:diagramma2ab}
\begin{CD}
v^{\perp} @>{\lambda^{0}_{v}}>> H^{2}(K_{v},\mathbb{Z})\\
@V{\rm id}VV                   @VV{f_{p,[G]}^{0*}}V\\
(pw)^{\perp} @>>{\lambda^{0}_{pw}}> H^{2}(K_{pw},\mathbb{Z})
\end{CD}
\end{equation} 
is commutative, i. e. we have $\lambda^{0}_{pw}=f_{p,[G]}^{0*}\circ\lambda^{0}_{v}$.
 \item If we let $\pi_{i}^{0}$ for $i=1,2$ the projection of $K_{pw}\times K_{(m-p)w}$ to the $i$-th factor, for every $1\leq p<m$ the diagram 
\begin{equation}
\label{eq:diagramma5ab}
\begin{CD}
v^{\perp}                                                        @>{\lambda^{0}_{v}}>>             H^{2}(K_{v},\mathbb{Z})\\
@V{(\lambda^{0}_{pw},\lambda^{0}_{(m-p)w})}VV                                                 @VV{f_{p}^{0*}}V\\
H^{2}(K_{pw},\mathbb{Z})\oplus H^{2}(K_{(m-p)w},\mathbb{Z})  @>>{\pi_{1}^{0*}+\pi_{2}^{0*}}>  H^{2}(K_{pw}\times K_{(m-p)w},\mathbb{Z})
\end{CD}
\end{equation}
is commutative, i.e. $f_{p}^{0*}\circ\lambda^{0}_{v}=\pi_{1}^{0*}\circ\lambda^{0}_{pw}+\pi_{2}^{0*}\circ\lambda^{0}_{(m-p)w}$.
\end{enumerate}  
\end{prop}

\proof For point (1) notice that by definition of $f_{p,[G]}^{0}$ the diagram 
\begin{equation}
\label{eq:diagpk}
\begin{CD}
K_{pw} @>{f_{p,[G]}^{0}}>> K_{v}\\
@V{\iota_{pw}}VV                   @VV{\iota_{v}}V\\
M_{pw} @>>{f_{p,[G]}}> M_{v}
\end{CD}
\end{equation} 
is commutative, where $\iota_{u}:K_{u}\longrightarrow M_{u}$ is the inclusion. Since by definition $\lambda^{0}_{pw}=\iota_{pw}^*\circ\lambda_{pw}$ and $\lambda^{0}_{v}=\iota_{v}^*\circ \lambda_{v}$, the result follows from point (1) of Proposition \ref{prop:diag1}.

For point (2), recall by Lemma \ref{lem:firstcohom} that $$H^{1}(K_{pw},\mathbb{Z})=H^{1}(K_{(m-p)w},\mathbb{Z})=0.$$This implies that $Tor(H^{j}(K_{pw},\mathbb{Z}),H^{3-j}(K_{(m-p)w},\mathbb{Z}))=0$ for $j\leq 3$, so by the K\"unneth formula for integral cohomology we get $$H^{2}(K_{pw},\mathbb{Z})\oplus H^{2}(K_{(m-p)w},\mathbb{Z})\simeq H^{2}(K_{pw}\times K_{(m-p)w},\mathbb{Z})$$where the isomorphism is given by $\pi_{1}^{0*}+\pi_{2}^{0*}$. The commutative diagram
$$\begin{CD}
K_{pw}\times K_{(m-p)w} @>{f_{p}^{0}}>> K_{v}\\
@V{(\iota_{pw},\iota_{(m-p)w})}VV                   @VV{\iota_{v}}V\\
M_{pw}\times M_{(m-p)w} @>>{f_{p}}> M_{v}
\end{CD}$$
then implies that $$(f^{0}_{p})^{*}\circ\lambda^{0}_{v}=(f^{0}_{p})^{*}\circ\iota_{v}^{*}\circ\lambda_{v}=(\iota_{pw},\iota_{(m-p)w})^{*}\circ f^{*}_{p}\circ\lambda_{v}=$$ $$=(\iota_{pw},\iota_{(m-p)w})^{*}\circ(\pi_{1}^{*}\circ\lambda_{pw}+\pi_{2}^{*}\circ\lambda_{(m-p)w}),$$where this last equality follows from point (2) of Proposition \ref{prop:diag1}. Using now the following commutative diagram
$$\begin{CD}
K_{pw} @<{\pi_{1}^{0}}<< K_{pw}\times K_{(m-p)w} @>{\pi^{0}_{2}}>> K_{(m-p)w}\\
@V{\iota_{pw}}VV       @V{(\iota_{pw}\times\iota_{(m-p)w})}VV                   @VV{\iota_{(m-p)w}}V\\
M_{pw} @<<{\pi_{1}}<     M_{pw}\times M_{(m-p)w} @>>{\pi_{2}}>     M_{v}
\end{CD}$$
we see that $$(\iota_{pw},\iota_{(m-p)w})^{*}\circ(\pi_{1}^{*}\circ\lambda_{pw}+\pi_{2}^{*}\circ\lambda_{(m-p)w})=$$ $$=\iota_{pw}^{*}\circ\pi_{1}^{*}\circ\lambda_{pw}+\iota_{(m-p)w}^{*}\circ\pi_{2}^{*}\circ\lambda_{(m-p)w}=$$ $$=(\pi^{0}_{1})^{*}\circ\iota_{pw}^{*}\circ\lambda_{pw}+(\pi_{2}^{0})^{*}\circ\iota_{(m-p)w}^{*}\circ\lambda_{(m-p)w}=$$ $$=(\pi^{0}_{1})^{*}\circ\lambda^{0}_{pw}+(\pi^{0}_{2})^{*}\circ\lambda^{0}_{(m-p)w},$$and we are done.\endproof

\section{\label{sec4}The morphisms $\lambda_{v}$ and $\lambda^{0}_{v}$ are isomorphisms}

In this section we show that if $(S,v,H)$ is an $(m,k)-$triple, then if $S$ is K3 we have that $\lambda_{v}:v^{\perp}\longrightarrow H^{2}(M_{v},\mathbb{Z})$ is an isomorphism of $\mathbb{Z}-$modules, and if $S$ is Abelian and $(m,k)\neq(1,1),(1,2)$ then $\lambda^{0}_{v}:v^{\perp}\longrightarrow H^{2}(K_{v},\mathbb{Z})$ is an isomorphism of $\mathbb{Z}-$modules. 

As an important preliminary result we compute the second Betti numbers of these varieties and show that $b_{2}(M_{v})=23$ if $S$ is K3, and that $b_{2}(K_{v})=7$ if $S$ is Abelian and $(m,k)\neq(1,1),(1,2)$. Since this is know to hold if $m=1$ (see \cite{OG1} and \cite{Y1}) and if $(m,k)=(2,1)$ (see \cite{PR}), we suppose from now on that $(m,k)\neq(2,1)$, and that $m\neq 1$. 

We first show that $b_{2}(M_{v})\geq 23$ and that $b_{2}(K_{v})\geq 7$ as an easy consequence of the previous section:

\begin{lem}
\label{lem:almeno23}
Let $(S,v,H)$ be an $(m,k)-$triple, and if $S$ is Abelian suppose that $(m,k)\neq(1,1)$.
\begin{enumerate}
 \item If $S$ is K3, the morphism $\lambda_{v}:v^{\perp}\longrightarrow H^{2}(M_{v},\mathbb{Z})$ is injective. In particular $b_{2}(M_{v}(S,H))\geq 23$.
 \item If $S$ is Abelian, the morphism $\lambda^{0}_{v}:v^{\perp}\longrightarrow H^{2}(K_{v},\mathbb{Z})$ is injective. In particular $b_{2}(K_{v}(S,H))\geq 7$.
\end{enumerate}
\end{lem}

\proof The statement holds for $(1,k)-$triples (see \cite{OG1}, \cite{Y1}) and $(2,1)-$triples (see \cite{PR}). For the remaining cases, we start by considering $S$ to be a K3 surface. 

We write $v=mw$, and we choose an $H-$semistable sheaf $G$ of Mukai vector $(m-1)w$. As seen in the previous section we have a morphism $f_{1,[G]}:M_{w}\longrightarrow M_{v}$ (mapping $[F]$ to $[F\oplus G]$), and by Proposition \ref{prop:diag1} we have $f_{1,[G]}^{*}\circ\lambda_{v}=\lambda_{w}$. 

Since $(S,w,H)$ is a $(1,k)-$triple (by Lemma \ref{lem:genvw}), it follows that $\lambda_{w}$ is an isomorphism. In particular it is injective, so $\lambda_{v}$ is injective too. Since $rk(v^{\perp})=23$, this completes the proof for K3 surfaces.

Suppose now that $S$ is an Abelian surface. If $k\geq 2$ and $[G]\in K_{(m-1)w}$, by Proposition \ref{prop:diag1ab} we have $f_{1,[G]}^{0*}\circ\lambda^{0}_{v}=\lambda^{0}_{w}$. As again $(S,w,H)$ is a $(1,k)-$triple and $k\geq 2$, it follows that $\lambda^{0}_{w}$ is injective, so the morphism $\lambda^{0}_{v}$ is injective too.

The only remaining case is when $k=1$ and $m\geq 3$. In this case $K_{w}$ consists of a single point, so $\lambda^{0}_{w}=0$. However, by Theorem 1.7 of \cite{PR} the morphism $\lambda^{0}_{2w}$ is injective. Since $m\geq 3$ we have $m-2\geq 1$, so $K_{(m-2)w}\neq\emptyset$ and we may choose $[G]\in K_{(m-2)w}$. By Proposition \ref{prop:diag1ab} we know that $f_{2,[G]}^{0*}\circ\lambda^{0}_{v}=\lambda^{0}_{2w}$, so the injectivity of $\lambda^{0}_{v}$ follows.\endproof

The core of this section is devoted to prove for every $(m,k)-$triple $(S,v,H)$ the inequalities $b_{2}(M_{v})\leq 23$ (if $S$ is K3) and $b_{2}(K_{v})\leq 7$ (if $S$ is Abelian). By point (1) of Theorem \ref{thm:mio}, if $S$ is K3 the topological type of $M_{v}$  only depends on $m$ and $k$, and the same holds for $K_v$ if $S$ is Abelian: we then only need to show that $b_{2}(M_{v})\leq 23$ for one particular $(m,k)-$triple $(S,v,H)$ where $S$ is a K3, and that $b_{2}(K_{v})\leq 7$ for one particular $(m,k)-$triple $(S,v,H)$ where $S$ is an Abelian surface. 

We now review the definition of $K_{v}$ and fix the setting in order to deal with both cases at the same time.

\begin{oss}
\label{oss:defkv}
{\rm If $S$ is an Abelian surface, the definition of $K_v$ as the fiber over $(0_S,\mathcal{O}_{S})$ of the morphism 
$a_v: M_v\longrightarrow S\times\hat{S}$ depends on the choice of a sheaf $G$ whose S-equivalence class is in $M_v$ (see Section 2.2 of \cite{PR3}). Since $a_v$ is an isotrivial fibration, different choices of $G$ produce isomorphic copies of $K_v$.}

{\rm In the proof of the inequality $b_{2}(K_{v})\leq 7$ we will only consider those Abelian surfaces $S$ whose N\'eron-Severi group
is generated by the class of an ample line bundle $H$, and we need an explicit description of the S-equivalence classes in $M_v$ that belong to $K_v$. For this purpose we choose $G$ so that for every $[F]\in M_v$ we have
\begin{equation}
[F]\in K_v \;\; \Longleftrightarrow \;\; \det(F)\in\mathbb{Z}H\subset Pic(S)\;\;{\rm and}  \;\;\sum \left({\mathbf c}_2(F)\right)=0_{S}\in S 
\end{equation}
where $\det(F)$ is the determinant bundle of $F$, the cycle  ${\mathbf c}_2(F)\in CH_0(S)$ is the second Chern class of $F$ and 
$\sum:CH_0(S)\longrightarrow S$ is the Albanese morphism, i.e. the group morphism mapping $[p]$ to $p$ for all $p\in S$. The existence of such a $G$ follows from Lemma 2.10 of \cite{PR3}.}
\end{oss}

\begin{sett}
\label{sett:notazioni} 
{\rm Let $m,k\in\mathbb{N}$ such that $m\geq 2$, $k\geq 1$ and $(m,k)\neq(2,1)$. Let $S$ be  either a  projective K3 surface or an Abelian surface  such that $NS(S)=\mathbb{Z}\cdot h$, where $h$ is the first Chern class in cohomology of an ample line bundle $H$ and $h^{2}=2k$. If $S$ is Abelian further assume that the self-intersection $H^{2}$ of $H$ in the Chow ring of $S$ belongs to the kernel of the Albanese morphism $\sum:CH_0(S)\longrightarrow S$, i. e. $$\sum (H^2)=0_{S}\in S.$$}

{\rm We will use the following notation: we let $$M:=\left\{\begin{array}{ll} M_{m(0,h,0)}(S,H), & {\rm if}\,\,S\,\,{\rm is}\,\,{\rm K3}\\ K_{m(0,h,0)}(S,H), & {\rm if}\,\,S\,\,{\rm is}\,\,{\rm Abelian}\end{array}\right.$$and $$M':=\left\{\begin{array}{ll} M_{(0,mh,1-m^{2}k)}(S,H), & {\rm if}\,\,S\,\,{\rm is}\,\,{\rm K3}\\ K_{(0,mh,1-m^{2}k)}(S,H), & {\rm if}\,\,S\,\,{\rm is}\,\,{\rm Abelian}\end{array}\right.$$}

{\rm In both cases we set $P:=|mH|$, and we let $P^0\subset P$ be the open subvariety parameterizing smooth curves and $P^1\subset P$ the open subvariety parameterizing those curves which are either smooth or nodal with at most one node. Notice that $P^{0}\subseteq P^{1}$.}

{\rm We let $\phi:M\longrightarrow P$ and $\phi':M'\longrightarrow P$ be the morphisms mapping an S-equivalence class of semistable sheaves to the corresponding Fitting subscheme (see section 2.2 of \cite{Lp1} and section 2.3 of \cite{Lp2}). We moreover set $$M^0=\phi^{-1}(P^0),\,\,\,\,\,\,\,\,\,M^1=\phi^{-1}(P^1),$$ $$(M')^0=(\phi')^{-1}(P^0),\,\,\,\,\,\,\,\,\,(M')^1=(\phi')^{-1}(P^1).$$Notice again that $M^{0}\subseteq M^{1}\subseteq M$ and $(M')^{0}\subseteq(M')^{1}\subseteq M'$. Finally, we let $J\subset M^1$ and $J'\subset(M')^1$ be the open subvarieties parameterizing sheaves that are push-forward on S of line bundles on curves in $P^1$.}
\end{sett}

Since the square (with respect to the Mukai pairing) of the Mukai vector defining $M$ and $M'$ is $2km^{2}$, we have $$\dim(M)=\dim (M')=2\dim(P),$$and since $(0,mh,1-m^{2}k)$ is primitive $M'$ is an irreducible holomorphic symplectic manifold. It follows in particular that the open subvarieties $M'^0$, $M'^1$ and $J'$ of $M'$ are smooth.

Although $M$ is only an irreducible symplectic variety, since a pure one-dimensional rank-one torsion-free coherent sheaf on an irreducible curve has no one-dimensional quotients, the open subvarieties $M^0$, $M^1$ and $J$ are contained in the stable locus $M^s$ of $M$, which is moreover the smooth locus of $M$: as a consequence $M^0$, $M^1$ and $J$ are smooth as well.

We now need to recall some properties of the objects introduced in Setting \ref{sett:notazioni}, and we collect them in the following three remarks. The first one concerns the linear system $P$.

\begin{oss}
\label{oss:sistlin}
{\rm Since $m\geq 2$, $k\geq 1$ and $(m,k)\neq (2,1)$, the complete linear system $P=|mH|$ is very ample: this follows from Theorem 5.2 and Theorem 6.1 (iii) of \cite{SD} if $S$ is a K3, and from the Main Theorem of \cite{Ram} if $S$ is an Abelian surface. As a consequence, the locus $P^0$ parametrizing smooth curves is a non empty open subset and its complement $P\setminus P^0$ is an irreducible divisor containing as a dense, smooth, locally closed subvariety the locus $P^1\setminus P^0$ parametrizing nodal curves with exactly one node. In particular $P\setminus P^1$ is a closed subvariety of codimension at least $2$ in $P$.}      
\end{oss} 

The second remark describes the sheaves corresponding to the points of $J$ and $J'$.

\begin{oss}
\label{oss:relj} 
{\rm The arithmetic genus of the curves in $P$ is $g:=km^2+1$, independently on the fact that $S$ is K3 or Abelian.} 

{\rm If $S$ is K3, the smooth varieties $J'$ and $J$ are identified with the relative generalized Jacobians of degree $1$ and $g-1$, respectively, on the family of irreducible curves parameterized by $P^1$. The varieties $M^1$ and $(M')^1$ are the corresponding relative compactified Jacobians.}

{\rm If $S$ is an Abelian surface, the definition of $K_v$ (see \ref{oss:defkv}) and the assumption $\sum(H^2)=0_{S}\in S$ imply that a sheaf $F$ belongs to $J$ (resp. to $J'$) if and only if there are a curve $C\in P^1$ with an embedding $\iota:C\longrightarrow S$ and a Cartier divisor $D=\sum n_i p_i$ of degree $g-1=km^{2}$ (resp. of degree 1) on $C$ such that $F\simeq\iota_*\mathscr{O}_{C}(D)$ and $\sum n_i\iota(p_i)=0_{S}\in S$.} 
\end{oss}

The last one concerns a relation between $J$ and $J'$.

\begin{oss}
\label{oss:tensj} 
{\rm Let $C$ be a curve in $P^{1}$, $\iota:C\longrightarrow S$ the inclusion and $L\in Pic(C)$ a line bundle of degree 1. If $S$ is Abelian suppose furthermore that there is a Cartier divisor $D=\sum n_{i}p_{i}$ on $C$ such that $\sum n_{i}\iota(p_{i})=0_{S}$ and $L\simeq\mathscr{O}_{C}(D)$. By the previous remark we know that $\iota_{*}L\in J'$, and that every point of $J$ may be represented in this way.} 

{\rm Notice that $L^{\otimes km^{2}}$ is a line bundle of degree $g-1$ on $C$, and if $S$ is Abelian then $L^{\otimes km^{2}}=\mathscr{O}_{C}((g-1)D)$, where $(g-1)D=\sum(g-1)n_{i}p_{i}$ is such that $\sum(g-1)n_{i}\iota(p_{i})=0_{S}$. It then follows that $\iota_{*}L^{\otimes km^{2}}\in J$. We have then constructed a map $$\psi:J'\longrightarrow J$$over $P^{1}$ that will be of the utmost importance in what follows. This map is in fact a $P^{1}-$morphism of varieties by the universal property of the relative Jacobian.}
\end{oss}

The previous Remark is the reason why we consider the Mukai vectors we have chosen in Setting \ref{sett:notazioni}: the morphism $\psi:J'\longrightarrow J$ will be used to compare the second Betti numbers of $M$ and $M'$, and it turns out to be quasi-finite. This is the content of the following:

\begin{prop}
\label{prop:psi}
We use the notation introduced in Setting \ref{sett:notazioni}.  
\begin{enumerate}
 \item The $P^{1}-$morphism $\psi:J'\longrightarrow J$ is \'etale.
 \item The restriction $\psi^0:(M')^0\longrightarrow M^0$ of $\psi$ to $(M')^{0}$ is a proper \'etale morphism of degree $(g-1)^{2g}$ if $S$ is a K3, and of degree $(g-1)^{2g-4}$ if $S$ is Abelian. 
\end{enumerate}
\end{prop}

\proof In $S$ is K3, since the relative generalized Jacobians $J'$ and $J$ have the same dimension, to prove that $\psi$ is \'etale it is sufficient to prove that the differential of $\psi$ is an isomorphism at each point of $J'$. 

To show this, notice that if $C\in P^{1}$ the fibers $\phi^{-1}(C)$ and $(\phi')^{-1}(C)$ are the generalized Jacobians of degree $g-1$ and $1$, respectively, of $C$. The restriction $\psi_{|(\phi')^{-1}(C)}:(\phi')^{-1}(C)\longrightarrow\phi^{-1}(C)$ of the morphism $\psi$ to the fibers is the $(g-1)-$th tensor power, hence it is \'etale. Since $J'$ and $J$ are smooth over $P^1$ and have the same dimension, it follows that the differential of $\psi$ is an isomorphism at each point of $J'$.

If $S$ is Abelian, then $J'$ is a closed subvariety of the degree $1$ relative generalized Jacobian of the family of curves parameterized by $P^{1}$, while $J$ is a closed subvariety of the degree $g-1$ relative generalized Jacobian of the family of curves parametrized by $P^1$. 

The morphism $\psi$ is \'etale because it is the restriction of the \'etale morphism between the relative generalized Jacobians to smooth subvarieties having the same dimension. It follows that the differential of $\psi$ is injective at every point, so $\psi$ is \'etale.

For point (2), we immediately notice that the fact that $\psi^0$ is \'etale is a consequence of the first item of the statement. As $P^0$ parametrizes smooth curves, the restrictions $(\phi')^0:(M')^0\longrightarrow P^0$
and $\phi^0:M^0\longrightarrow P^0$ of $(\phi')^0$ and $\phi^0$ respectively are proper. Since $\phi'^0=\phi^0\circ \psi^0$, it follows that the morphism $\psi^0$ is proper as well. 

In order to determine the degree of $\psi^{0}$, recall that if $L$ and $L'$ are line bundles on a smooth curve $C\in P^{0}$, the tensor powers $L^{\otimes r}$ and $(L')^{\otimes r}$ are isomorphic if and only if $L$ and $L'$ differ for a $r-$torsion point of the degree $0$ Jacobian of $C$. Since the dimension of the Jacobian of $C$ is $g$, it has exactly $(g-1)^{2g}$ torsion points of order $g-1$.

Now, if $S$ is K3 it follows that the degree of $\psi^{0}$ is $\left(g-1\right)^{2g}$. If $S$ is Abelian and $L$ and $L'$ represent points of $J'$, by Remark \ref{oss:relj} we have $$L\otimes(L')^{\vee}\simeq\mathcal{O}_{C}(\sum n_i p_i)$$for a divisor $\sum n_i p_i$ on $C$ such that $\sum n_i \iota(p_i)=0_{S}\in S$, where $\iota: C\rightarrow S$ is the closed embedding.
 
Equivalently, $L\otimes(L')^{\vee}\in\ker(Alb(\iota))$, where $Alb(\iota):Jac(C)\longrightarrow S$ is the morphism induced by $\iota$
between the Albanese varieties $Jac(C)$ and $S$ of $C$ and $S$. It follows that the degree of $\psi^0$ is the number of the $(g-1)$-torsion points of $\ker(Alb(\iota))$.
  
But since $C$ is a smooth ample divisor of $S$, the Lefschetz Hyperplane Section Theorem implies that $\iota^*(H^1(S,\mathbb{Z}))$ is a saturated submodule of $H^1(C,\mathbb{Z})$. It follows that $\ker(Alb(\iota))$ is a connected subgroup of codimension $2$ in $Jac(C)$ and, as $\ker(Alb(\iota))$ is an Abelian variety of dimension $g-2$, the cardinality of its $(g-1)$-torsion subgroup is $\left(g-1\right)^{2g-4}$.\endproof 

To prove that $b_2(M)\leq b_2(M')$ we will first prove that $b_2(M)=b_2(J)$ and $b_2(M')=b_2(J')$ (by using that $M\setminus J$ and $M'\setminus J'$ have codimension at least two in $M$ and $M'$ respectively), and then that $b_2(J)\leq
b_2(J')$ (by using the morphism $\psi:J'\longrightarrow J$). These two facts are the content of the two following Lemmas.

\begin{lem}
\label{lem:b2eb2aperti}
We use the notation introduced in Setting \ref{sett:notazioni}.
\begin{enumerate}
 \item The open embedding $J\subset M$ induces an isomorphism $H^{2}(M,\mathbb{Z})\simeq H^{2}(J,\mathbb{Z})$. In particular we have $b_{2}(M)=b_{2}(J)$.
 \item The open embedding $J'\subset M'$ induces an isomorphism $H^{2}(M',\mathbb{Z})\simeq H^{2}(J',\mathbb{Z})$. In particular we have $b_{2}(M')=b_{2}(J')$.
\end{enumerate}
\end{lem}

\proof We only prove (1), the proof of (2) is similar and a bit easier since $M'$ is smooth.

We have a chain of open embeddings $J\subset M^{1}\subset M^{s}\subset M$: we show that each of these embeddings induces an isomorphism on the second integral cohomology groups. 

We first show that the inclusion $J\subseteq M^{1}$ induces an isomorphism $H^{2}(M^{1},\mathbb{Z})\simeq H^{2}(J,\mathbb{Z})$. As $M^{1}$ is smooth, this holds if $M^{1}\setminus J$ has codimension at least two in $M^{1}$.

To show this, recall that $M^{1}$ is irreducible, and has dimension $2g$ if $S$ is K3 and dimension $2g-4$ if $S$ is Abelian. Moreover we have $J\supset M^{0}=\phi^{-1}(P^0)$: as a consequence, it will be enough to show that for every $C\in P^1\setminus P^0$
the intersection $\phi^{-1}(C)\cap (M^{1}\setminus J)$ has dimension $g-1$ if $S$ is K3 and dimension $g-3$ if $S$ is Abelian.

Now, as $C\in P^{1}\setminus P^{0}$ implies that the singular locus of $C$ consists of a unique node, it follows that every point in $\phi^{-1}(C)\cap M^{1}\setminus J$ represents a sheaf of the form $\iota_*\nu_* (L)$ where $\nu:\tilde{C}\longrightarrow C$ is the normalization of $C$, the morphism $\iota:C\longrightarrow S$ is the closed embedding of $C$ in $S$ and $L$ is a line bundle  of degree $g-2$ on the curve $\tilde{C}$ (which is a curve of genus $g-1$).

If $S$ is K3 this implies that $\phi^{-1}(C)\cap (M^{1}\setminus J)$ is isomorphic to the Jacobian $J(\tilde{C},g-2)$ of degree $g-1$ over $\tilde{C}$, so its dimension is $g-1$ as desired.

If $S$ is Abelian, this together with Remark \ref{oss:defkv} imply that $\phi^{-1}(C)\cap (M^{1}\setminus J)$ is isomorphic to the subvariety of $J(\tilde{C},g-2)$ consisting of the line bundles $L$ on $\widetilde{C}$ of degree $g-2$ and such that $\sum(\mathbf{c}_2(\iota_*\nu_* (L)))=0_{S}\in S$.

Consider then the non-constant regular morphism $h:J(\tilde{C},g-2)\longrightarrow S$ mapping $L\in J(\tilde{C},g-2)$ to $\sum(\mathbf{c}_2(\iota_*\nu_* (L)))$. Since both $J(\tilde{C},g-2)$ and $S$ are Abelian varieties and the surface $S$ does not contain any elliptic curve, it follows that $h$ is surjective. Since $\phi^{-1}(C)\cap (M^{1}\setminus J)$ is isomorphic to a fiber of $h$, its dimension is then $g-3$. This concludes the proof that the inclusion $J\subseteq M^{1}$ induces an isomorphism $H^{2}(M^{1},\mathbb{Z})\simeq H^{2}(J,\mathbb{Z})$. 

We now prove that the inclusion $M^{1}\subseteq M^{s}$ induces an isomorphism $H^{2}(M^{s},\mathbb{Z})\simeq H^{2}(M^{1},\mathbb{Z})$. Again, this follows if we show that $M^s\setminus M^1$ has codimension at least two in $M^s$, i.e. $M^s\setminus M^1$ does not contain any divisor of $M^s$.

To prove this, suppose that $D$ is an irreducible divisor of $M^s$ contained in $M^s\setminus M^1$ and let $\overline{D}\subset M$ be its closure in $M$. Since $M$ is a locally factorial irreducible symplectic variety (see Proposition A.2 of \cite{PR3}), by point (4) of Corollary 1.4 of \cite{Ma3} its image $\phi(\overline{D})$ in $P$ has codimension at most one in $P$. But $D\subset\overline{D}\subset (M\setminus M^1)=\phi^{-1}(P\setminus P^1)$, so $\phi(\overline{D})\subset (P\setminus P^1)$ and, as seen in Remark \ref{oss:sistlin}, we have that $P\setminus P^1$ has codimension at least two in $P$, so we get a contradiction.

We are only left with the proof of the fact that the inclusion $M^{s}\subseteq M$ induces an isomorphism $H^{2}(M,\mathbb{Z})\simeq H^{2}(M^{s},\mathbb{Z})$: this is the content of Proposition \ref{prop:extmr}.\endproof

As previously announced, we now prove the relation between the second Betti numbers of $J$ and $J'$, which is precisely as follows:

\begin{lem}
\label{lem:b2U}
We have that $b_{2}(J)\leq b_{2}(J')$. 
\end{lem}

\proof The proof is divided in two main steps: we first show that the divisor $J\setminus M^0=\phi^{-1}((P^1\setminus P^0))$ of $J$ and the divisor and $J'\setminus(M')^0=(\phi')^{-1}((P^1\setminus P^0))$ of $J'$ are both are irreducible and smooth. We will then prove the desired inequality.

\textit{Step 1.} We only prove that $J\setminus M^{0}$ is an irreducible, smooth divisor of $J$, the case of $J'$ being similar. To do so, notice that since divisors have pure dimension and $P^1\setminus P^0$ is smooth and irreducible by Remark \ref{oss:sistlin}, it will be enough to show that for every $C\in P^{1}\setminus P^{0}$ we have that $\phi^{-1}(C)\cap J$ is smooth and irreducible of dimension equal to the dimension of $P$ (which is $g$ if $S$ is K3 and $g-2$ if $S$ is Abelian).

If $S$ is K3 we know that $\phi^{-1}(C)\cap J$ is isomorphic to the generalized Jacobian $J(C,g-1)$ of line bundles of degree $g-1$ on $C$. Since $C$ is an irreducible, nodal curve with a unique node, it follows that $J(C,g-1)$ is irreducible, smooth and its dimension is $g$.

If $S$ is Abelian then $\phi^{-1}(C)\cap J$ is isomorphic to $a^{-1}(0_{S})$, where $a$ is the morphism $$a:J(C,g-1)\longrightarrow S,\,\,\,\,\,\,\,a(\mathscr{O}_{C}(D)):=\sum\iota_{*}(D)\in S,$$where $D$ is a Cartier divisor on $C$ and $\iota:C\longrightarrow S$ is the closed embedding of $C$ in $S$ (see Remark \ref{oss:relj}).

Consider now the normalization $\nu:\tilde{C}\longrightarrow C$ of $C$, and let $J(\tilde{C},g-1)$ be the Jacobian of line bundles of degree $g-1$ on the smooth curve $\tilde{C}$. The morphism $\nu^{*}:J(C,g-1)\longrightarrow J(\tilde{C},g-1)$
makes $J(C,g-1)$ a $\mathbb{C}^{*}$-bundle over $J(\tilde{C},g-1)$. Moreover, if we let $$\tilde{a}:J(\tilde{C},g-1)\longrightarrow S,\,\,\,\,\,\,\,\tilde{a}(\mathscr{O}_{\tilde{C}}(\tilde{D})):=\sum\iota_{*}\nu_{*}(\tilde{D})$$for every divisor $\tilde{D}$ on $\tilde{C}$, we have $\widetilde{a}\circ\nu^{*}=a$. 

Up to translation, the morphism $\tilde{a}$ may be identified with the morphism $Alb(\iota\circ\nu):Jac(\tilde{C})\longrightarrow S$
induced by $\iota\circ\nu$ on the Albanese varieties of $\tilde{C}$ and $S$. It follows that $\phi^{-1}(C)\cap J$ is a $\mathbb{C}^{*}-$bundle over a fiber of the group morphism $Alb(\iota\circ\nu)$: as a consequence, it is smooth.

To conclude this first step it only remains to show that $\phi^{-1}(C)\cap J$ is irreducible of dimension $g-2$, or equivalently that the kernel of $Alb(\iota\circ\nu)$ is connected of dimension $g-3$. To do so, it is enough to prove that $$(\iota\circ\nu)_*: H_{1}(\tilde{C},\mathbb{Z})\longrightarrow H_{1}(S,\mathbb{Z})$$is surjective.

To prove the surjectivity of $(\iota\circ\nu)_{*}$, let $b:\tilde{S}\longrightarrow S$ be the blow up of $S$ at the singular point of $C$, and let $i:\tilde{C}\longrightarrow\tilde{S}$ be the embedding of $\tilde{C}$ in $\tilde{S}$ as strict transform of $C$. Notice that $\iota\circ\nu =b\circ i$ and that $b_*:H_{1}(\tilde{S},\mathbb{Z})\longrightarrow H_{1}(S,\mathbb{Z})$ is surjective. To conclude we then only need to prove that $i_*: H_{1}(\tilde{C},\mathbb{Z})\longrightarrow H_{1}(\tilde{S},\mathbb{Z})$ is surjective. 

But since $C\in |mH|$ and the node of $C$ is the center of the blow-up morphism $b$, the curve $\tilde{C}$ on $\tilde{S}$ is linearly equivalent to $b^{*}(mH)-2E$, where $E$ is the exceptional divisor of $b$. Since the Seshadri constant of $H$ is strictly bigger than 1 by point (B) of Theorem A.1 of \cite{BS}, the divisor $b^{*}(mH)-2E$ is ample. The Lefschetz Hyperplane Section Theorem implies then that $i_*$ is surjective, concluding the proof of Step 1.
 
\textit{Step 2}. We now compare the second Betti numbers of $J$ and $J'$ by using the commutative diagram 
$$\begin{CD}
(M')^0 @>{h'}>> J'\\
@VV{\rm \psi^0}V                   @VV{\rm \psi}V\\
M^{0} @>{h}>> J
\end{CD}$$
where $h$ and $h'$ are the natural inclusions, and the commutative diagram with exact rows it induces in rational cohomology by considering the pairs $(J,M^{0})$ and $(J',(M')^{0})$: 
\begin{equation} 
\label{eq:diacha}
\begin{CD}
H^{2}(J,M^0) @>{g}>> H^{2}(J) @>{h^*}>> H^{2}(M^0) \\
@VV{\rm \overline{\psi}^*}V @VV{\rm \psi^*}V                   @VV{\rm \psi^{0*}}V\\
H^{2}(J',(M')^0) @>{g'}>> H^{2}(J') @>{h'^*}>> H^{2}((M')^0)
\end{CD}
\end{equation}

To prove that $b_{2}(J)\leq b_{2}(J')$ we will prove that $\psi^*$ is injective. By Step 1 of the proof we know that the divisors $J\setminus M^0$ of $J$ and $J'\setminus(M')^0$ of $J'$ are both smooth and irreducible, and we moreover know that $J$ and $J'$ are smooth varieties. Excision and Thom isomorphism then give $$H^{2}(J,M^0)\simeq H^{0}(J\setminus M^0)\simeq \mathbb{Q},\,\,\,\,\,\,H^{2}(J',(M')^0)\simeq H^{0}(J'\setminus(M')^0)\simeq \mathbb{Q}.$$As  $\psi$ is \'etale by point (1) of Proposition \ref{prop:psi}, the morphism $\overline{\psi}^*$ is an isomorphism. Moreover, by point (2) of Proposition \ref{prop:psi} the morphism $\psi^{0}$ is a topological covering: it follows that $\psi^{0*}:H^{2}(M^0)\longrightarrow H^{2}((M')^0)$ is injective.

We claim that $g':H^{2}(J',(M')^0)\longrightarrow H^{2}(J')$ is injective as well. To prove this, recall that if $\gamma$ is a generator of $H^{2}(J',(M')^0)$, then $g'(\gamma)$ is the first Chern class (in rational cohomology) of the line bundle associated with the irreducible divisor $J'\setminus (M')^0$ of the manifold $J'$. But then $g'(\gamma)$ extends to the non-zero class of an effective divisor on the projective manifold $M'$. By point (2) of Lemma \ref{lem:b2eb2aperti} we know that the open embedding $J'\subset M'$ induces an isomorphism $H^{2}(M')\simeq H^{2}(J')$: it follows that $g'(\gamma)\ne 0$, so $g'$ is injective.

A diagram chasing argument on the diagram (\ref{eq:diacha}) will now conclude the proof of the injectivity of $\psi^*$.\endproof 

We are now in the position to prove the following, which is Theorem \ref{thm:b2} of the Introduction.

\begin{thm}
\label{thm:b22}
Let $(S,v,H)$ be an $(m,k)-$triple. 
\begin{enumerate}
 \item If $S$ is K3, then $b_{2}(M_{v})=23$.
 \item If $S$ is Abelian and $(m,k)\neq(1,1),(1,2)$, then $b_{2}(K_{v})=7$.
\end{enumerate}
\end{thm}

\proof We already know this for $(1,k)-$triples (see \cite{OG1}, \cite{Y1}) and for $(2,1)-$triples (see \cite{PR}). Hence we suppose that either $m>2$ or $m=2$ and $k>1$.

By point 1 of Theorem \ref{thm:mio} we know that $M_{v}=M_{v}(S,H)$ is deformation equivalent to the moduli space $M$, hence $b_{2}(M_{v})=b_{2}(M)$. Now, by point (1) of Lemma \ref{lem:b2eb2aperti} we see that $b_{2}(M)=b_{2}(J)$, and by Lemma \ref{lem:b2U} we see that $b_{2}(J)\leq b_{2}(J')$. By point (2) of Lemma \ref{lem:b2eb2aperti} we have $b_{2}(J')=b_{2}(M')$. 

Now, recall that $M'$ is the moduli space associated with a $(1,m^{2}k)-$triple: if $S$ is K3 then $b_{2}(M')=23$, while if $S$ is Abelian then $b_{2}(M')=7$. By Lemma \ref{lem:almeno23} the statement follows.\endproof

\begin{oss}
\label{oss:satur}
{\rm If $S$ is Abelian, Theorem \ref{thm:b22} holds only for $(m,k)-$triples $(S,v,H)$ with $(m,k)\neq(1,1),(1,2)$. Indeed if $(m,k)=(1,1)$ then $K_{v}$ is a point, while if $(m,k)=(1,2)$ we have that $b_{2}(K_{v})=22$ as $K_{v}$ is a K3 surface, and more precisely it is a Kummer surface by Theorem 3.2 of \cite{Y6}.}

{\rm If $(m,k)=(1,2)$ the morphism $\lambda^{0}_{v}:v^{\perp}\longrightarrow H^{2}(K_{v},\mathbb{Z})$ is injective and its image is saturated in $H^{2}(K_{v},\mathbb{Z})$. To prove this, notice that it is enough to prove it for a particular $(1,2)-$triple, since this property is invariant under deformation.}

{\rm Choose then $v=(1,0,-2)$, so that $K_{v}(S,H)$ is the Kummer surface $Kum(S)$ of $S$. Notice that we have a natural embedding $$f:H^{2}(S,\mathbb{Z})\longrightarrow H^{2}(Kum(S),\mathbb{Z})$$whose image is saturated, which realizes $H^{2}(S,\mathbb{Z})$ as the orthogonal complement of the sublattice of $H^{2}(Kum(S),\mathbb{Z})$ spanned by the 16 $(-2)-$curves $E_{1},\cdots,E_{16}$ of $Kum(S)$(see \cite{B3}, Exp. VIII).}

{\rm We have $v^{\perp}\simeq H^{2}(S,\mathbb{Z})\oplus^{\perp}\mathbb{Z}(1,0,2)$. A direct computation shows moreover that $\lambda^{0}_{v}((0,\alpha,0))=f(\alpha)$ for every $\alpha\in H^{2}(S,\mathbb{Z})$, and that $\lambda^{0}_{v}(1,0,2)=e:=\frac{1}{2}\sum_{i=1}^{16}E_{i}$. If now $\beta\in Im(\lambda^{0}_{v})$, there are $\alpha\in H^{2}(S,\mathbb{Z})$ and $p\in\mathbb{Z}$ such that $n\beta=f(\alpha)+pe$. By intersecting with $E_{i}$ for some $1\leq i\leq 16$ we see that $p=qn$ for some $q\in\mathbb{Z}$. It follows that $n(\beta-qe)\in Im(f)$, and as $f$ is saturated it follows that $\beta-qe\in Im(f)$. But this implies that $\beta\in Im(\lambda^{0}_{v})$, completing the proof.}
\end{oss}

The previous Theorem allows us to prove that the morphism $\lambda_{v}$ (resp. $\lambda^{0}_{v}$) is an isomorphism of $\mathbb{Z}-$modules (preserving the Hodge structures):

\begin{prop}
\label{prop:b2leq}
Let $(S,v,H)$ be an $(m,k)-$triple, and if $S$ is Abelian suppose that $(m,k)\neq(1,1),(1,2)$. 
\begin{enumerate}
 \item If $S$ is K3, the morphism $\lambda_{v}:v^{\perp}\longrightarrow H^{2}(M_{v},\mathbb{Z})$ is an isomorphism of $\mathbb{Z}-$modules which preserves the Hodge structures.
 \item If $S$ is Abelian, the morphism $\lambda^{0}_{v}:v^{\perp}\longrightarrow H^{2}(K_{v},\mathbb{Z})$ is an isomorphism of $\mathbb{Z}-$modules which preserves the Hodge structures.
\end{enumerate}
\end{prop}

\proof The fact that $\lambda_{v}$ and $\lambda_{v}^{0}$ are Hodge morphisms is already known by Proposition \ref{prop:lambdav}, so we are left with showing that they are isomorphisms of $\mathbb{Z}-$modules. We already know this for $(1,k)-$triples (see \cite{OG1}, \cite{Y1}) and for $(2,1)-$triples (see \cite{PR}). Hence we suppose that either $m>2$ or $m=2$ and $k>1$. We start with the case of K3 surfaces.

Write $v=mw$ and choose $G\in M_{(m-1)w}$. In section 4.2 we constructed a morphism $f_{1,[G]}:M_{w}\longrightarrow M_{v}$ which makes the following diagram commutative
\begin{equation}
\begin{CD}
v^{\perp} @>{\lambda_{v}}>> H^{2}(M_{v},\mathbb{Z})\\
@V{\rm id}VV                   @VV{f_{1,[G]}^{*}}V\\
w^{\perp} @>>{\lambda_{w}}> H^{2}(M_{w},\mathbb{Z})
\end{CD}
\end{equation}
(see point (1) of Proposition \ref{prop:diag1}). All the $\mathbb{Z}-$modules in this diagram are free of rank 23 by Theorem \ref{thm:b22}, and $\lambda_{w}$ is an isomorphism of $\mathbb{Z}-$modules: it follows that $\lambda_{v}$ and $f_{1,[G]}^{*}$ are isomorphisms of $\mathbb{Z}-$modules as well.

The case of Abelian surfaces can be treated in a very similar way. If $k>2$, we proceed as before and use the morphism $f^{0}_{1,[G]}:K_{w}\longrightarrow K_{v}$ constructed in section 4.2, so that the following diagram is commutative
\begin{equation}
\label{eq:ddabab}
\begin{CD}
v^{\perp} @>{\lambda^{0}_{v}}>> H^{2}(K_{v},\mathbb{Z})\\
@V{\rm id}VV                   @VV{(f^{0}_{1,[G]})^{*}}V\\
w^{\perp} @>>{\lambda^{0}_{w}}> H^{2}(K_{w},\mathbb{Z})
\end{CD}
\end{equation}
(see point (1) of Proposition \ref{prop:diag1ab}). All the $\mathbb{Z}-$modules in this diagram are free of rank 7 by Theorem \ref{thm:b22} and $\lambda^{0}_{w}$ is an isomorphism of $\mathbb{Z}-$modules, so we conclude as before. 

If $k=1$, then $m\geq 3$ and as in the proof of Lemma \ref{lem:almeno23} we choose $G\in K_{(m-2)w}$. Using the morphism $f^{0}_{2,[G]}:K_{2w}\longrightarrow K_{v}$ we have by point (1) of Proposition \ref{prop:diag1ab} a commutative diagram
$$\begin{CD}
v^{\perp} @>{\lambda^{0}_{v}}>> H^{2}(K_{v},\mathbb{Z})\\
@V{\rm id}VV                   @VV{(f^{0}_{2,[G]})^{*}}V\\
(2w)^{\perp} @>>{\lambda^{0}_{2w}}> H^{2}(K_{2w},\mathbb{Z})
\end{CD}$$
where all the $\mathbb{Z}-$modules are free of rank 7 by Theorem \ref{thm:b22} and $\lambda^{0}_{w}$ is an isomorphism of $\mathbb{Z}-$modules. We then conclude as before even in this case.

The only remaining case is when $k=2$, so that $m\geq 2$. The $\mathbb{Z}-$modules in diagram (\ref{eq:ddabab}) are all free, and all of them have rank 7 but $H^{2}(K_{w},\mathbb{Z})$, which has rank 22. Anyway we know that $\lambda^{0}_{w}$ is injective by point (2) of Lemma \ref{lem:almeno23} and its image is saturated in $H^{2}(K_{w},\mathbb{Z})$ by Remark \ref{oss:satur}. 

The fact that $\lambda^{0}_{w}$ is injective implies that $\lambda^{0}_{v}$ is injective. We now show that it is surjective, so take $\beta\in H^{2}(K_{v},\mathbb{Z})$. As $v^{\perp}$ and $H^{2}(K_{v},\mathbb{Z})$ are free of the same rank by Theorem \ref{thm:b22}, there is $n\geq 1$ such that $n\beta=\lambda^{0}_{v}(\alpha)$ for some $\alpha\in v^{\perp}$. Hence $$n(f^{0}_{m-1,[G]})^{*}(\beta)=(f^{0}_{m-1,[G]})^{*}(\lambda^{0}_{v}(\alpha))=\lambda^{0}_{w}(\alpha),$$where the last equality follows from the commutativity of diagram (\ref{eq:ddabab}). 

As the image of $\lambda^{0}_{w}$ is saturated in $H^{2}(K_{w},\mathbb{Z})$ it follows that there is $\alpha'\in v^{\perp}$ such that $(f^{0}_{m-1,[G]})^{*}(\beta)=\lambda^{0}_{w}(\alpha')$. By injectivity of $\lambda^{0}_{w}$ we see that $\alpha=n\alpha'$, so that $\lambda^{0}_{v}(\alpha')=\beta$.\endproof

\begin{cor}
\label{cor:abh2}
Let $(S,v,H)$ an $(m,k)-$triple where $S$ is Abelian. The morphism $$(\lambda_{v},a_{v}^{*}):v^{\perp}\oplus H^{2}(S\times\widehat{S},\mathbb{Z})\longrightarrow H^{2}(M_{v},\mathbb{Z})$$is an isomorphism of Hodge structures.
\end{cor}

\proof The Leray spectral sequence associated with $a_{v}:M_{v}\longrightarrow S\times\widehat{S}$ gives an exact sequence $$0\longrightarrow H^{2}(S\times\widehat{S},\mathbb{Z})\stackrel{a_{v}^{*}}\longrightarrow H^{2}(M_{v},\mathbb{Z})\stackrel{\phi}\longrightarrow H^{0}(S\times\widehat{S},R^{2}a_{v*}\mathbb{Z}).$$The fibers of the local system $R^{2}a_{v*}\mathbb{Z}$ are isomorphic to $H^{2}(K_{v},\mathbb{Z})$. Moreover, the restriction morphism $\iota_{v}^{*}:H^{2}(M_{v},\mathbb{Z})\longrightarrow H^{2}(K_{v},\mathbb{Z})$ is surjective: indeed by Proposition \ref{prop:b2leq} we know that $\lambda_{v}^{0}:v^{\perp}\longrightarrow H^{2}(K_{v},\mathbb{Z})$ is an isomorphism, and $\lambda_{v}^{0}=\iota_{v}^{*}\circ\lambda_{v}$, so $\iota_{v}^{*}$ is surjective: this implies that $H^{0}(S\times\widehat{S},R^{2}a_{v*}\mathbb{Z})$ identifies with $H^{2}(K_{v},\mathbb{Z})$ and the previous exact sequence becomes $$0\longrightarrow H^{2}(S\times\widehat{S},\mathbb{Z})\stackrel{a_{v}^{*}}\longrightarrow H^{2}(M_{v},\mathbb{Z})\stackrel{\iota_{v}^{*}}\longrightarrow H^{2}(K_{v},\mathbb{Z})\longrightarrow 0.$$As $\lambda^{0}_{v}=\iota_{v}^{*}\circ\lambda_{v}$, by Proposition \ref{prop:b2leq} the statement follows.\endproof

Proposition \ref{prop:b2leq} proves part of Theorem \ref{thm:main}: what is left to prove is that the morphism $\lambda_{v}$ (resp. $\lambda_{v}^{0}$) is an isometry, and this is the content of the next section.

\section{The morphisms $\lambda_{v}$ and $\lambda^{0}_{v}$ are isometries}

If $(S,v,H)$ is an $(m,k)-$triple, on $v^{\perp}$ we have a pure weight two Hodge structure induced by the one we have on $\widetilde{H}(S,\mathbb{Z})$ and a lattice structure with respect to the Mukai pairing $(\cdot,\cdot)$. Moreover, as seen in section 2.4 if $S$ is K3 on $H^{2}(M_{v},\mathbb{Z})$ we have a pure weight two Hodge structure and a lattice structure with respect to the Beauville-Namikawa form $b_{v}$, and if $S$ is Abelian the same holds for $H^{2}(K_{v},\mathbb{Z})$. We will denote by $B_{v}$ the bilinear form associated with the quadratic form $b_{v}$.

The aim of this section is to complete the proof of Theorem \ref{thm:main}, i. e. that if $(S,v,H)$ is an $(m,k)-$triple then $\lambda_{v}:v^{\perp}\longrightarrow H^{2}(M_{v},\mathbb{Z})$ is a Hodge isomorphism of $\mathbb{Z}-$modules which is an isometry if $S$ is K3, and $\lambda_{v}^{0}:v^{\perp}\longrightarrow H^{2}(K_{v},\mathbb{Z})$ is a Hodge isomorphism of $\mathbb{Z}-$modules which is an isometry if $S$ is Abelian and $(m,k)\neq(1,1),(1,2)$.

By Proposition \ref{prop:b2leq} we just need to prove that $\lambda_{v}$ is an isometry if $S$ is K3, and that $\lambda_{v}^{0}$ is an isometry if $S$ is Abelian and $(m,k)\neq(1,1),(1,2)$. To prove this we will need the following result, which is known to hold for irreducible symplectic manifolds:

\begin{lem}
\label{lem:huyb}Let $X$ be a $\mathbb{Q}-$factorial irreducible symplectic variety of dimension $2n$ and whose singular locus is of codimension at least 4. Let $SH^{2}(X,\mathbb{C})$ be the subalgebra of $H^{*}(X,\mathbb{C})$ spanned by $H^{2}(X,\mathbb{C})$. Then $$SH^{2}(X,\mathbb{C})\simeq\frac{S^{*}H^{2}(X,\mathbb{C})}{\langle\alpha^{n+1}\,|\,\alpha\in H^{2}(X,\mathbb{C}),\,\,q_{X}(\alpha)=0\rangle},$$where $S^{*}H^{2}(X,\mathbb{C})$ is the symmetric algebra of $H^{2}(X,\mathbb{C})$, and $q_{X}$ is the Beauville-Namikawa form on $H^{2}(X,\mathbb{C})$.
\end{lem} 

\proof By Proposition 1.10 of \cite{PR3} we know that if $X$ is an irreducible symplectic variety then $H^{1}(X,\mathscr{O}_{X})=0$ and $h^{0}(X^{s},\Omega_{X^{s}}^{2})=1$. If moreover $X$ is $\mathbb{Q}-$factorial and its singular locus has codimension at least 4, then by point (3) of Theorem 8 of \cite{N1} the period map $p_{X}:Def(X)\longrightarrow Q_{X}$ is an open embedding, where $Def(X)$ is the base of a Kuranishi family of $X$ and $Q_{X}\subseteq\mathbb{P}(H^{2}(X,\mathbb{C}))$ is the quadric associated with the quadratic form naturally induced on $H^{2}(X,\mathbb{C})$ by the Beauville-Namikawa form on $H^{2}(X,\mathbb{Z})$. One may then follow the proof of Proposition 24.1 of \cite{GHJ} to get the statement.\endproof

We are now in the position to complete the proof of the following:

\begin{prop}
\label{prop:hodgeisoproof}
Let $(S,v,H)$ be an $(m,k)-$triple, and if $S$ is Abelian suppose that $(m,k)\neq(1,1)$. 
\begin{enumerate}
 \item If $S$ is K3, the morphism $\lambda_{v}:v^{\perp}\longrightarrow H^{2}(M_{v},\mathbb{Z})$ is an isometry.
 \item If $S$ is Abelian, the morphism $\lambda^{0}_{v}:v^{\perp}\longrightarrow H^{2}(K_{v},\mathbb{Z})$ is an isometry.
\end{enumerate}
\end{prop}

\proof The proof in the case of K3 surfaces and in the case of Abelian surfaces follows the same pattern, but there are important differencies. We will start by proving the statement for K3 surfaces, and then we proceed with Abelian surfaces.

\textbf{The case of K3 surfaces}. Suppose that $S$ is K3, we will proceed by induction on $m$. As, by \cite{Y1} and of \cite{PR2} the statement holds for $m=1$ and $(m,k)=(2,1)$, we suppose that it holds for $m-1$, and we show that it holds for $m$. We need to show that if $\alpha,\beta\in v^{\perp}$, then $B_{v}(\lambda_{v}(\alpha),\lambda_{v}(\beta))=(\alpha,\beta)$.

In the following section we will loosely use the same notation for a quadratic form (or a bilinear form) on a $\mathbb{Z}-$module $\Lambda$ and the one naturally induced by it on $\Lambda\otimes\mathbb{C}$. Similarily, we use the same notation for a morphism between two $\mathbb{Z}-$modules $\Lambda$ and $\Lambda'$ and the one induced by it between $\Lambda\otimes\mathbb{C}$ and $\Lambda'\otimes\mathbb{C}$. 

Since the Beauville-Namikawa pairing $B_{v}$ and the Mukai pairing $(\cdot,\cdot)$ are integral, primitive bilinear forms of signature $(3,20)$, if $Q_{M}\subseteq H^{2}(M_{v},\mathbb{C})$ and $Q_{v}\subseteq v^{\perp}\otimes\mathbb{C}$ are the quadrics defined by $B_{v}$ and $(\cdot,\cdot)$ respectively, it will be enough to prove that $Q_{M}=\lambda_{v}(Q_{v})$. Moreover, as $Q_{M}$ and $Q_{v}$ are irreducible and have the same dimension, it will be enough to show that $Q_{M}\subseteq\lambda_{v}(Q_{v})$, i. e. that for every $\alpha\in v^{\perp}\otimes\mathbb{C}$ if $b_{v}(\lambda_{v}(\alpha))=0$ then $(\alpha,\alpha)=0$.

Hence, let $\alpha\in v^{\perp}\otimes\mathbb{C}$ and suppose that $b_{v}(\lambda_{v}(\alpha))=0$. Since $M_{v}$ is a $\mathbb{Q}-$factorial irreducible symplectic variety of dimension $2d$ (where $d=m^{2}k+1$) whose singular locus has codimension at least 4, by Lemma \ref{lem:huyb} we have $\lambda_{v}(\alpha)^{d+1}=0$ in the cohomology of $M_{v}$.

Now, let us consider the morphism $f_{m-1}:M_{(m-1)w}\times M_{w}\longrightarrow M_{v}$ constructed in section 4.2. As $\lambda_{v}(\alpha)^{d+1}=0$, it follows that $$(f_{m-1}^{*}(\lambda_{v}(\alpha)))^{d+1}=f_{m-1}^{*}(\lambda_{v}(\alpha)^{d+1})=0.$$But since by point (2) of Proposition \ref{prop:diag1} we have $$f_{m-1}^{*}(\lambda_{v}(\alpha))=\pi_{1}^{*}(\lambda_{(m-1)w}(\alpha))+\pi_{2}^{*}(\lambda_{w}(\alpha)),$$where $\pi_{1}:M_{(m-1)w}\times M_{w}\longrightarrow M_{(m-1)w}$ and $\pi_{2}:M_{(m-1)w}\times M_{w}\longrightarrow M_{w}$ are the two projections, it follows that 
\begin{equation}
\label{equation:sum}
\sum_{j=0}^{d+1}\binom{d+1}{j}(\pi_{1}^{*}\lambda_{(m-1)w}(\alpha))^{j}\wedge(\pi_{2}^{*}\lambda_{w}(\alpha))^{d+1-j}=0.
\end{equation}

Now, the K\"unneth decomposition of $H^{2d+2}(M_{(m-1)w}\times M_{w})$ is $$H^{2d+2}(M_{(m-1)w}\times M_{w})=\oplus_{p=0}^{2d+2}H^{p}(M_{(m-1)w})\otimes H^{2d+2-p}(M_{w}),$$and for every $j\in\{0,\cdots,d+1\}$ we have $$(\pi_{1}^{*}\lambda_{(m-1)w}(\alpha))^{j}\wedge(\pi_{2}^{*}\lambda_{w}(\alpha))^{d+1-j}\in H^{2j}(M_{(m-1)w})\otimes H^{2d+2-2j}(M_{w}).$$As a consequence, equation (\ref{equation:sum}) gives that $$(\pi_{1}^{*}\lambda_{(m-1)w}(\alpha))^{j}\wedge(\pi_{2}^{*}\lambda_{w}(\alpha))^{d+1-j}=0$$for every $0\leq j\leq d+1$. 

It then follows that for every $0\leq j\leq d+1$ we have either $(\pi_{1}^{*}\lambda_{(m-1)w}(\alpha))^{j}=0$ or $(\pi_{2}^{*}\lambda_{w}(\alpha))^{d+1-j}=0$, and so that either $\lambda_{(m-1)w}(\alpha)^{j}=0$ or that $\lambda_{w}(\alpha)^{d+1-j}=0$.

Let $2d':=\dim(M_{(m-1)w})$ and $j\leq 2d'$: if $\lambda_{(m-1)w}(\alpha)^{j}=0$, then $\lambda_{(m-1)w}(\alpha)^{2d'}=0$, so $b_{(m-1)w}(\lambda_{(m-1)w}(\alpha))=0$ by the Fujiki relations. Similarly, if $2d'':=\dim(M_{w})$ and $d+1-j\leq 2d''$, if $\lambda_{w}(\alpha)^{d+1-j}=0$, then $b_{w}(\lambda_{w}(\alpha))=0$ by the Fujiki relations. 

We deduce that if there is a positive integer $j$ such that $j\leq 2d'$ and $d+1-j\leq 2d''$, then either $b_{(m-1)w}(\lambda_{(m-1)w}(\alpha))=0$ or $b_{w}(\lambda_{w}(\alpha))=0$: by induction hypothesis, this implies $(\alpha,\alpha)=0$, concluding the proof. We are then left with proving that an integer $j$ such that $j\leq 2d'$ and $d+1-j\leq 2d''$ exists.

To do so, notice that $j\leq 2d'$ if and only if $j\leq 2k(m-1)^{2}+2$, and $d+1-j\leq 2d''$ if and only if $j\geq 2k(m^{2}/2-1)$, so we just need to prove that there is an integer $j$ such that $$2k\bigg(\frac{m^{2}}{2}-1\bigg)\leq j\leq 2k(m-1)^{2}+2.$$But notice that the difference between the right and the left hand side of this inequality is $$2k(m-1)^{2}+2-2k\bigg(\frac{m^{2}}{2}-1\bigg)=k(m-2)^{2}+2\geq 2$$for every $m$: the existence of $j$ is proved, and we are done with point (1).

\textbf{The case of Abelian surfaces}. If $S$ in an Abelian surface, recall that $\lambda_{v}^{0}=\iota_{v}^{*}\circ\lambda_{v}$ where $\iota_{v}:K_{v}\longrightarrow M_{v}$ is the inclusion. As before, we just need to prove that for every $\alpha\in v^{\perp}\otimes\mathbb{C}$, if $b_{v}(\lambda_{v}^{0}(\alpha))=0$, then $(\alpha,\alpha)=0$.

The proof will rest on the following Lemma, which provides an equivalent condition to $b_{v}(\lambda_{v}^{0}(\alpha))=0$, and a necessary condition for the same equality:

\begin{lem}
\label{lem:intqvk} Let $S$ be an Abelian surface, let $(S,v,H)$ be an $(m,k)$-triple, asssume that $(m,k)\ne(1,1)$ and
let $\alpha\in v^{\perp}\otimes\mathbb{C}$.
\begin{enumerate}
 \item We have $b_{v}(\lambda^{0}_{v}(\alpha))=0$ if and only if for every $\beta\in H^{2}(S\times\widehat{S},\mathbb{C})$ we have $$\int_{M_{v}}(\lambda_{v}(\alpha)+a_{v}^{*}(\beta))^{2m^{2}k+2}=0.$$
 \item If $b_{v}(\lambda_{v}^{0}(\alpha))=0$, then for every $\beta\in H^{2}(S\times\widehat{S},\mathbb{C})$ we have $$(\lambda_{v}(\alpha)+a_{v}^{*}\beta)^{m^{2}k+4}=0.$$
\end{enumerate}
\end{lem}

\proof We start by proving point (1), and let $s:=2m^{2}k+2$. Recall that we have an \'etale cover $$\tau_{v}:K_{v}\times S\times\widehat{S}\longrightarrow M_{v},\,\,\,\,\,\,\,\,\tau_{v}(\mathcal{E},p,L):=\tau_{p}^{*}\mathcal{E}\otimes L,$$where $\tau_{p}$ is the translation by $p$ on $S$, which fits into a commutative diagram
\begin{equation}
\begin{CD}
K_{v}\times S\times\widehat{S} @>{\tau_{v}}>> M_{v}\\
@V{p_{2}}VV                   @VV{a_{v}}V\\
S\times\widehat{S} @>>{\rho}> S\times\widehat{S}
\end{CD}
\end{equation}
where $p_{2}$ is the projection and $\rho$ is the multiplication by a nonzero scalar $d'$. The morphism $\tau_{v}$ induces an isomorphism $$\tau_{v}^{*}:H^{2}(M_{v},\mathbb{C})\longrightarrow H^{2}(K_{v}\times S\times\widehat{S},\mathbb{C}).$$

If $d$ is the degree of $\tau_{v}$, then for every $\beta\in H^{2}(S\times\widehat{S},\mathbb{C})$ we have $$\int_{K_{v}\times S\times\widehat{S}}(\tau_{v}^{*}\lambda_{v}(\alpha)+\tau_{v}^{*}a_{v}^{*}\beta)^{s}=d\int_{M_{v}}(\lambda_{v}(\alpha)+a_{v}^{*}\beta)^{s}.$$Hence the right hand side of this equality is zero for every $\beta\in H^{2}(S\times\widehat{S},\mathbb{C})$ if and only if the left hand side is zero for every such $\beta$.

Now, the K\"unneth decomposition gives $$H^{2}(K_{v}\times S\times\widehat{S})\simeq(H^{2}(K_{v})\otimes H^{0}(S\times\widehat{S}))\oplus(H^{0}(K_{v})\otimes H^{2}(S\times\widehat{S}))$$since $K_{v}$ is an irreducible symplectic manifold. If we let $p_{1}$ and $p_2$ be the projections of $K_{v}\times S\times\widehat{S}$ to $K_{v}$ and $S\times\widehat{S}$ respectively, we then have that $\tau_{v}^{*}\lambda_{v}(\alpha)=p_{1}^{*}\lambda_{v}^{0}(\alpha)+p_{2}^{*}\gamma$ for $\kappa\in H^{2}(K_{v},\mathbb{C})$ and $\gamma\in H^{2}(S\times\widehat{S},\mathbb{C})$, and $\tau_{v}^{*}a_{v}^{*}\beta=p_{2}^{*}(d'\beta)$.

As a consequence we have 
\begin{equation}
\label{eq:tauv}
\tau_{v}^{*}\lambda_{v}(\alpha)+\tau_{v}^{*}a_{v}^{*}\beta=p_{1}^{*}\lambda_{v}^{0}(\alpha)+p_{2}^{*}(\gamma+d'\beta).
\end{equation}
Hence $$(\tau_{v}^{*}\lambda_{v}(\alpha)+\tau_{v}^{*}a_{v}^{*}\beta)^{s}=\sum_{j=0}^{s}\binom{s}{j}(p_{1}^{*}\lambda_{v}^{0}(\alpha))^{j}\wedge(p_{2}^{*}(\gamma+d'\beta))^{s-j}.$$Since $\dim(K_{v})=2m^{2}k-2$ and $\dim(S\times\widehat{S})=4$, we see that $$\sum_{j=0}^{s}\binom{s}{j}(p_{1}^{*}\lambda_{v}^{0}(\alpha))^{j}\wedge(p_{2}^{*}(\gamma+d'\beta))^{s-j}=\binom{s}{s-4}(p_{1}^{*}\kappa)^{s-4}\wedge(p_{2}^{*}(\gamma+d'\beta))^{4}.$$

We thus get $$\int_{K_{v}\times S\times\widehat{S}}(\tau_{v}^{*}\lambda_{v}(\alpha)+\tau_{v}^{*}a_{v}^{*}\beta)^{s}=$$ $$=\binom{s}{s-4}\bigg(\int_{S\times\widehat{S}}(\gamma+d'\beta)^{4}\bigg)\cdot\bigg(\int_{K_{v}}\lambda_{v}^{0}(\alpha)^{2m^{2}k-2}\bigg).$$As a consequence, the left hand side in this equality is 0 for every $\beta\in H^{2}(S\times\widehat{S},\mathbb{C})$ if and only if $$\int_{K_{v}}\lambda_{v}^{0}(\alpha)^{2m^{2}k-2}=0.$$But $$\int_{K_{v}}\lambda_{v}^{0}(\alpha)^{2m^{2}k-2}=D_{v}b_{v}(\lambda_{v}^{0}(\alpha))^{m^{2}k-1},$$where $D_{v}$ is the Fujiki constant of $K_{v}$, concluding the proof of point (1).

For point (2), suppose that $b_{v}(\lambda_{v}^{0}(\alpha))=0$. As $\tau_{v}:K_{v}\times S\times\widehat{S}\longrightarrow M_{v}$ is \'etale, we have that $$\tau_{v}^{*}:H^{n}(M_{v},\mathbb{C})\longrightarrow H^{n}(K_{v}\times S\times\widehat{S},\mathbb{C})$$is injective for every $n$, hence in order to prove that $(\lambda_{v}(\alpha)+a_{v}^{*}\beta)^{m^{2}k+4}=0$ we just need to prove that $$(\tau_{v}^{*}(\lambda_{v}(\alpha)+a_{v}^{*}\beta))^{m^{2}k+4}=0.$$By equation (\ref{eq:tauv}), letting $r:=m^{2}k+4$ it then follows that $$(\tau_{v}^{*}(\lambda_{v}(\alpha)+a_{v}^{*}\beta))^{r}=\sum_{j=0}^{r}\binom{r}{j}(p_{1}^{*}\lambda_{v}^{0}(\alpha))^{j}\wedge(p_{2}^{*}(\gamma+d'\beta))^{r-j}.$$ 

Now, as $K_{v}$ is a $\mathbb{Q}-$factorial irreducible symplectic variety of dimension $2m^{2}k-2$ and as $b_{v}(\lambda^{0}_{v}(\alpha))=0$, by Lemma \ref{lem:huyb} we have that $\lambda_{v}^{0}(\alpha)^{m^{2}k}=0$. As $m^{2}k=r-4$, it follows that for every $j\geq r-4=m^2k$ we have $(p_{1}^{*}\lambda_{v}^{0}(\alpha))^{j}=0$. Moreover, as $S\times\widehat{S}$ has dimension 4, for every $j<r-4$ we have $r-j>4$ and $(p_{2}^{*}(\gamma+d\beta))^{r-j}=0$. It follows that $$\sum_{j=0}^{r}\binom{r}{j}(p_{1}^{*}\lambda_{v}^{0}(\alpha))^{j}\wedge(p_{2}^{*}(\gamma+d'\beta))^{r-j}=0,$$ because every summand has a trivial factor. This concludes the proof.\endproof

The proof of (2) of Proposition \ref{prop:hodgeisoproof} is again by induction on $m$ for fixed $k$. The base of the induction is $m=1$ for $k>1$ and $m=2$ for $k=1$: in these cases the statement holds by \cite{Y2} and \cite{PR2} .

Suppose now that $\alpha\in v^{\perp}\otimes\mathbb{C}$ is such that $b_{v}(\lambda_{v}^{0}(\alpha))=0$. By point (2) of Lemma \ref{lem:intqvk} we then have $(\lambda_{v}(\alpha)+a_{v}^{*}\beta)^{m^{2}k+4}=0$ for every $\beta\in H^{2}(S\times\widehat{S},\mathbb{C})$. Now, consider the morphism $$f_{m-1}:M_{(m-1)w}\times M_{w}\longrightarrow M_{v}$$constructed in section 4.2. For every $\beta\in H^{2}(S\times\widehat{S},\mathbb{C})$ we have 
\begin{equation}\label{a1}0=f_{m-1}^{*}(\lambda_{v}(\alpha)+a_{v}^{*}\beta)^{m^{2}k+4}=(f_{m-1}^{*}\lambda_{v}(\alpha)+f_{m-1}^{*}a_{v}^{*}\beta)^{m^{2}k+4}=\end{equation} $$=(\pi_{1}^{*}\lambda_{(m-1)w}(\alpha)+\pi_{2}^{*}\lambda_{w}(\alpha)+f_{m-1}^{*}a_{v}^{*}\beta)^{m^{2}k+4},$$where the last equality is a consequence of point (2) of Proposition \ref{prop:diag1}. By definition of $a_v$ (see Section 4.1 of \cite{Y2}) there are nonzero integral constants $c'$ and $c''$ such that $$a_{v}\circ f_{m-1}=(c'a_{(m-1)w}\circ\pi_{1}+c''a_{w}\circ\pi_{2})$$ and letting $r:=m^{2}k+4$, as a consequence of (\ref{a1}) we see that $$(\pi_{1}^{*}(\lambda_{(m-1)w}(\alpha)+c'a_{(m-1)w}^{*}\beta)+\pi_{2}^{*}(\lambda_{w}(\alpha)+c''a_{w}^{*}\beta))^{r}=0,$$so that $$\sum_{j=0}^{r}\binom{r}{j}(\pi_{1}^{*}(\lambda_{(m-1)w}(\alpha)+c'a_{(m-1)w}^{*}\beta))^{j}\wedge(\pi_{2}^{*}(\lambda_{w}(\alpha)+c''a_{w}^{*}\beta))^{r-j}=0.$$As in the case of K3 surfaces, this implies that for every $0\leq j\leq r$ we have $$(\pi_{1}^{*}(\lambda_{(m-1)w}(\alpha)+c'a_{(m-1)w}^{*}\beta))^{j}\wedge(\pi_{2}^{*}(\lambda_{w}(\alpha)+c''a_{w}^{*}\beta))^{r-j}=0,$$so either $(\lambda_{(m-1)w}(\alpha)+c'a_{(m-1)w}^{*}\beta)^{j}=0$ or $(\lambda_{w}(\alpha)+c''a_{w}^{*}\beta)^{r-j}=0$.

Now, notice that the subsets $$V_{j}:=\{\beta\in H^{2}(S\times\widehat{S},\mathbb{C})\,|\,(\lambda_{(m-1)w}(\alpha)+c'a_{(m-1)w}^{*}\beta)^{j}=0\}$$and $$W_{j}:=\{\beta\in H^{2}(S\times\widehat{S},\mathbb{C})\,|\,(\lambda_{w}(\alpha)+c''a_{w}^{*}\beta)^{r-j}=0\}$$are two affine subvarieties of $H^{2}(S\times\widehat{S},\mathbb{C})$, and by the previous discussion for every $0\leq j\leq r$ we have $V_{j}\cup W_{j}=H^{2}(S\times\widehat{S},\mathbb{C})$. 

It follows that either $V_{j}=H^{2}(S\times\widehat{S},\mathbb{C})$ or $W_{j}=H^{2}(S\times\widehat{S},\mathbb{C})$, so for every $0\leq j\leq r$ we have either that $(\lambda_{(m-1)w}(\alpha)+c'a_{(m-1)w}^{*}\beta)^{j}=0$ for every $\beta\in H^{2}(S\times\widehat{S},\mathbb{C})$, or that $(\lambda_{w}(\alpha)+c''a_{w}^{*}\beta)^{r-j}=0$ for every $\beta\in H^{2}(S\times\widehat{S},\mathbb{C})$.

Let now $2d'$ be the dimension of $M_{(m-1)w}$ and $2d''$ the dimension of $M_{w}$. If there is $j\leq 2d'$ such that we have $(\lambda_{(m-1)w}(\alpha)+c'a_{(m-1)w}^{*}\beta)^{j}=0$ for every $\beta\in H^{2}(S\times\widehat{S},\mathbb{C})$, then $(\lambda_{(m-1)w}(\alpha)+c'a_{(m-1)w}^{*}\beta)^{2d'}=0$, so Lemma \ref{lem:intqvk} implies $b_{(m-1)w}(\lambda_{(m-1)w}^{0}(\alpha))=0$. 

Similarly, if $j$ is such that $m^{2}k+4-j\leq 2d''$ and for every $\beta\in H^{2}(S\times\widehat{S},\mathbb{C})$ we have $(\lambda_{w}(\alpha)+c''a_{w}^{*}\beta)^{m^{2}k+4-j}=0$, then $(\lambda_{w}(\alpha)+c''a_{w}^{*}\beta)^{2d''}=0$, hence by Lemma \ref{lem:intqvk} we get that $b_{w}(\lambda_{w}^{0}(\alpha))=0$. 

As a consequence, the existence of an integer $j$ such that $j\leq 2d'$ and $m^{2}k+4-j\leq 2d''$ implies that either $b_{(m-1)w}(\lambda_{(m-1)w}^{0}(\alpha))=0$ or else $b_{w}(\lambda_{w}^{0}(\alpha))=0$. To prove the existence of such a $j$ we notice that $j\leq 2d'$ if and only if $j\leq 2k(m-1)^{2}+2$, and $m^{2}k+4-j\leq 2d''$ if and only if $j\geq 2k(m^{2}/2-1)+2$, so we just need to prove that there is an integer $j$ such that $$2k\bigg(\frac{m^{2}}{2}-1\bigg)+2\leq j\leq 2k(m-1)^{2}+2.$$ Since the difference between the right and the left hand side of this inequality is $$2k(m-1)^{2}-2k\bigg(\frac{m^{2}}{2}-1\bigg)=k(m-2)^{2}\geq 0$$for every $m$, the desired $j$ exists. 

If $k>1$, since the induction starts from $m=1$, by inductive hypothesis, both the conditions $b_{(m-1)w}(\lambda_{(m-1)w}^{0}(\alpha))=0$ and $b_{w}(\lambda_{w}^{0}(\alpha))=0$ imply $(\alpha,\alpha)=0$ and we are done.

If $k=1$, we need to have $m\geq 2$, and we know that the result holds for $m=2$, so we may suppose $m\geq 3$. The same procedure as before tells us that either $b_{(m-1)w}(\lambda^{0}_{(m-1)v}(\alpha))=$ or that $b_{w}(\lambda^{0}_{w}(\alpha))=0$. As $m\geq 3$ we then have $m-1\geq 2$ and by induction we conclude that $(\alpha,\alpha)=0$. However, since $K_w$ is a point in this case, the condition $b_{w}(\lambda_{w}^{0}(\alpha))=0$ is empty. 

Anyway we may proceed in this way. By what we proved before, for every $0\leq j\leq r$ we have either that $(\lambda_{(m-1)w}(\alpha)+c'a_{(m-1)w}^{*}\beta)^{j}=0$ for every $\beta\in H^{2}(S\times\widehat{S},\mathbb{C})$, or that $(\lambda_{w}(\alpha)+c''a_{w}^{*}\beta)^{r-j}=0$ for every $\beta\in H^{2}(S\times\widehat{S},\mathbb{C})$. 

Suppose now that $r-j\leq 4$. Since the class $(\lambda_{w}(\alpha)+c''a_{w}^{*}\beta)^{r-j}$ cannot be zero for every $\beta\in H^{2}(S\times\widehat{S},\mathbb{C})$, it follows that if $r-j\leq 4$ we need to have $(\lambda_{(m-1)w}(\alpha)+c'a_{(m-1)w}^{*}\beta)^{j}=0$ for every $\beta\in H^{2}(S\times\widehat{S},\mathbb{C})$. 

We now remark that as $r=m^{2}k+4$ and $k=1$, we have $r-j\leq 4$ if and only if $j\geq m^{2}$. As a consequence, letting $2d'$ be the dimension of $M_{(m-1)w}$, if there is $m^{2}\leq j\leq 2d'$ we may conclude as before thanks to point (1) of Lemma \ref{lem:intqvk}, obtaining $b_{(m-1)w}(\lambda_{(m-1)w}^{0}(\alpha))=0$, and hence $(\alpha,\alpha)=0$ as already remarked before.

In conclusion, we just need to prove that there is an integer $j$ such that $m^{2}\leq j\leq 2d'$. As $k=1$ we have $2d'=2(m-1)^{2}+2$, and the difference $$2(m-1)^{2}+2-m^{2}=(m-2)^{2}\geq 0,$$so the existence of such a $j$ is proved, concluding the proof.\endproof

Theorem \ref{thm:main} is now a consequence of Propositions \ref{prop:b2leq} and \ref{prop:hodgeisoproof}. As an immediate corollary of Theorem \ref{thm:main} we get the following about the analitically locally trivial deformation classes of the irreducible symplectic varieties we get from $(m,k)-$triples.

\begin{cor}
\label{cor:defomk}
For $i=1,2$ let $(S_{i},v_{i},H_{i})$ be an $(m_{i},k_{i})-$triple. 
\begin{enumerate}
 \item If $S_{1}$ and $S_{2}$ are both K3 surfaces, the moduli spaces $M_{v_{1}}(S_{1},H_{1})$ and $M_{v_{2}}(S_{2},H_{2})$ are analytically locally trivially deformation equivalent if and only if $(m_{1},k_{1})=(m_{2},k_{2})$.
 \item If $S_{1}$ and $S_{2}$ are both Abelian surfaces, the moduli spaces $K_{v_{1}}(S_{1},H_{1})$ and $K_{v_{2}}(S_{2},H_{2})$ are analytically locally trivially deformation equivalent if and only if $(m_{1},k_{1})=(m_{2},k_{2})$.
 \item If $S_{1}$ is K3 and $S_{2}$ is Abelian, the moduli spaces $M_{v_{1}}(S_{1},H_{1})$ and $K_{v_{2}}(S_{2},H_{2})$ are not deformation equivalent.
\end{enumerate}
\end{cor} 

\proof The last point of the statement is clear since the two Betti numbers of $M_{v_{1}}$ and $K_{v_{2}}$ are different by Theorem \ref{thm:b22}.

For the remaining part of the statement, we present a proof only for the case where $S_{1}$ and $S_{2}$ are both K3 surfaces, the other being similar. If $(m_{1},k_{1})=(m_{2},k_{2})$, by point (1) of Theorem \ref{thm:mio} the moduli spaces are deformation equivalent.

Conversely, suppose that $M_{v_{1}}(S_{1},H_{1})$ and $M_{v_{2}}(S_{2},H_{2})$ are deformation equivalent. As the Beauville form is deformation invariant it follows that $H^{2}(M_{v_{1}},\mathbb{Z})$ and $H^{2}(M_{v_{2}},\mathbb{Z})$ are isometric. By Theorem \ref{thm:main} it follows that $v_{1}^{\perp}$ and $v_{2}^{\perp}$ have to be isometric, and hence they have the same discriminant. 

As $v_{i}^{\perp}=w_{i}^{\perp}$, the discriminant of $v_{i}^{\perp}$ is $w_{i}^{2}=2k_{i}$: this implies $k_{1}=k_{2}$. As the moduli spaces have the same dimension (since they are deformation equivalent), we finally get $m_{1}=m_{2}$.\endproof

The deformation class of the moduli space $M_{v}(S,H)$ associated with an $(m,k)-$triple $(S,v,H)$ is then completely determined by $m$ and $k$. Notice that if $d\in\mathbb{N}$, $d>2$, then there can be different pairs $(m,k)\in\mathbb{N}^{2}$ such that $2m^{2}k\pm 2=2d$: by Corollary \ref{cor:defomk} each of these pairs corresponds to a unique deformation class of moduli spaces of sheaves of the same dimension.

\section{The Fujiki constants of $M_{v}$ and $K_{v}$}

The aim of this section is to calculate the Fujiki constant of the irreducible symplectic varieties $M_{v}$ and $K_{v}$. We recall that by \cite{S}, if $X$ is an $2n-$dimensional irreducible symplectic variety, then there is a positive real number $C_{X}$, called the Fujiki constant of $X$, such that $$C_{X}b_{X}(\alpha)^{n}=\int_{X}\alpha^{2n}$$for every $\alpha\in H^{2}(X,\mathbb{R})$, where $b_{X}$ is the Beauville-Namikawa form of $X$. 

As if $(S,v,H)$ is an $(m,k)-$triple then $M_{v}(S,H)$ (if $S$ is K3) and $K_{v}(S,H)$ (if $S$ is Abelian and $(m,k)\neq(1,1)$) are irreducible symplectic varieties, they have a Fujiki constant. We will let $C_{v}$ be the Fujiki constant of $M_{v}(S,H)$ (if $S$ is K3), and $D_{v}$ be the Fujiki constant of $K_{v}(S,H)$ (if $S$ is Abelian). 

For $(1,k)-$triples we have $C_{v}=\frac{(2k+2)!}{(k+1)2^{k+1}}$ and $D_{v}=\frac{(2k-2)!}{(k-1)2^{k-1}}$ (see \cite{B}), while for $(2,1)-$triples it is $C_{v}=945=\frac{10!}{5!2^{5}}$ and $D_{v}=15=\frac{6!}{3!2^{3}}$ (see \cite{R1}). 

In this section we prove Theorem \ref{thm:fujiki}, where we calculate the Fujiki constant for all other $(m,k)$-triples. Theorem \ref{thm:fujiki} is an immediate consequence of the following Proposition:

\begin{prop}
\label{prop:cucv1}Let $(S,v,H)$ be an $(m,k)-$triple and $(T,u,L)$ a $(1,m^{2}k)-$triple. 
\begin{enumerate}
 \item If $S$ is K3 we have $C_{v}=C_{u}$.
 \item If $S$ is Abelian and $(m,k)\neq(1,1)$, we have $D_{v}=D_{u}$.
\end{enumerate}
\end{prop}

We describe the idea of the proof in the case of K3 surfaces, the one for Abelian surfaces being similar. 

First, recall that the Fujiki constant is deformation invariant, and that for two $(m,k)-$triples the corresponding moduli spaces are deformation equivalent by \cite{PR3}. We then reduce to the case of an $(m,k)-$triple $(S,v,H)$ where $H^{2}=2k$, the N\'eron-Severi group of $S$ is generated by the first Chern class $h$ of $H$ and $v=m(0,h,0)$. As in the computation of the second Betti number we consider the Mukai vector $u=(0,mh,1-m^{2}k)$ on $S$ and relate the Fujiki constants of the $(m,k)-$triple $(S,v,H)$ and of the $(1,m^{2}k)-$triple $(S,u,H)$.

In order to prove Proposition \ref{prop:cucv1} for K3 surfaces, we will make use of the two Hodge isometries $$\lambda_{u}:u^{\perp}\longrightarrow H^{2}(M_{u},\mathbb{Z}),\,\,\,\,\,\,\,\,\lambda_{v}:v^{\perp}\longrightarrow H^{2}(M_{v},\mathbb{Z}).$$If we let $\sigma\in H^{2,0}(S)$ be a nontrivial class, observe that $\alpha:=(0,\sigma+\overline{\sigma},0)\in u^{\perp}\cap v^{\perp}$, as $u,v\in\widetilde{H}^{1,1}(S)$, $\sigma+\overline{\sigma}\in\widetilde{H}^{2,0}(S)\oplus\widetilde{H}^{0,2}(S)$, and $\widetilde{H}^{1,1}(S)$ is orthogonal to $\widetilde{H}^{2,0}(S)\oplus\widetilde{H}^{0,2}(S)$ with respect to the Mukai form on $\widetilde{H}(S,\mathbb{C})$. 

As a consequence $\lambda_{u}(\alpha)\in H^{2}(M_{u},\mathbb{Z})$ and $\lambda_{v}(\alpha)\in H^{2}(M_{v},\mathbb{Z})$, and if $2n$ is the dimension of $M_{u}$ and $M_{v}$, then we have $$\int_{M_{u}}\lambda_{u}(\alpha)^{2n}=C_{u}b_{u}(\lambda_{u}(\alpha))^{n}=C_{u}(\alpha,\alpha)^{n},$$and $$\int_{M_{v}}\lambda_{v}(\alpha)^{2n}=C_{v}b_{v}(\lambda_{u}(\alpha))^{n}=C_{v}(\alpha,\alpha)^{n}.$$Since $(\alpha,\alpha)\neq 0$, to show that $C_{u}=C_{v}$ it will be enough to show that $$\int_{M_{u}}\lambda_{u}(\alpha)^{2n}=\int_{M_{v}}\lambda_{v}(\alpha)^{2n}.$$

Since for a coherent sheaf $\mathscr{F}$ of Mukai vector $(0,\xi,a)$ we have that $a=\chi(\mathscr{F})$, the general point of $M_{v}$ represents the push-forward to $S$ of a degree $km^2$ line bundle on a smooth curve $C$ of genus $g=km^2+1$ in the linear system $|mH|$. Analogously, the general point of $M_{u}$ represents the push-forward of a degree $1$ line bundle on a smooth curve $C$ of $|mH|$. As in Proposition \ref{prop:psi}, the $km^2$ tensor power gives a dominant, generically finite, rational map $\psi:M_{u}\dashrightarrow M_{v}$ of degree $(g-1)^{2g}$. 

Letting $\psi^*:H^{2}(M_{v},\mathbb{Z})\rightarrow H^{2}(M_{u},\mathbb{Z})$  be the morphism induced by $\psi$ in cohomology, we will show Proposition \ref{prop:cucv1} as a consequence of the equality $$\psi^{*}\lambda_{v}(\alpha)=km^2\lambda_{u}(\alpha).$$ 

This equality will be proved in Proposition \ref{prop:relation} by constructing a smooth surface $Y$ together with two families of stable sheaves of Mukai vector $u$ and $v$ parameterized by $Y$: the families induce modular morphisms from $Y$ to $M_{u}$ and to $M_{v}$ and the equality will be obtained by comparing the pull-backs of $\lambda_{u}(\alpha)$ and $\lambda_{v}(\alpha)$ to $Y$ under these modular maps. 

\subsection{The proof of Proposition \ref{prop:cucv1}}

In order to deal with the K3 case and the Abelian case simultaneously we use the notation introduced in Setting \ref{sett:notazioni}. We then let $S$ be either a projective K3 surface or an Abelian surface such that $NS(S)=\mathbb{Z}\cdot h$, where $h$ is the first Chern class in cohomology of an ample line bundle $H$ with $h^{2}=2k$, and if $S$ is Abelian we will suppose furthormore that the self-intersection $H^{2}$ of $H$ in the Chow ring of $S$ belongs to the kernel of the Albanese morphism $\sum:CH_0(S)\longrightarrow S$. Moreover, we consider the Mukai vectors $v=m(0,h,0)$ and $u=(0,mh,1-m^{2}k)$ and, as in  Setting \ref{sett:notazioni}, we set
 $$M:=\left\{\begin{array}{ll} M_{v}(S,H), & {\rm if}\,\,S\,\,{\rm is}\,\,{\rm K3}\\ K_{v}(S,H), & {\rm if}\,\,S\,\,{\rm is}\,\,{\rm Abelian}\end{array}\right.$$and $$M':=\left\{\begin{array}{ll} M_{u}(S,H), & {\rm if}\,\,S\,\,{\rm is}\,\,{\rm K3}\\ K_{u}(S,H), & {\rm if}\,\,S\,\,{\rm is}\,\,{\rm Abelian}\end{array}\right.$$

We let $P:=|mH|$, which is a projective space of dimension $km^{2}+1$, and we let $P^{0}$ be the subset of $P$ parameterizing smooth curves and $P^{1}$ the subset of $S$ parameterizing nodal curves having at most one node. As in Proposition \ref{prop:psi}, if $C\in P^{0}$ and we let $\iota:C\longrightarrow S$ be its embedding, if $L\in Pic^{1}(C)$ then $\iota_{*}L\in(M')^{s}$, while if $L\in Pic^{km^{2}}(C)$ then $\iota_{*}L\in M^{s}$.

By Remark \ref{oss:tensj} and Proposition \ref{prop:psi} we have a dominant, generically finite, rational morphism from $M'$ to $M$ which maps $\iota_{*}L$ to $\iota_{*}L^{\otimes km^{2}}$ (for every $C\in P^{1}$ whose embedding in $S$ is $\iota$, and for every $L\in Pic^1(C)$): this was the definition of the morphism $\psi:J'\longrightarrow J$, which by a slight abuse of notation we will write $\psi:M'\dashrightarrow M$. 

By smoothness of $M'$, the rational map $\psi$ induces a pull back in integral or complex cohomology $\psi^{*}:H^{2}(M)\rightarrow H^{2}(M')$ which is defined as follows: $$H^{2}(M,\mathbb{Z})\longrightarrow H^{2}(J,\mathbb{Z})\stackrel{\psi^{*}}\longrightarrow H^{2}(J',\mathbb{Z})\longrightarrow H^{2}(M',\mathbb{Z}),$$where the first and the last morphisms are the restriction morphisms (which are isomorphisms by Lemma \ref{lem:b2eb2aperti}), and the map in the middle is well-defined since $\psi:J'\longrightarrow J$ is a well-defined morphism of quasi-projective varieties.

We conclude the notation with the following: if $S$ is K3 we set $$\lambda:=\lambda_v:v^{\perp}\rightarrow H^{2}(M,\mathbb{Z}),\,\,\,\,\,\,\,\,\lambda':=\lambda_u:u^{\perp}\rightarrow H^{2}(M',\mathbb{Z}),$$and if $S$ is Abelian we set $$\lambda:=\lambda^{0}_v:v^{\perp}\rightarrow H^{2}(M,\mathbb{Z}),\,\,\,\,\,\,\,\,\lambda':=\lambda^{0}_u:u^{\perp}\rightarrow H^{2}(M',\mathbb{Z}).$$

As a first step towards the proof of Proposition \ref{prop:cucv1}, we use $\psi^{*}$ to compare holomorphic two forms on $M$ and $M'$.

\begin{prop}
\label{prop:relation}
Let $\sigma\in H^{2,0}(S)$ be a nontrivial class and consider the class $\alpha:=(0,\sigma,0)\in\widetilde{H}(S,\mathbb{C})$. Then $\alpha\in v^{\perp}\cap u^{\perp}$, and we have $$\psi^{*}\lambda(\alpha)=km^2\lambda'(\alpha).$$
\end{prop}

\proof The fact that $\alpha\in v^{\perp}\cap u^{\perp}$ is immediate since $\alpha\in\widetilde{H}^{2,0}(S)$ while $v,w\in\widetilde{H}^{1,1}(S)$, and recall that $\widetilde{H}^{2,0}(S)\oplus\widetilde{H}^{0,2}(S)$ is orthogonal to $\widetilde{H}^{1,1}(S)$ with respect to the Mukai pairing. Notice that $\alpha\in(v^{\perp})^{2,0}$ and $\alpha\neq 0$, hence by Proposition \ref{prop:lambdav} we have $\lambda'(\alpha)\in H^{2,0}(M')$ and by Lemma \ref{lem:almeno23} we have $\lambda'(\alpha)\neq 0$.
 
Let $Y$ be a smooth algebraic surface equipped with a regular morphism $f':Y\rightarrow M'$ such that the pull back to $Y$ of the symplectic form on $M'$ is not zero. Shrinking $Y$ and composing with an \'etale morphism if necessary, we may also assume that there exists a smooth family of curves $g:\mathcal{C}\longrightarrow Y$ of the linear system $|mH|$ with an inclusion $\iota:\mathcal{C}\longrightarrow S\times Y$ over $Y$ and a line bundle $\mathcal{L}$ having degree one on the fibers of $g$. This implies that $\mathcal{F}:=\iota_{*}\mathcal{L}$ is a $Y-$flat family of stable sheaves on $S$ of Mukai vector $u$, hence we have a modular morphism $f_{\mathcal{F}}$. We may assume that $f'=f_{\mathcal{F}}$. 

Let $Y\subset\overline{Y}$ be a smooth algebraic compactification so that $f'$ admits an extension $\overline{f'}:\overline{Y}\rightarrow M'$. As $\overline{f}$ is a morphism between projective manifolds, and as $\lambda'(\alpha)\in H^{2,0}(M')$ is not zero, it follows that $\overline{f}^{*}(\alpha)\neq 0$. As the restriction of classes of holomorphic two forms from $\overline{Y}$ to $Y$ is always injective, we have $\overline{f'}^{*}(\lambda'(\alpha))_{|Y}\neq 0$, so $$(f')^{*}(\lambda'(\alpha))=\overline{f'}^{*}(\lambda'(\alpha))_{|Y}\neq 0\in H^2(Y,\mathbb{C}).$$

We now claim that
\begin{equation}
\label{eq:primopull}
(f')^{*}(\lambda'(\alpha))=-p_{Y*}(p_{S}^{*}(\sigma)\cdot\iota_{*}(c_1(\mathcal{L})))\in H^2(Y,\mathbb{C}).
\end{equation}
Indeed, since $f'$ is the modular map associated with $\mathcal{F}=\iota_{*}(\mathcal{L})$, by item (3) Proposition \ref{prop:descent} and using Grothendieck-Riemann-Roch we obtain $$(f')^{*}(\lambda'(\alpha))=[p_{Y*}(p_{S}^{*}(\alpha^{\vee}\cdot\sqrt{td(S)})\cdot ch(\iota_{*}(\mathcal{L}))]_{2}=$$ $$=[p_{Y*}(p_{S}^{*}(\alpha^{\vee}\cdot\sqrt{td(S)})\cdot \iota_{*}(ch(\mathcal{L})td(T_{\iota}))]_{2}=$$ $$=-p_{Y*}(p_{S}^{*}(\sigma)\cdot \iota_{*}(c_1(\mathcal{L})+td_1(T_{\iota}))=$$ $$=-p_{Y*}(p_{S}^{*}(\sigma)\cdot\iota_{*}(td_1(T_{\iota})))-p_{Y*}(p_{S}^{*}(\sigma)\cdot\iota_{*}(c_1(\mathcal{L})))$$where $T_{\iota}$ is the relative tangent bundle of $\iota$, i.e. $T_{\iota}=T_{\mathcal{C}}-\iota^*(T_{S\times Y})$ in the Grothendieck group of $\mathcal{C}$, and $td_{1}(T_{\iota})$ is the degree one part of its Todd class.

Since $S$ has trivial canonical bundle and every fiber of $g:\mathcal{C}\longrightarrow Y$ belongs to the linear system $|mH|$, shrinking $Y$ if necessary, we may assume that the relative canonical bundle of $g$ is represented by $\iota^{*}(p_{S}^{*}(mH))$. As a consequence, $td_1(T_{\iota})=-\frac{m}{2}\iota^{*}(p_{S}^{*}(H))$ and $$p_{Y*}(p_{S}^{*}(\sigma)\cdot\iota_{*}(td_1(T_{\iota})))=-\frac{m}{2}p_{Y*}(p_{S}^{*}(\sigma)\cdot\iota_{*}(\iota^{*}(p_{S}^{*}(H))))=$$ $$=-\frac{m}{2}p_{Y*}(p_{S}^{*}(\sigma)\cdot p_{S}^{*}(H)\cdot\iota_{*}([\mathcal{C}]))=0$$since $H$ is orthogonal to $\sigma$. This completes the proof of equation (\ref{eq:primopull}).

We now consider the flat family of sheaves on $S$ with Mukai vector $v$ given by $\iota_{*}(\mathcal{L}^{\otimes km^2})$ and the associated modular map $f:Y\rightarrow M$. The proof we just presented for equation (\ref{eq:primopull}) may be used to prove that
\begin{equation}
\label{eq:secondopull}
f^{*}(\lambda(\alpha))=p_{Y*}(p_{S}^{*}(\sigma)\cdot \iota_{*}(c_1(\mathcal{L}^{\otimes km^2}))).
\end{equation}
Equations (\ref{eq:primopull}) and (\ref{eq:secondopull}) give then that $f^{*}(\lambda(\alpha))=km^{2}(f')^{*}(\lambda'(\alpha))$. Now, notice that by construction we have $f=\psi\circ f'$, so $$(f')^{*}(\psi^{*}(\lambda(\alpha)))=f^{*}(\lambda(\sigma))=km^{2}(f')^{*}(\lambda'(\alpha)).$$

Since $\psi^{*}(\lambda(\alpha))=p\lambda'(\alpha)$ for some $p\neq 0$, the previous equation gives $p(f')^{*}(\lambda'(\alpha))=km^{2}(f')^{*}(\lambda'(\alpha))$, and since $(f')^{*}(\lambda'(\alpha))\neq 0$ we get $p=km^{2}$, concluding the proof.\endproof

We are now ready to prove Proposition \ref{prop:cucv1}, which implies Theorem \ref{thm:fujiki}.
\par\bigskip
{\em Proof of Proposition \ref{prop:cucv1}}. We use the definitions and notation introduced before in this section, and we denote by $C$ (resp. $C'$) and $b$ (resp. $b'$) the Fujiki constant and the Beauville-Namikawa form of $M$ (resp. $M'$). To prove Proposition \ref{prop:cucv1} we just need to prove that $C=C'$.

We let $2n$ be the dimension of $M'$ and $M$, i.e. $n=m^{2}k+1$ if $S$ is a K3 surface and $n=m^{2}k-1$ if $S$ is Abelian, and we choose a nontrivial class $\sigma\in H^{2,0}(S)$. 

Let $\alpha:=(0,\sigma+\overline{\sigma},0)\in u^{\perp}\cap v^{\perp}$. As a consequence $\lambda(\alpha)\in H^{2}(M,\mathbb{C})$ and $\lambda'(\alpha)\in H^{2}(M',\mathbb{C})$: by definition of the Fujiki constant and by Proposition \ref{prop:hodgeisoproof} we have $$\int_{M'}\lambda'(\alpha)^{2n}=C'b'(\lambda'(\alpha))^{n}=C'(\alpha,\alpha)^{n},$$and $$\int_{M}\lambda(\alpha)^{2n}=Cb(\lambda(\alpha))^{n}=C(\alpha,\alpha)^{n}.$$

Notice that $(\alpha,\alpha)\neq 0$, hence in order to show that $C=C'$ it is enough to show that 
\begin{equation}
\label{eq:u0}
\int_{M'}\lambda'(\alpha)^{2n}=\int_{M}\lambda(\alpha)^{2n}.
\end{equation}
By Proposition \ref{prop:psi} the degree of the generically finite, dominant rational morphism $\psi:M'\dashrightarrow M$ is $(km^2)^{2n}$ and, since 
$\lambda(\alpha)\in H^{2,0}(M)\oplus H^{2,0}(M)$,  we have 
\begin{equation}
\label{eq:u1}
\int_{M'}(\psi^{*}\lambda(\alpha))^{2n}=(km^2)^{2n}\int_{M}\lambda(\alpha)^{2n}.
\end{equation}

But by Proposition \ref{prop:relation} we have $\psi^{*}\lambda(\alpha)=(km^2)\lambda'(\alpha)$, hence we get
\begin{equation}
\label{eq:u2} 
\int_{M'}(\psi^{*}\lambda(\alpha))^{2n}=(km^2)^{2n}\int_{M'}\lambda'(\alpha)^{2n}.
\end{equation}
Comparing equations (\ref{eq:u1}) and (\ref{eq:u2}) we obtain equation (\ref{eq:u0}), and the result follows.\endproof

\section{Appendix: codimension of the complement of $R^{s}_{v}$ in $R^{ss}_{v}$}

Let $(S,v,H)$ be an $(m,k)-$triple where $m\geq 2$, and write $v=mw$. In this Appendix we aim to prove a technical result about the codimension of $R^{ss}_{v}\setminus R^{s}_{v}$ in $R^{ss}_{v}$ that is used in the proof of Proposition \ref{prop:extmr}.

As a first step we bound the dimension of the set of strictly semistable sheaves.

\begin{lem}
\label{lem:uppdimynp}
Let $(S,v,H)$ be an $(m,k)-$triple with $m\neq 1$ and $(m,k)\neq(2,1)$. There exists an algebraic variety $Y_{v}$ and a $Y_{v}$-flat coherent sheaf $\mathcal{E}_{v}$ on $S\times Y_{v}$ for which the two following properties hold: 
\begin{enumerate}
 \item for every strictly $H-$semistable sheaf $E$ on $S$ with Mukai vector $v$ there is $y\in Y_{v}$ such that the restriction $E_y$ of $\mathcal{E}_{v}$ to $S\times\{y\}$ is isomorphic to $E$,
\item we have $\dim(Y_{v})\leq\dim(M_{v})-3$.
\end{enumerate}
\end{lem}

\proof Let us first consider the case $m=2$, and write $v=2w$. We show that there exists an algebraic variety $Y_{v}$ satisfying (1) and (2) if $k\geq 2$, and satisfying (1) and $\dim(Y_{v})\leq\dim(M_{v})-1$ if $k=1$.

A strictly $H-$semistable sheaf $E$ on $S$ with Mukai vector $2w$ fits into an exact sequence of the form
\begin{equation} 
\label{eq:esattaapp} 
0\longrightarrow K\longrightarrow E\longrightarrow Q\longrightarrow 0.
\end{equation}
Notice that: trivial extensions of this form are parametrized by $M_w\times M_w$; isomorphism classes of non-trivial extensions with $K\nsimeq E$ are parametrized by a $\mathbb{P}^{w^2-1}$-bundle $P_1$ over the complement of the diagonal in $M_w\times M_w$; isomorphism classes of non-trivial extensions with $K\simeq E$ are parametrized by a $\mathbb{P}^{w^2+1}$-bundle $P_2$ over $M_w$.

We notice that as $\dim(M_{w})=w^{2}+2$ and $\dim(M_{2w})=4w^{2}+2$, we get $$\dim(M_{w}\times M_{w})=2w^2+4\leq\left\{\begin{array}{ll}\dim(M_{2w})-4, & {\rm if}\,\,k\geq 2\\ \dim(M_{2w})-2, & {\rm if}\,\,k=1\end{array}\right.$$ $$\dim(P_1)=3w^2+3\leq\left\{\begin{array}{ll} \dim(M_{2w})-3, & {\rm if}\,\,k\geq 2\\ \dim(M_{2w})-1, & {\rm if}\,\,k=1\end{array}\right.$$ and 
$$\dim(P_2)=2w^2+3\leq\dim(M_{2w})-3.$$

Notice that a tautological family exists on an \'etale cover of $M_w$, so a tautological family of trivial extensions exists on an \'etale cover of $M_w\times M_w$. Moreover, by Corollary 4.4. of \cite{L} and the considerations following it, there exist tautological families of non-trivial extensions on \'etale covers by affine varieties of $P_1$ and $P_2$.

As a consequence, there is an \'etale cover $Y_{v}$ of $(M_w\times M_w)\cup P_1 \cup P_2$ over which there is a tautological family of extensions. By the previous considerations about the dimension of $M_{w}\times M_{w}$, $P_{1}$ and $P_{2}$ it moreover follows that if $k\geq 2$ we have $\dim(Y_{v})\leq\dim(M_{v})-3$, while if $k=1$ we have $\dim(Y_{v})\leq\dim(M_{v})-1$. This concludes the proof of the statement in the case $m=2$.

It remains to deal with the case $m>2$ and, by induction on $m$, it is sufficient to show that, if for every $2\le i\le m-1$ there exist $Y_{iw}$ and $\mathcal{E}_{iw}$ satisfying (1) and such that $\dim(Y_{iw})\le\dim(M_{iw})$, then there exist $Y_{mw}$ and $\mathcal{E}_{mw}$ satisfying (1) such that $$\dim(Y_{mw})\le\dim(M_{mw})-3.$$

As before, a strictly $H-$semistable sheaf $E$ on $S$ with Mukai vector $mw$ always fits in the exact sequence (\ref{eq:esattaapp}) where $Q$ is an $H-$stable sheaf with Mukai vector $iw$ for $1\le i\le m-1$ and $K$ is an $H-$semistable sheaf with Mukai vector $(m-i)w$. Notice that either $K$ is an $H-$stable sheaf represented by a point in $M^{s}_{(m-i)w}$ or it is represented by a point of $Y_{(m-i)w}$.

We distinguish three cases:
\begin{enumerate}
 \item the extension is trivial,
 \item the extension is non-trivial and $Q$ is not a summand of the $H-$polystable sheaf in the same S-equivalence class of $K$, 
 \item the extension is non-trivial and $Q$ is a summand of the $H-$polystable sheaf in the same S-equivalence class of $K$.
\end{enumerate}

A parameter space for the extensions falling in the first case is $$V_{1,i}:=(M_{(m-i)w}^{s}\cup Y_{(m-i)w})\times M_{iw}^{s}.$$We notice that by induction we have $\dim(Y_{(m-i)w})\leq\dim(M^{s}_{(m-i)w})$, hence $$\dim(V_{1,i})=2k(m-i)^2+2+2ki^2+2=2km^{2}+4-4ki(m-i)=$$ $$=\dim(M_{mw})+2-4ki(m-i)\leq\dim(M_{mw})-6$$where we use that $m\ge 3$ implies $i(m-i)\ge2$.

In the second case, by hypothesis $\dim(ext^1(Q,K))=(m-i)iw^2$ and a parameter space for the extensions is $\mathbb{P}^{(m-i)iw^2-1}$ bundle $V_{2,i}$ over $V_{1,i}$: as a consequence we obtain $$\dim(V_{2,i})=2k(m-i)^2+2+2ki^2+2+2k(m-i)i-1=2km^{2}-2ki(m-i)+3=$$ $$=\dim(M_{mw})+1-2ki(m-i)\leq\dim(M_{mw})-3$$where again we use that since $m\geq 3$ then $i(m-i)\geq 2$.

In the remaining case we have $$\dim(Ext^1(Q,K))=(m-i)iw^2+\dim(Hom(Q,K))+\dim(Hom(K,Q)),$$and $\dim(Hom(Q,K))$ and $\dim(Hom(K,Q))$ are bounded above by the number of simple direct summands of the $H-$polystable sheaf in the S-equivalence class of $K$: as $K\in M_{(m-i)w}$, it follows that $$\dim(Ext^{1}(Q,K))\leq 2ki(m-i)+2(m-i).$$
 
The set of the extensions in the third case are then parametrized by a reducible variety $V_{3,i}$ whose components are projective bundles with fibers of dimension at most $2ki(m-i)+2(m-i)-1$ over suitable components of $Y_{(m-i)w}$ or $M_{(m-i)w}^s$. Since by induction we have $\dim(Y_{(m-i)w})\leq\dim(M_{(m-i)w})$, it follows that $$\dim(V_{3,i})\leq 2k(m-i)^2+2+2ki(m-i)i+2(m-i)-1=$$
$$=2km^2-2kmi+2ki^2(1-i)+2m-2i-1=$$ $$=2km^{2}+2-(m-i)(2ki-2)-2ki^{3}-3=$$ $$=\dim(M_{mw})-(m-i)(2ki-2)-2ki^{3}-3\le \dim(M_{mw})-3.$$

Finally, arguing as for $m=2$, for every $1\le i\le m-1$ and $1\le j\le3$ the variety $V_{j,i}$ admits an \'etale cover by affine varieties equipped with a tautological family of extensions. 

By letting $Y_{v}$ be the disjoint union of these affine varieties, we see that $Y_{v}$ is an algebraic variety equipped with the desired family $\mathcal{E}_{v}$, and such that $\dim(Y_{v})\leq\dim(M_{v})-3$, concluding the proof.\endproof

Using Lemma \ref{lem:uppdimynp} we are now able to prove the lower bound for the codimension of $R^{ss}_{v}\setminus R^{s}_{v}$ in $R^{ss}_{v}$ which is used in section 3.3.

\begin{prop}
\label{lem:codimrsv}
Let $(S,v,H)$ be an $(m,k)-$triple where $m\geq 2$ and $(m,k)\neq(2,1)$. Then $R^{ss}_{v}$ is a locally factorial, quasi-projective variety and $R^{s}_{v}$ is a smooth open subset of $R^{ss}_{v}$ whose complement $R^{ss}_{v}\setminus R^{s}_{v}$ has codimension at least 3 in $R^{ss}_{v}$.
\end{prop}

\proof By part 1 of Proposition 3.10 of \cite{KLS} we know that $R^{ss}_{v}$ is locally of complete intersection and its singular locus has codimension at least 3. By Corollary 3.14 in Exp. XI of \cite{SGA2} it follows that $R^{ss}_{v}$ is locally factorial.

As $R^{s}_{v}$ parameterizes stable quotients, it is a smooth open subset of $R^{ss}_{v}$, and since $q^{s}_{v}:R^{s}_{v}\longrightarrow M^{s}_{v}$ is a $PGL(N_{v})-$bundle, we see that $R^{s}_{v}$ has dimension equal to $\dim(M^{s}_{v})+\dim(PGL(N_{v}))$. 

As, by Theorem 4.1 of \cite{LS} $M_{v}$ is irreducible and it is the quotient of $R^{ss}_{v}$ under the action of the connected group $PGL(N_{v})$, we see that $R^{ss}_{v}$ is irreducible: it follows that $R^{s}_{v}$ is a dense open subset of $R^{ss}_{v}$, so $$\dim(R^{ss}_{v})=\dim(R^{s}_{v})=\dim(M_{v})+\dim(PGL(N_{v})).$$

Recall from Section 4.3 of \cite{HL} that $R^{ss}_{v}$ is an open subset of a Quot-scheme parametrizing quotients of a vector bundle $\mathcal{H}_v$ on $S$ such that $\dim(Hom(\mathcal{H}_v,F))$ is constant on the set of $H-$semistable sheaves with Mukai vector $v$ on $S$. 

This implies the following: consider a family $\mathcal{E}$ of $H-$semistable sheaves on $S$ with Mukai vector $v$ parametrized by a variety $Y$, and choose a point $y\in Y$. If $E_{y}$ denotes the restriction of $\mathcal{E}$ to $S\times\{y\}$, let $\phi_{y}:\mathcal{H}_v\longrightarrow E_y$ be a surjective morphism: then there is a Zariski open neighborhood $U$ of $y$ in $Y$ such that $\phi_{y}$ extends to a surjective morphism from the trivial family over $U$ with fiber $\mathcal{H}_v$ to the restriction of $\mathcal{E}$ to $S\times U$.

As a consequence, consider an algebraic variety $Y_{v}$ and the family $\mathcal{E}_{v}$ verifying points (1) and (2) of Lemma \ref{lem:uppdimynp}. Up to replace $Y_{v}$ by a suitable affine covering we may assume that the sheaf $\mathcal{E}_{v}$ is obtained 
from the universal quotient on $R^{ss}_{v}$ by base change under a modular morphism $f:Y_{v}\rightarrow R^{ss}_{v}$.

By property (1) of Lemma \ref{lem:uppdimynp} every strictly $H-$semistable sheaf with Mukai vector $v$ is represented by some point in the image of $f$. Since points of $R^{ss}_{v}\setminus R^{s}_{v}$ represent strictly $H-$semistable sheaves, and since two quotients $[\mathcal{H}_{v}\longrightarrow F_{1}]$ and $[\mathcal{H}_{v}\longrightarrow F_{2}]$ lie in the same $PGL(N_{v})-$orbit if and only if $F_{1}$ and $F_{2}$ are isomorphic (as abstract coherent sheaves), it follows that $$\dim(R^{ss}_{v}\setminus R^{s}_{v})\le \dim(Y_{v})+\dim(PGL(N_{v}))$$ and by (2) of Lemma \ref{lem:uppdimynp} we conclude that $$\dim(R^{ss}_{v}\setminus R^{s}_{v})\le \dim(M_v)+\dim(PGL(N_{v}))-3=\dim(R^{ss}_{v})-3,$$concluding the proof.\endproof

We conclude with a statement about Abelian surfaces which is similar to Proposition \ref{lem:codimrsv}, where we replace $R^{ss}_{v}$ by the inverse image $R^{ss,0}_{v}$ of $K_v$ in $R^{ss}_{v}$, and $R^{s}_{v}$ by the inverse image $R^{s,0}_{v,0}$ of $K_v^{s}$ in $R^{ss}_{v}$.

\begin{cor}
\label{cor:codimrsv}
Let $(S,v,H)$ be an $(m,k)-$triple where $S$ is Abelian, $m\geq 2$ and $(m,k)\neq(2,1)$. Then $R^{ss,0}_{v}$ is locally factorial, $R^{s,0}_{v}$ is a smooth open subset of $R^{ss,0}_{v}$ and $R^{ss,0}_{v}\setminus R^{s,0}_{v}$ has codimension at least 3 in $R^{ss,0}_{v}$.
\end{cor}

\proof If $F$ is an $H-$polystable sheaf with Mukai vector $v$ on $S$, $U$ is an open neighborhood of the quotient representing $F$ in $R^{ss}_{v}$ and $D(F)$ is a universal deformation space for $F$, by Luna's \'Etale Slice Theorem the codimension in $U$ of the locus parametrizing strictly $H-$semistable quotients equals the codimension in $D(F)$ of the locus defined by the same property.
 
Since $S\times\widehat{S}$ acts transitively on the fibers of $a_v:M_v\longrightarrow S\times\widehat{S}$ (by translation and tensorization), the codimension of $R^{ss,0}_{v}\setminus R^{s,0}_{v}$ in $R^{ss,0}_{v}$ equals the codimension of $R^{ss}_{v}\setminus R^{s}_{v}$ in $R^{ss}_{v}$, which by Proposition \ref{lem:codimrsv} it is at least $3$. 

Since $R^{s,0}_{v}$ is smooth, it follows that the singular locus of $R^{ss,0}_{v}$ has codimension at least $3$ and, since 
$R^{ss,0}_{v}$ is a fiber of the expected dimension of the surjective morphism $a_{v}\circ q_{v}:R^{ss}_{v}\rightarrow S\times\widehat{S}$, the subscheme $R^{ss,0}_{v}$ of the local complete intersection $R^{ss}_{v}$ is a local complete intersection too. As in the proof of Proposition \ref{lem:codimrsv}, by Corollary 3.14 in Exp. XI of \cite{SGA2} we deduce that $R^{ss,0}_{v}$ is locally factorial, completing the proof.\endproof

\end{document}